\documentclass[10pt]{amsart}
\usepackage{amscd, amssymb}
\usepackage{array}
\def\today{\number\day\space\ifcase\month\or January\or February\or
March\or April\or May\or June\or July\or August\or September\or
October\or November\or December\fi \space\number\year}

\newcommand{\ind}{\operatorname{ind}}
\newcommand{\edge}{\operatorname{edge}}
\newcommand{\comsquare}[8]
{\begin{CD}
#1 @>#2>> #3 \\
@V{#4}VV @VV{#5}V \\
#6 @>>#7> #8
\end{CD}
}

\newcommand{\af}{\alpha}
\newcommand{\bt}{\beta}

\newcommand{\dt}{\delta}
\newcommand{\ep}{\varepsilon}
\newcommand{\zt}{\zeta}

\newcommand{\te}{\theta}
\newcommand{\ld}{\lambda}
\newcommand{\sm}{\sigma}

\newcommand{\ph}{\varphi}
\newcommand{\ps}{\psi}
\newcommand{\rh}{\rho}
\newcommand{\om}{\omega}
\newcommand{\ta}{\tau}

\newcommand{\Te}{\Theta}

\newcommand{\Om}{\Omega}

\newcommand{\ov}{\overline}

\newcommand{\I}{\infty}

\newcommand{\ts}[1]{{\textstyle{#1}}}

\newcommand{\andeqn}{\quad {\text{and}} \quad}

\newcommand{\E}{{\mathcal{E}}}

\newcommand{\cH}{{\mathcal{H}}}
\newcommand{\K}{{\mathcal{K}}}
\newcommand{\PU}{{\mathcal{P U}}}
\newcommand{\cS}{{\mathcal{S}}}
\newcommand{\U}{{\mathcal{U}}}

\newcommand{\Q}{{\mathbb{Q}}}
\newcommand{\Z}{{\mathbb{Z}}}
\newcommand{\R}{{\mathbb{R}}}
\newcommand{\C}{{\mathbb{C}}}
\newcommand{\N}{{\mathbb{N}}}

\newcommand{\CC}{{\mathbb{C}}}

\newcommand{\QQ}{{\mathbb{Q}}}
\newcommand{\RR}{{\mathbb{R}}}
\newcommand{\TT}{{\mathbb{T}}}
\newcommand{\ZZ}{{\mathbb{Z}}}

\newcommand{\card}{{\operatorname{card}}}
\newcommand{\dist}{{\operatorname{dist}}}
\newcommand{\Aut}{\operatorname{Aut}}

\newcommand{\ab}{\operatorname{ab}}
\newcommand{\pr}{\operatorname{pr}}
\newcommand{\Id}{\operatorname{Id}}
\newcommand{\SL}{\operatorname{SL}}

\newcommand{\pt}{\operatorname{pt}}
\newcommand{\res}{\operatorname{res}}

\newcommand{\Top}{\operatorname{top}}

\newcommand{\Ad}{\operatorname{Ad}}
\newcommand{\id}{\operatorname{id}}

\newcommand{\ev}{\operatorname{ev}}

\newcommand{\EG}{{\underline{E}} G}

\newcommand{\lk}{\langle}
\newcommand{\rk}{\rangle}

\newcommand{\ca}{C*-algebra}
\newcommand{\suca}{simple unital \ca}
\newcommand{\ct}{continuous}
\newcommand{\pj}{projection}

\newcommand{\hm}{homomorphism}
\newcommand{\wolog}{without loss of generality}
\newcommand{\Wolog}{Without loss of generality}

\newcommand{\ifo}{if and only if}

\newcommand{\mops}{mutually orthogonal \pj s}

\newcommand{\hsa}{hereditary subalgebra}

\newcommand{\mvnt}{Murray-von Neumann equivalent}

\newcommand{\tRp}{tracial Rokhlin property}

\newcommand{\fd}{finite dimensional}
\newcommand{\uct}{Universal Coefficient Theorem}

\newcommand{\tst}{tracial state}  %

\theoremstyle{plain}
\newtheorem{thm}{Theorem}[section]
\newtheorem{cor}[thm]{Corollary}
\newtheorem{lem}[thm]{Lemma}
\newtheorem{prop}[thm]{Proposition}
\theoremstyle{definition}
\newtheorem{defn}[thm]{Definition}
\newtheorem{remark}[thm]{Remark}
\newtheorem{ex}[thm]{Example}
\newtheorem{ntn}[thm]{Notation}
\numberwithin{equation}{section}
\title[Crossed products of irrational rotation algebras by finite
groups]{The structure of crossed products of irrational rotation
algebras by finite subgroups
    of $\SL_2 (\ZZ)$}

\date{25~Sept.\  2006}

\author[Echterhoff]{Siegfried Echterhoff}
\address{Westf\"{a}lische Wilhelms-Universit\"{a}t M\"{u}nster,
   Mathematisches Institut, Einsteinstr.~62,
   D-48149 M\"{u}nster, Germany}
\email{echters@math.uni-muenster.de}
\author[L\"{u}ck]{Wolfgang L\"{u}ck}
\address{Westf\"{a}lische Wilhelms-Universit\"{a}t M\"{u}nster,
   Mathematisches Institut, Einsteinstr.~62,
   D-48149 M\"{u}nster, Germany}
\email{lueck@math.uni-muenster.de}
\author[Phillips]{N.~Christopher Phillips}
\address{Department of Mathematics, University  of Oregon,
        Eugene OR 97403-1222, USA.}
\email[]{ncp@darkwing.uoregon.edu}
\author[Walters]{Samuel Walters}
\address{Department of Mathematics and Computer Science,
The University of Northern British Columbia,
Prince George, B.C. V2N 4Z9, Canada}
\email{walters@hilbert.unbc.ca or walters@unbc.ca }
\thanks{Research of the first and second authors partially
supported by the Deutsche Forschungsgemeinschaft (SFB~478).
Research of the third author partially supported by
NSF grants DMS~0070776 and DMS~0302401.}

\begin{document}

\begin{abstract}
Let $F \subseteq \SL_2 (\ZZ)$ be a finite subgroup
(necessarily isomorphic to one of
$\ZZ_2,$ $\ZZ_3,$ $\ZZ_4,$ or $\ZZ_6$),
and let $F$ act on the irrational rotational algebra $A_{\te}$ via the
restriction of the canonical action of $\SL_2 (\ZZ).$
Then the crossed product
$A_{\theta} \rtimes_{\alpha} F$ and the
fixed point algebra $A_{\te}^F$ are AF~algebras.
The same is true
for the crossed product and fixed point algebra of the
flip action of $\ZZ_2$ on any simple $d$-dimensional
noncommutative torus $A_{\Theta}.$
Along the way, we prove a
number of general results which should have
useful applications in other situations.
\end{abstract}
\maketitle

\setcounter{section}{-1}

\section{Introduction and statement of the main results}

For $\te \in \R$ let $A_{\te}$ be the rotation algebra, which is the
universal \ca\  generated by unitaries $u$ and $v$
satisfying $v u = \exp (2 \pi i \te) u v.$
The group $\SL_2 (\ZZ)$ acts on $A_{\te}$ by sending the
matrix
\[
n = \left( \begin{array}{cc} n_{1, 1} & n_{1, 2} \\
            n_{2, 1} & n_{2, 2} \end{array} \right)
\]
to the automorphism determined by
\begin{equation}\label{eq-action}
\af_n (u)
= \exp (\pi i n_{1, 1} n_{2, 1} \te) u^{n_{1, 1}} v^{n_{2, 1}}\;
{\text{and}} \;
\af_n (v)
= \exp (\pi i n_{1, 2} n_{2, 2} \te) u^{n_{1, 2}} v^{n_{2, 2}}.
\end{equation}
In this paper we are concerned with the structure
of the crossed products $A_{\te} \rtimes_{\alpha} F$ for
the finite subgroups $F \subseteq \SL_2 (\ZZ),$
which are, up to conjugacy, the groups
$F = \ZZ_2, \ZZ_3, \ZZ_4, \ZZ_6$ (where $\ZZ_m$
stands for the cyclic group of order $m$)
generated by the matrices
\begin{equation}\label{eq-gen}
\begin{split}
&\left(\begin{smallmatrix} -1 & 0 \\ 0 & - 1 \end{smallmatrix} \right)
\;\; {\text{(for $\ZZ_2$)}}, \quad \quad\quad \quad\quad
\left(\begin{smallmatrix} -1 & -1 \\1 & 0 \end{smallmatrix} \right)
\;\; {\text{(for $\ZZ_3$)}},
\\
&\left(\begin{smallmatrix} 0 & -1 \\ 1 & 0 \end{smallmatrix} \right)
\;\; \;\;{\text{(for $\ZZ_4$)}}, \quad \quad {\text{and}} \quad \quad
\left(\begin{smallmatrix} 0 & -1 \\ 1 & 1 \end{smallmatrix} \right)
\;\; {\text{(for $\ZZ_6$)}}.
\end{split}
\end{equation}
We refer to Proposition~21 of~\cite{P_FW4} for a
proof that one obtains essentially all interesting actions
of the cyclic groups $\ZZ_2, \ZZ_3, \ZZ_4, \ZZ_6$
on $A_{\te}$ by restricting the action $\af$ of
$\SL_2 (\ZZ)$ to these subgroups.
The main results
of this paper culminate in a proof of the following theorem:

\begin{thm}[Theorems \ref{thm-generators}, \ref{T:CrPrdIsAF},
           and~\ref{T:IsomOfCrPrd}]\label{thm-culmain}
Let $F$ be any of the finite subgroups
$\ZZ_2, \ZZ_3, \ZZ_4, \ZZ_6 \subseteq \SL_2 (\ZZ)$ with generators
given as above and let $\theta \in \RR \smallsetminus \QQ.$
Then the crossed product $A_{\te} \rtimes_{\alpha} F$ is
an AF~algebra. For all $\theta \in \RR$ we have
\[
K_0 (A_{\theta} \rtimes_{\af} \ZZ_2) \cong \ZZ^6, \quad
K_0 (A_{\theta} \rtimes_{\af} \ZZ_3) \cong \ZZ^8,
\]
\[
K_0 (A_{\theta} \rtimes_{\af} \ZZ_4) \cong \ZZ^9,
\quad {\text{and}} \quad
K_0 (A_{\theta} \rtimes_{\af} \ZZ_6) \cong \ZZ^{10}.
\]
For $F = \ZZ_k$ for $k = 2, 3, 4, 6,$
the image of $K_0 (A_{\theta} \rtimes_{\af} F)$
under the canonical \tst\  on $A_{\theta} \rtimes_{\alpha} F$
(which is unique) is equal to
$\frac{1}{k} (\ZZ + \theta \ZZ).$
As a consequence, $A_{\te} \rtimes_{\alpha} \ZZ_k$
is isomorphic to $A_{\te'} \rtimes_{\alpha} \ZZ_l$ if and
only if $k = l$ and  $\theta' = \pm \theta$ mod $\ZZ.$
\end{thm}

If $\theta\in \RR \smallsetminus \QQ$ then the fixed point
algebra $A_{\te}^F$ is Morita equivalent
to $A_{\te} \rtimes_{\alpha} F.$
(This follows from the Proposition in~\cite{Rs}.)
Thus, as a direct corollary
of Theorem~\ref{thm-culmain} we get:

\begin{cor}[Corollary~\ref{C:FixIsAF}]\label{cor-fix}
Let $\alpha:F\to \Aut(A_{\te})$ be
as above with  $F = \ZZ_2, \ZZ_3, \ZZ_4, \ZZ_6.$
Then the fixed point algebras
$A_{\te}^F$ are~AF for all $\theta \in \RR \smallsetminus \QQ.$
\end{cor}

The proof of Theorem~\ref{thm-culmain} has three independent steps:
\begin{enumerate}
\item\label{Step1}
Computation of the $K$-theory of the crossed product.
\item\label{Step2}
Proof that the crossed product satisfies the
Universal Coefficient Theorem.
\item\label{Step3}
Proof that the action of the group has the \tRp\  of~\cite{PhtRp1a}.
\end{enumerate}
Given these steps, one uses Theorem~2.6 of~\cite{PhtRp1a}
to prove that the crossed product has tracial rank zero
in the sense of~\cite{LnTAF, LnTTR},
and then Huaxin Lin's
classification theorem (Theorem~5.2 of~\cite{Ln15})
to conclude that the crossed product is~AF.

Our results have been motivated by previous
studies of the algebras $A_{\te} \rtimes_{\af} F$ and
$A_{\te}^F$ by several authors.
See, for example, \cite{BEEK, P_BK, P_FW1, FW1, P_FW2, P_FW6, P_FW5,
P_FW3, P_FW4, Kj, Pl, W1, Wl1, W3}.
The most was
known about the crossed product by the flip:
its (unordered) $K$-theory has been computed in~\cite{Kj},
and the crossed product has been proved to
be an AF algebra in~\cite{P_BK}.
The next best understood case is the case
of $\ZZ_4,$ which has been intensively
studied in a series of papers culminating in~\cite{W3}.
It is proved there that for most irrational $\te$
(that is, a dense $G_{\dt}$-set of them)
the crossed product  $A_{\te} \rtimes_{\alpha} \ZZ_4$
has tracial rank zero in the sense of~\cite{LnTAF},
and that for most of those values of $\te,$ it is an AF~algebra.
Step~(\ref{Step1}) is done for most $\te$ in~\cite{Wl1},
Step~(\ref{Step2}) is done for all $\te$ in~\cite{W3},
and, in place of Step~(\ref{Step3}),
a direct proof is given in~\cite{W3}
that the crossed product has tracial rank zero for most $\te.$
That $A_{\te} \rtimes_{\alpha} \ZZ_4$ is
an AF~algebra for most $\te$ then follows from Lin's
classification theorem.

Using completely different methods, we
prove in this paper several general results which
will imply, as special cases, the three steps
described above
for {\emph{all}} $\theta \in \RR \smallsetminus \QQ$
and for {\emph{all}}
choices for $F = \ZZ_2, \ZZ_3, \ZZ_4, \ZZ_6.$
For the computation of the $K$-theory (Step~(\ref{Step1}) above),
we observe that the crossed products
$A_{\te} \rtimes_{\alpha} F$ can be realized as
(reduced) twisted group algebras
$C_r^* (\ZZ^2 \rtimes F, \, \widetilde{\om}_{\theta}),$
where $\ZZ^2 \rtimes F$ is the semidirect product
by the obvious action of $F$ on $\ZZ^2$
and $\widetilde{\om}_{\te}\in Z^2(\ZZ^2 \rtimes F, \, \TT)$
is a suitable circle valued $2$-cocycle on
$\ZZ^2 \rtimes F.$
Since $\widetilde{\om}_{\theta}$
is homotopic to the
trivial cocycle in a sense explained below,
the following result will imply that
the $K$-groups of $A_{\te} \rtimes_{\alpha} F$
are isomorphic to the
$K$-groups of the untwisted group
algebra $C_r^*(\ZZ^2 \rtimes F) \cong C (\TT^2) \rtimes F,$
which we compute (together with an explicit basis) by using
some general methods from~\cite{Davis-Lueck (2003)}.

\begin{thm}[Theorem~\ref{thm-twist2}]\label{thm-twist}
Suppose that $\om_0, \om_1\in Z^2(G, \TT)$ are
homotopic Borel $2$-cocycles
on the second countable locally compact group $G,$
that is, there exists a Borel $2$-cocycle
$\Om \in Z^2 (G, \, C ([0,1], \, \TT))$ such that
$\om_j = \Om (\, \cdot \, , \, \cdot \, ) (j)$ for $j = 0, 1.$
Suppose further
that $G$ satisfies the Baum-Connes
conjecture with coefficients.
(If $G$ is amenable or a-T-menable,
this is automatic by Theorem~1.1 of~\cite{HK}.)
Then
$K_* (C_r^* (G, \om_0)) \cong K_* (C_r^* ( G, \om_1)).$
\end{thm}

Indeed, using the homotopy $\Om \in Z^2 (G, C ([0,1], \, \TT))$
one can construct a \ca\  %
$C_r^* (G, \Om)$ which is a continuous field of \ca s
over $[0, 1]$ with fibers isomorphic
to the twisted group algebras $C_r^* (G, \om_{\theta}),$
with $\om_{\theta} = \Om (\, \cdot \, , \, \cdot \, ) (\theta)$ for
$\theta \in [0,1].$
We show, under the conditions on $G$ in Theorem~\ref{thm-twist}
(in fact, slightly more generally),
that the quotient maps
$q_{\theta} \colon C_r^* (G, \Om) \to C_r^* (G, \om_{\theta})$
induce isomorphisms on $K$-theory.
This special structure of the isomorphism will be used
for our description of an explicit basis
of $A_{\te} \rtimes_{\af} F$:  we extend the basis
of $K_0 (C_r^* (\ZZ^2 \rtimes F))$ to a basis
of $K_0 (C_r^* (\ZZ^2 \rtimes F, \, \Om))$ and then evaluate
at the other fibers.

The $K_0$-groups of the
crossed products of the irrational
rotation algebras by the actions of $\Z_2,$ $\Z_3,$ $\Z_4,$
and $\Z_6$ have been computed independently in~\cite{Pl},
using completely different methods.
These methods give no information about the $K_1$-groups,
and therefore do not suffice for the proof that the
crossed products are AF~algebras.
We should also mention that Georges Skandalis
pointed out to us an alternative way to
compute the $K$-theory of the crossed products
$A_{\te} \rtimes_{\alpha} F,$ also
using the Baum-Connes conjecture.
However, the method described above
seems best suited for computing
an explicit basis for $K_0 (A_{\te} \rtimes_{\alpha} F).$

Step~(\ref{Step2}) (the \uct)
will be obtained from the Baum-Connes conjecture
by slightly extending ideas first used in
Section~10.2 of~\cite{Tu} and Lemma~5.4 of~\cite{CEO2}.

Step~(\ref{Step3})
(the \tRp) is obtained via Theorem~\ref{OuterImpTRP},
according to which an action of a finite
group on a simple separable unital nuclear \ca\  with unique \tst\  %
has the \tRp\  \ifo\  the
induced action is outer on the type~II$_1$ factor
obtained from the trace using the Gelfand-Naimark-Segal construction.
Along the way,
we give a characterization of the \tRp\  in terms of trace norms
which does not require uniqueness of the \tst.

The methods developed here can also be applied
to the crossed products of all simple $d$-dimensional
noncommutative tori $A_{\Theta},$ for $d \geq 2,$
by the flip action of $\ZZ_2,$ which
sends the unitary generators $u_1, u_2, \ldots, u_d$ to
their adjoints $u_1^*, u_2^*, \ldots, u_d^*.$
Using  the general methods explained above together
the main result of~\cite{PhtRp2},
which shows that every simple higher
dimensional $d$-torus is an AT~algebra, we prove:

\begin{thm}[Corollary~\ref{KThOfNCTorus}, Theorem~\ref{T:HDFlipIsAF},
               and Corollary~\ref{C:FlipFixIsAF}]\label{them-Theta}
Let $A_{\Theta}$ be the noncommutative
$d$-torus corresponding to a real $d \times d$
skew symmetric matrix $\Theta.$
Let
$\alpha \colon \ZZ_2 \to \Aut (A_{\Theta})$ denote the flip action.
Then
\[
K_0 (A_{\Theta} \rtimes_{\alpha} \ZZ_2) \cong \ZZ^{3 \cdot 2^{d-1}}
\andeqn
K_1 (A_{\Theta} \rtimes_{\alpha} \ZZ_2) = \{ 0 \}.
\]
If $A_{\Theta}$ is simple, then
$A_{\Theta} \rtimes_{\alpha} \ZZ_2$ and the
fixed point algebra $A_{\Theta}^{\ZZ_2}$ are AF~algebras.
\end{thm}

This result generalizes Theorem~3.1 of~\cite{Bc}
and completely answers a question raised
in the introduction of~\cite{P_FW2}.

In Section~\ref{sec-twist}, we give the necessary background
on twisted group algebras
and the proof of Theorem~\ref{thm-twist}.
The realization
of the crossed products $A_{\te} \rtimes_{\alpha} F$
as twisted group algebras of the group
$\ZZ^2 \rtimes F$ is done in Section~\ref{sec-action}.
The $K$-theory computations
for $C_r^* (\ZZ^2 \rtimes F)$ are done in
Sections \ref{sec-untwist} and~\ref{sec-gen}.
Section~\ref{Sec:TRP} contains the proof that the relevant actions
have the \tRp.
In Section~\ref{Sec:CP},
we prove the \uct\  and put everything together.
We also prove Theorem~\ref{them-Theta}.

This paper contains the main result of Section~10
of the unpublished long preprint~\cite{PhW},
Theorem~\ref{them-Theta} here.
Although the three steps,
as described at the beginning of the introduction,
are the same, the proofs of all three of them differ
substantially from the proofs given in~\cite{PhW}.
This paper also supersedes Sections~8 and~9 of~\cite{PhW}.
Again, the proofs of Steps~(\ref{Step2}) and~(\ref{Step3})
differ substantially from those in~\cite{PhW}.
We improve on~\cite{PhW} by calculating the $K$-theory of the
crossed products $A_{\te} \rtimes_{\alpha} F.$
No $K$-theory calculations were done in~\cite{PhW},
so that the crossed products were merely shown to be
simple AH~algebras with real
rank zero for $F = \ZZ_3, \ZZ_4, \ZZ_6$ and all irrational $\te,$
and to be AF~algebras for $F = \ZZ_4$ and a dense $G_{\dt}$-set
of values of~$\te.$
(The case $F = \ZZ_2$ was already known,
and in any case it is covered by Theorem~\ref{them-Theta}.)

The third author would like to thank Masaki Izumi for
useful discussions concerning the \tRp.

\section{Twisted group algebras and actions on $\K$}\label{sec-twist}

We begin by recalling some basic facts about group \ca s
twisted by cocycles.
This is a special case of the theory of crossed product \ca s
twisted by cocycles of~\cite{PR}.
Also see~\cite{ZM} for the case in which the group is discrete.

Let $G$ be a second countable locally compact group,
with modular function $\Delta_G \colon G \to (0, \infty).$
Let $\om \colon G \times G \to \TT$
be a Borel $2$-cocycle on $G.$
(Recall the algebraic conditions:
$\om (s, t) \om (r, s t) = \om (r, s) \om (r s, t)$
and $\om (1, s) = \om (s, 1) = 1$ for $r, s, t \in G.$)
Then the twisted convolution algebra $L^1 (G, \om)$
is defined to be the vector space of all
integrable complex functions on $G$ with convolution and
involution given by
\[
(f *_{\om} g) (s)
= \int_G f (t) g (t^{-1} s) \om (t, \, t^{-1} s) \, d t \quad
{\text{and}} \quad
f^* (s) = \Delta_G (s^{-1}) {\overline{\om (s, s^{-1}) f (s^{-1})}}.
\]
An {\emph{$\om$-representation}} of $G$ on a Hilbert space
$\cH$ is a Borel map $V \colon G \to \U (\cH),$
the unitary group of $\cH,$
with respect to the strong operator topology on $\U (\cH),$
satisfying
\[
V (s) V (t) = \om (s, t) V (s t)
\]
far all $s, t \in G.$
The {\emph{regular $\om$-representation}} of $G$ is the
representation $L_{\om} \colon G \to \U (L^2 (G))$
given by the formula
\[
(L_{\om} (s) \xi) (t)
   = \om (s, \, s^{-1} t) \xi (s^{-1} t)
\]
for $\xi \in L^2 (G).$
Every $\om$-representation $V \colon G \to \U (\cH)$
determines a contractive $*$-homomorphism,
also denoted $V,$ from $L^1 (G, \om)$
to $B (\cH)$ via the formula
\[
V (f) = \int_G f (s) V (s) \, d s
\]
for $f \in L^1 (G, \om),$
and every nondegenerate representation of
$L^1 (G, \om)$ appears in this way.
The full twisted group algebra $C^*(G, \om)$
is defined to be the enveloping
\ca\  of $L^1 (G, \om)$
and the reduced twisted group algebra
$C_r^* (G, \om)$ is defined to be the image
of $C^* (G, \om)$ under the regular
$\om$-representation $L_{\om}.$
If $G$ is amenable, then both algebras coincide
(Theorem~3.11 of~\cite{PR}).
In this case we simply write $C^* (G, \om).$

Recall that two cocycles $\om, \om' \in Z^2 (G, \TT)$
are cohomologous
if there exists a Borel function $u \colon G \to \TT$
such that,
with $\partial u (s, t) = u (s) u (t) \overline{u (s t)},$
we have $\om' = \partial u \cdot \om.$
(See Chapter~1 of~\cite{mo1}; the action of $G$ on $\TT$ is trivial.)
It is then easy to verify that
the map $L^1 (G, \om') \to L^1 (G, \om),$
given by $f \mapsto u \cdot f,$
extends to isomorphisms of the
full and reduced twisted group \ca s.

In this section, let $\K$ denote the \ca\  of compact operators
on a separable infinite dimensional Hilbert space $\cH,$
let $\U$ be the group of unitary operators on $\cH,$
and let $\PU = \U / \TT 1$ denote the projective unitary group.
If $V \colon G \to \U$ is an $\om$-representation,
then $V$ determines an action
$\alpha \colon G \to \PU \cong \Aut (\K)$ by
$\alpha_s = \Ad (V (s))$ for $s \in G.$
(Continuity follows from Proposition~5(a) of~\cite{mo3}.)
With ${\overline{\om}} (s, t) = {\overline{\om (s, t)}},$
the full and reduced
crossed products $\K \rtimes_{\alpha} G$ and
$\K \rtimes_{\alpha, r} G$ are
then isomorphic to $C^* (G, {\overline{\om}}) \otimes \K$
and $C_r^* (G, {\overline{\om}}) \otimes \K.$
On the level of
$L^1$-algebras the isomorphisms are given by
the map $\Phi \colon L^1 (G, {\overline{\om}}) \odot \K \to L^1 (G, \K)$
determined by
\begin{equation}\label{eq-iso}
\Phi (f \otimes k) (s) = f (s) k V (s)^*
\end{equation}
for $f \in L^1 (G, {\overline{\om}}),$ $k \in \K,$ and $s \in G.$
(See for example Theorem~1.4.15 and Example~1.1.4 of~\cite{ech-mem},
but this is easily proved directly.)
Conversely, if $\alpha \colon G \to \Aut (\K)$
is any action of $G$ on $\K,$
then by choosing a Borel section $c \colon \PU \to \U$ we obtain a
Borel map $V^{\alpha} = c \circ \alpha \colon G \to \U$
such that
$\alpha_s = \Ad (V^{\alpha} (s))$ for all $s \in G.$
It is then easy to check that there exists a
Borel cocycle $\om_{\alpha} \in Z^2 (G, \TT)$
such that
$V^{\alpha} (s) V^{\alpha} (t) = \om_{\alpha} (s, t) V^{\alpha} (s t)$
for all $s, t \in G.$
Thus, any action $\af \colon G \to \Aut(\K)$
is given by an $\om_{\alpha}$-representation
for some cocycle $\om_{\alpha} \in Z^2 (G, \TT).$

If $A$ is any \ca, then two actions
$\alpha, \alpha' \colon G \to \Aut (A)$ are called
{\emph{exterior equivalent}} if there exists a strictly
continuous map $v \colon G \to \U M (A)$
(the unitary group in the multiplier algebra of $A$) such that
$\alpha'_s = \Ad (v_s) \circ \alpha_s$ and
$v_{s t} = v_s \alpha_s (v_t)$ for all $s, t \in G.$
(Compare with~8.11.3 of~\cite{Pd1}.)
It is easily seen that exterior equivalent actions have isomorphic
crossed products and isomorphic reduced crossed products,
with isomorphisms
given on $L^1 (G, A)$ by $f \mapsto f \cdot v.$

Following~\cite{EW1},
if $X$ is any locally compact space
we denote by $\E_G (X)$ the set of all
exterior equivalence classes of $C_0 (X)$-linear actions of $G$
on $C_0 (X, \K),$ that is, actions
$\alpha \colon G \to \Aut \big( C_0 (X, \K) \big)$
such that $\alpha_s (f \cdot g) = f \cdot \alpha_s (g)$ for
all $f \in C_0 (X)$ and $g \in C_0 (X, \K).$
Any such action is completely determined by
its evaluations $\alpha^x \colon G \to \Aut (\K)$ for $x \in X,$
defined by $\alpha^x_s (g (x)) = \big( \alpha_s (g) \big) (x)$
for $g \in C_0 (X, \K)$ and $x \in X.$
Thus, we should understand
a $C_0 (X)$-linear action as a continuous
family of actions on $\K.$
It is shown in~\cite{CKRW} (see Lemma~3.1 and Theorem~3.6)
and~\cite{EW1} that $\E_G (X)$ is
a group with multiplication
$[\alpha] \cdot [\beta] = [\alpha \otimes_X \beta],$ where
$(\alpha \otimes_X \beta)^x = \alpha^x \otimes \beta^x$
for $x \in X.$
The following result,
basically due to Mackey, is well known;
for example see~6.3 in~\cite{CKRW}:

\begin{prop}\label{prop-action}
Let $G$ be a second countable locally compact group.
Two actions $\alpha$ and $\alpha'$ of $G$ on $\K$
are exterior equivalent if and only if
$[\om_{\alpha}] = [\om_{\alpha'}] \in H^2 (G, \TT).$
Moreover, the mapping $[\alpha] \mapsto [\om_{\alpha}]$
is a bijective homomorphism from $\E_G (\pt)$ to $H^2 (G, \TT).$
\end{prop}

\begin{proof}
A calculation, which we omit, shows that if $\af$ and $\af'$
are two actions with corresponding cocycles $\om$ and $\om',$
then $\af$ is exterior equivalent to $\af'$ \ifo\  %
$\om$ is cohomologous to $\om'.$
This shows that $[\af] \mapsto [\om_{\af}]$
is well defined and injective.
For surjectivity,
given $\om,$ take $\af$ to be conjugation by
the regular $\om$-representation of $G$
on $L^2 (G)$ (or $L^2 (G) \otimes \cH$ if $G$ is finite).
(Continuity of $\af$ follows from Proposition~5(a) of~\cite{mo3}.)
\end{proof}

We write $\K_{\om}$ if we consider $\K$
equipped with an action corresponding
to $\om$ as above.

For an action $\af \colon G \to \Aut (A),$
we denote by $K_*^{\Top} (G; A)$ the left hand side
of the Baum-Connes Conjecture with coefficients~$A.$
See Section~9 of~\cite{BCH},
where it is called $K_*^G (\EG; A).$
As there, we let
$\mu \colon K_*^{\Top} (G; A) \to K_* (A \rtimes_{\af, r} G)$
be the assembly map.
Versions of the following definition
have appeared in~\cite{math, CHMM}.

\begin{defn}\label{defn-twistBC}
Let $[\om] \in H^2 (G, \TT).$
Then the twisted topological $K$-theory of $G$ with respect to
$[\om]$ is defined to be the topological $K$-theory
$K_*^{\Top} (G; \om) = K_*^{\Top} (G; \K_{{\overline{\om}}}).$
The twisted assembly map for $G$ with respect to $\om$
is then defined to be
\[
\mu_{\om} \colon K_*^{\Top} (G; \K_{{\overline{\om}}})
\to K_* (\K_{{\overline{\om}}} \rtimes_r G)
\stackrel{\cong}{\longrightarrow} K_* (C_r^* (G, \om)).
\]
\end{defn}

Up to obvious isomorphisms,
the definition
does not depend on the choice of the representative $\om$
of the class $[\om] \in H^2 (G, \TT)$ and the representative $\alpha$
of the class in $\E_G (\pt)$ corresponding to ${\overline{\om}}.$

Of course, if $G$ is a group which satisfies
the Baum-Connes conjecture with coefficients, then
the twisted assembly map is an isomorphism
for all $\om \in Z^2 (G, \TT).$
But it is sometimes
possible to show that the Baum-Connes conjecture holds
for $\K$ (with respect to any action of $G$ on $\K$)
without knowing that the conjecture holds for all coefficients.
For example, it is shown
in~\cite{CEN} that every almost connected group
satisfies the conjecture for $\K$
but the Baum-Connes conjecture with arbitrary coefficients is not known
in general for those groups.

We now consider homotopies between cocycles.
The definition
of the twisted assembly map suggests
that it is actually useful to consider homotopies
for actions on $\K.$
The first part of the following definition is a special case,
in different language,
of Definition~3.1 of~\cite{Ph_AS}.

\begin{defn}\label{defn-homotopK}
A {\emph{homotopy}} of actions on $\K$
is a $C \big([0, 1] \big)$-linear action
of $G$ on $C \big([0, 1], \, \K \big).$
A homotopy between two classes $[\beta^0], [\beta^1] \in \E_G (\pt)$
is a class $[\beta] \in \E_G \big([0, 1] \big)$ with evaluations
$[\beta^0]$ and $[\beta^1]$ at the points $0$ and $1.$
\end{defn}

\begin{remark}\label{rem-homotopy}
The notion of homotopy in Definition~\ref{defn-homotopK} is,
at least in principle, weaker than the
notion of homotopy between classes
in $H^2 (G, \TT) \cong \E_G (\pt)$ used in
Theorem~\ref{thm-twist}.
Indeed, let $\Om \colon G \times G \to C \big([0, 1], \, \TT \big)$
be a $2$-cocycle with evaluations
$\om_x = \Om ( \, \cdot \, ,  \, \cdot \, ) (x)$ for $x \in [0, 1].$
By Proposition~3.1 of~\cite{HORR},
and with $L_{\om_x}$ being the regular $\om_x$-representation
on $\K (L^2 (G)),$
the formula
\[
\alpha_s (g) (x) = L_{\om_x} (s) g (x) L_{\om_x} (s)^*,
\]
for $g \in C \big([0, 1], \, \K \big),$
defines a \ct\  action
$\alpha \colon G \to \Aut \big( C \big([0, 1], \, \K \big) \big).$
Thus, if $[\alpha_0], [\alpha_1] \in \E_G (\pt)$
correspond to $[\om_0], [\om_1] \in H^2 (G, \TT)$
under the correspondence described in Proposition~\ref{prop-action},
we see that $[\alpha_0]$ and $[\alpha_1]$
are homotopic classes in $\E_G (\pt).$

Using the results of~\cite{EW1, EN} one can show that,
conversely, homotopy of classes
in $\E_G (\pt)$ implies homotopy
of the corresponding cocycles for a very large class
of groups.
This class includes all almost connected
groups by Theorem~1.4 of~\cite{EN} and all
smooth groups in the sense of~\cite{mo2}
by Theorem~5.4 of~\cite{EW1}; note
that it is shown in~\cite{mo2}
that all discrete groups and all compact groups are smooth.
However, it is not clear to us whether this direction holds
in general.

An even weaker notion
of homotopy of classes in $\E_G (\pt) \cong H^2 (G, \TT)$
is obtained by using the topology on $H^2 (G, \TT)$
defined by Moore in~\cite{mo3}
and defining two classes to be homotopic
if they can be connected via a continuous path
$\gamma \colon [0, 1] \to H^2 (G, \TT).$
Again, Theorem~1.4 of~\cite{EN} and Theorem~5.4 of~\cite{EW1}
imply that this notion of homotopy coincides
with the previous notions whenever $G$ is
smooth or almost connected, but we do not know this in general.
The problem is that the topology of $H^2 (G, \TT),$
which is defined via pointwise almost everywhere convergence
of cocycles, can have rather poor separation properties.
\end{remark}

The structure of the group $\E_G (X)$
of all exterior equivalence classes of $C_0 (X)$-linear
actions on $C_0 (X, \K)$ was extensively studied in~\cite{EW1, EN}.
As an easy corollary of the results obtained in~\cite{EW1} we get:

\begin{prop}\label{prop-compact}
Suppose that $G$ is compact and $X$ is contractible.
Then, for any $x \in X,$ the evaluation
map $\E_G (X) \to \E_G (\pt),$ given by $[\alpha] \mapsto[\alpha^x],$
is an isomorphism.
In particular, every $C_0 (X)$-linear action
$\alpha \colon G \to \Aut \big( C_0 (X, \K) \big)$
is exterior equivalent
to the diagonal action $\id_X \otimes \alpha^x$ on
$C_0 (X) \otimes \K \cong C_0 (X, \K).$
\end{prop}

\begin{proof}
Since compact groups are smooth in the sense of Moore
(see the discussion preceding Proposition~3.1
in~\cite{mo2} and also Remark~4.2 of~\cite{EW1}),
it follows from
Theorem~5.4 of~\cite{EW1} and its proof
that there is a short exact sequence
\[
1 \to \check{H}^1 (X; \widehat{G_{\ab}})
   \to \E_G (X) \to C (X, H^2 (G, \TT)) \to 1,
\]
where $\check{H}^1 (X; \widehat{G_{\ab}})$ denotes
\v{C}ech cohomology with coefficients
in the (discrete) dual group of the abelianization
$G_{\ab}$ of $G$ and
$H^2 (G, \TT)$ is topologized as in~\cite{mo3}.
The quotient map in this sequence
sends a class $[\alpha] \in \E_G (X)$ to the mapping
$x \mapsto [\alpha^x] \in \E_G (\pt) \cong H^2 (G, \TT).$
Since $X$ is contractible,
$\check{H}^1 (X; \widehat{G_{\ab}}) = 0.$
Since $G$ is compact,
$H^2 (G, \TT)$ is countable by Corollary~1, on page~56, of~\cite{mo1}.
Since it is locally compact Hausdorff
(see Remark~4.2 of~\cite{EW1}),
$H^2 (G, \TT)$ is discrete.
Since $X$ is connected, it follows that
$C (X, H^2 (G, \TT)) \cong H^2 (G, \TT).$
This completes the proof.
\end{proof}

The following result follows from Theorem~1.5 of~\cite{CEO2}.

\begin{prop}\label{prop-going-down}
Let $G$ be a second countable locally compact group
and let $A$ and $B$ be $G$-\ca s.
Let $z \in KK_0^G (A, B)$ have the property
that for all compact subgroups $L$ of $G,$
the Kasparov product with the restriction
$\res_L^G (z) \in KK_0^L (A, B)$
induces bijective homomorphisms
$KK_*^L (\CC, A) \to KK_*^L (\CC, B).$
Then the Kasparov product with
$z$ induces an isomorphism
$K_*^{\Top} (G; A) \cong K_*^{\Top} (G; B).$
\end{prop}

\begin{proof}
For a closed subgroup $H$ of $G$ and a proper $H$-space $Y$
we define ${\mathcal{F}}_H^n (C_0 (Y)) = KK_n^G (C_0 (Y), \, A)$ and
${\mathcal{G}}_H^n (C_0 (Y)) = KK_n^G (C_0 (Y), \, B).$
Let $\cS(G)$ be the collection of subgroups of $G$
defined in the introduction to Section~1 of~\cite{CEO2}.
Then, by Remark~1.2 of~\cite{CEO2},
the functors
${\mathcal{F}}^* = ({\mathcal{F}}^*_H)_{H \in \cS(G)}$ and
${\mathcal{G}}^* = ({\mathcal{G}}^*_H)_{H \in \cS(G)}$
are going-down functors in the sense of Definition~1.1 of~\cite{CEO2}.
For any
$z \in KK_0^G (A, B),$ the family of transformations
\[
{\mathcal{F}}_H^n (C_0 (Y))
= KK_n^H (C_0 (Y), \, A)
       \stackrel{\otimes_A \res_H^G (z)}{\longrightarrow}
   KK_n^H (C_0 (Y), \, B) = {\mathcal{G}}^n_H (C_0 (Y))
\]
is a going-down transformation in the sense of
Definition~1.4 of~\cite{CEO2}.
Indeed, the assumptions (1) and~(2) of that definition
follow from naturality of the Kasparov product.
Thus, the proof of the proposition will follow
from Theorem~1.5 of~\cite{CEO2} if we know that
\begin{equation}\label{P_GDIsom}
KK_n^L (C_0 (V), \, A)
\stackrel{\otimes_A \res_L^G (z)}{\longrightarrow}
KK_n^L (C_0 (V), \, B)
\end{equation}
is an isomorphism for all compact subgroups $L$ of $G$ and for all
real finite dimensional representation spaces $V$ of $L.$

Let $V$ be a real finite dimensional representation space of $L.$
Lemma~7.7(2) and Lemma~7.7(1) of~\cite{CE1},
in order, give isomorphisms
\[
KK_n^L (C_0 (V), \, D)
\stackrel{\cong}{\longrightarrow}
KK_n^L (C_0 (V) \otimes C_0 (V), \, C_0 (V) \otimes D) \\
\stackrel{\cong}{\longrightarrow}
KK_n^L (\CC, \, C_0 (V) \otimes D)
\]
for any $D,$
which respect Kasparov products.
To prove that~(\ref{P_GDIsom}) is an isomorphism,
it is therefore enough to show that
Kasparov product with $1_{C_0 (V)} \otimes_{\CC} \res_L^G (z)$
induces an isomorphism
\begin{equation}\label{P_GDIsom2}
KK_n^L (\CC, \, C_0 (V) \otimes A)
\cong KK_n^L (\CC, \, C_0 (V) \otimes B).
\end{equation}
We prove this by a bootstrap process.

First, let $L' \subseteq L \subseteq G$ be compact subgroups of $G.$
We claim that multiplication
with $1_{C (L/L')} \otimes_{\CC} \res_{L}^G (z)$ induces isomorphisms
$KK_n^L (\CC, \, C (L/L') \otimes A)
      \cong KK_n^L (\CC, C (L/L') \otimes B).$
For any $G$-\ca\  $D,$ let
$q_D \colon C (L/L') \otimes D \to D$ be the quotient map given by
evaluation at the identity coset $L'\in L/L'.$
By Theorem 20.5.5 of~\cite{Bl}, the chain of maps
\[
KK_n^L (\CC, \, C (L/L') \otimes D)
    \stackrel{\res_{L'}^L}{\longrightarrow}
    KK_n^{L'} (\CC, \, C (L/L') \otimes D)
    \stackrel{q_*}{\longrightarrow}
    KK_n^{L'} (\CC, D)
\]
is an isomorphism.
Moreover, when applied to $A$ and $B$ in place of $D,$ it transforms
Kasparov product with  $1_{C (L/L')} \otimes_{\CC} \res_{L}^G (z)$ to
Kasparov product with $\res_{L'}^G (z).$
The claim thus follows from our assumption on $z.$

Now let $V$ be any real finite dimensional representation space of $L.$
Let
\[
\varnothing = U_0 \subseteq U_1 \subseteq \cdots \subseteq U_l = V
\]
be a stratification as in the proof of Lemma~2.8 of~\cite{CEO2}.
Thus, the $U_j$ are $L$-invariant open subsets of $V$ such that the
difference sets $M_j = U_j \smallsetminus U_{j - 1}$
and their quotients $M_j / L$ are smooth manifolds,
and $M_j$
is locally $L$-homeomorphic to spaces of the
form $X \times L/L_j$ for suitable compact subgroups $L_j \subseteq L.$
Now $M_j / L$ is triangulable by Theorem~10.6 of~\cite{Mnk},
and refining the triangulation allows us to assume that
the inverse image of each open simplex $S$
is $L$-homeomorphic to $S \times L/L_j$
with the trivial action of $L$ on $S.$
By Bott periodicity and the previous paragraph,
(\ref{P_GDIsom2})~is an isomorphism when
$C_0 (V)$ is replaced by $C_0 (S \times L/L').$
Using six term exact sequences and the Five Lemma,
we conclude that~(\ref{P_GDIsom2}) is an isomorphism when
$C_0 (V)$ is replaced by $C_0 (X)$ for the inverse image $X$ in $M_j$
of any finite complex in the triangulation of $M_j / L,$
or the inverse image of the topological interior of one.
Using continuity of equivariant $K$-theory for compact group actions
(Theorem 2.8.3(6) of~\cite{Ph0} or Proposition 11.9.2 of~\cite{Bl}),
we find that~(\ref{P_GDIsom2}) is an isomorphism when
$C_0 (V)$ is replaced by $C_0 (M_j).$
Again using six term exact sequences and the Five Lemma,
we conclude that~(\ref{P_GDIsom2}) is an isomorphism for
our given real representation space~$V.$
\end{proof}

\begin{thm}\label{thm-homotop}
Let
$\alpha \colon G \to \Aut \big(C \big( [0, 1], \, \K \big) \big)$
be a homotopy of actions on $\K.$
Write $\K_x$ for
$\K$ when we consider $\K$ as the fiber at $x$ equipped with the
action $\alpha^x \colon G \to \Aut (\K),$
and let $q_x \colon C \big( [0, 1], \, \K \big) \to \K_x$
denote evaluation at $x.$
Then
\[
(q_x)_* \colon K_*^{\Top} \big( G; \, C \big([0, 1], \, \K \big) \big)
      \to K_*^{\Top} (G; \K_x)
\]
is an isomorphism for all $x \in [0, 1].$
\end{thm}

\begin{proof}
By Proposition~\ref{prop-going-down} it is enough to check that
\[
(q_x)_* \colon KK_*^L \big( \CC, \, C \big([0, 1], \, \K \big) \big)
    \to KK_*^L (\CC, \K_x)
\]
is an isomorphism for every compact subgroup $L \subseteq G.$
By the Green-Julg theorem
(Theorem 11.7.1 of~\cite{Bl} or Theorem 2.8.3(7) of~\cite{Ph0}),
this is equivalent to saying that
\[
(q_x \rtimes L)_*
\colon K_* \big(C \big([0, 1], \, \K \big) \rtimes_{\af} L \big)
\to K_* (\K_x \rtimes_{\af^x} L)
\]
is an isomorphism for all such $L.$
We know from Proposition~\ref{prop-compact} that
the restriction $\res_L^G (\alpha)$ of $\alpha$ to $L$
is exterior equivalent
to $\id_{[0, 1]} \otimes \res_L^G (\alpha^x).$
So
$C \big([0, 1], \, \K \big) \rtimes_{\af} L
\cong C \big([0, 1], \, \K_x \rtimes_{\af^x} L \big),$
and it
is easy to check that $q_x \rtimes L$ corresponds to the evaluation map
$\ev_x \colon C \big([0, 1],
       \, \K_x \rtimes_{\af^x} L \big) \to \K_x \rtimes_{\af^x} L.$
This map induces an isomorphism on $K$-theory because it is
a homotopy equivalence.
\end{proof}

\begin{cor}\label{cor-homotop}
Let $\alpha \colon G \to \Aut \big(C \big( [0, 1], \, \K \big) \big)$
be a homotopy for actions on $\K$
and assume that $G$ satisfies the Baum-Connes conjecture
with coefficients
$C \big([0, 1], \, \K \big)$ (with respect to $\alpha$)
and with coefficients $\K_x$ for
all $x \in [0, 1].$
Then, for any $x \in [0, 1],$
the map
\[
(q_x \rtimes G)_*
        \colon K_* (C \big([0, 1], \, \K \big) \rtimes_{\af, r} G)
       \to K_* (\K_x \rtimes_{\af^x, r} G)
\]
is an isomorphism.
\end{cor}

Corollary~\ref{cor-homotop} fails if $\K$
is replaced by a general AF~algebra,
or by a general commutative \ca,
even for $G = \Z_2.$
See Examples~3.3 and~3.5 of~\cite{Ph_AS}.
Sometimes something can be said; see Theorem~4.3 of~\cite{Ph_AS}.

\begin{thm}[Theorem~\ref{thm-twist}]\label{thm-twist2}
Suppose that $\om_0, \om_1\in Z^2(G, \TT)$ are
homotopic Borel $2$-cocycles
on the second countable locally compact group $G,$
that is, there exists a Borel $2$-cocycle
$\Om \in Z^2 (G, \, C ([0,1], \, \TT))$ such that
$\om_j = \Om (\, \cdot \, , \, \cdot \, ) (j)$ for $j = 0, 1.$
Suppose further
that $G$ satisfies the Baum-Connes
conjecture with coefficients.
Then
$K_* (C_r^* (G, \om_0)) \cong K_* (C_r^* ( G, \om_1)).$
\end{thm}

\begin{proof}
If $G$ satisfies the Baum-Connes conjecture
with arbitrary coefficients,
then $G$ certainly satisfies the hypotheses
of Corollary~\ref{cor-homotop}.
The proof follows by combining this corollary
with the isomorphism~(\ref{eq-iso}) and Remark~\ref{rem-homotopy}.
\end{proof}

\begin{remark}\label{rem-assumptions}
Let us briefly discuss the assumption on the Baum-Connes conjecture
as given in Corollary~\ref{cor-homotop}.
If $G$ satisfies the conjecture with coefficients
$C \big( [0, 1], \, \K \big)$
for {\emph{all}} homotopies
$\alpha \colon G \to \Aut \big( C \big( [0, 1], \, \K \big) \big)$
(we then say that $G$ satisfies the conjecture for homotopies),
then $G$ automatically satisfies the conjecture with coefficients $\K$
for any $G$-action on $\K.$
This follows because for any action $\beta$ of $G$ on $\K,$
the action $\id_{[0, 1]} \otimes \beta$
of $G$ on $C \big([0, 1], \, \K \big)$ is equivariantly
$KK$-equivalent to $\K_{\beta}.$
Conversely, if we just know that $G$ satisfies
the conjecture for $\K$ with respect to any given $G$-action,
then it follows from Proposition~3.1 of~\cite{CEN} that
$G$ satisfies the conjecture for homotopies
as soon as $G$ is exact
and has a $\gamma$-element in the sense of Kasparov.
\end{remark}

For later use it is convenient to restate Corollary~\ref{cor-homotop}
completely in terms of homotopies between cocycles.
We will need the reduced twisted crossed product \ca\  %
$C ([0, 1]) \rtimes_{\Om, r} G$ for a cocycle
$\Om \in Z^2 \big(G, \, C \big([0, 1], \, \TT \big) \big).$
It is a completion of the convolution algebra
$L^1 \big( G, \, C ([0, 1]), \, \Om \big)$ of all $C ( [0, 1] )$-valued
$L^1$-functions on $G$ with convolution given by
\[
f*_{\Om} g (s, x)
= \int_G f (t, x) g (t^{-1} s, \, x) \Om (t, \, t^{-1} s) (x) \, d x.
\]
See~\cite{PR}, and, for the case of a discrete group,~\cite{ZM}.

\begin{cor}\label{cor-twist}
Suppose that $G$ satisfies the Baum-Connes conjecture
with respect to any
homotopy of actions on~$\K.$
Let $G$ act trivially on $C \big([0, 1], \, \TT \big).$
Then for any $\Om \in Z^2 \big(G, \, C \big([0, 1], \, \TT \big) \big)$
and any $x \in [0, 1],$ the canonical quotient map
$p_x \colon C ([0, 1]) \rtimes_{\Om, r} G \to C_r^* (G, \om_x)$
induces an isomorphism on $K$-theory.
\end{cor}

\begin{proof}
As explained in Remark~\ref{rem-homotopy}, the cocycle
$\Om$ determines a homotopy
$\alpha \colon G \to \Aut \big( C \big([0, 1], \, \K \big) \big).$
We claim that, if
${\overline{\af}} \colon
      G \to \Aut \big(C \big([0, 1], \, \K \big) \big)$
is the homotopy
corresponding to the inverse ${\overline{\Om}}$ of $\Om,$
then there is an isomorphism
\[
\big(C ([0, 1]) \rtimes_{\Om, r} G \big) \otimes \K
      \cong C \big([0, 1], \, \K \big) \rtimes_{{\overline{\af}}, r} G.
\]
For the full twisted crossed products,
this follows from~Proposition 4.6 of~\cite{EW2}
(which is a special case of Theorem~3.4 of~\cite{PR}).
Remark~3.12 and the proof of Theorem~3.11 of~\cite{PR}
imply that this is correct for the reduced crossed products as well.

With $\om_x = \Om ( \, \cdot \, ,  \, \cdot \, ) (x)$
for $x \in [0, 1],$
there are obvious quotient maps
$p_x \colon C ([0, 1]) \rtimes_{\Om, r} G \to C_{r}^* (G, \om_x).$
Moreover, with the right vertical arrow
being the isomorphism of~(\ref{eq-iso}),
the diagram
\[
\begin{CD}
\big(C ([0, 1]) \rtimes_{\Om, r} G \big) \otimes \K
    @>p_x \otimes \Id_{\K} >> C_{r}^* (G, \om_x) \otimes \K \\
@V \cong VV     @VV\cong V\\
C \big([0, 1], \, \K \big) \rtimes_{{\overline{\af}}, r} G
@>> q_x \rtimes G>  \K \rtimes_{{\overline{\af}}^x, r} G
\end{CD}
\]
commutes.
Using the induced maps in $K$-theory for this diagram
together with Corollary~\ref{cor-homotop}
and Remark~\ref{rem-assumptions}, the result follows.
\end{proof}

\begin{defn}\label{P_real}
We call a class $[\om] \in H^2 (G, \TT)$ {\emph{real}}
if there exists a cocycle $c \in Z^2 (G, \RR)$ such that
$\om$ is equivalent to the cocycle
$(s, t) \mapsto \exp \big(i c (s, t) \big).$
\end{defn}

In this case we obtain a homotopy
$\Om \in Z^2 \big(G, \, C \big([0, 1], \, \TT \big) \big)$
between $[\om]$ and the trivial cocycle $[1]$ by the formula
\[
\Om (s, t) (x) = e^{i x \cdot c (s, t)}
\]
for $x \in [0, 1]$ and $s, t \in G.$
As a consequence we get the following corollary.
Several special cases of it appear
in~\cite{Ell1, PR1, CHMM, mm1, mm2}.
(In~\cite{math}, Mathai states a theorem
saying that for any discrete $G$ and any real
$2$-cocycle $\om$ on $G,$
one has $K_0 (C_r^* (G, \om)) \cong K_0 (C_r^* (G)).$
Unfortunately, the proof given in~\cite{math} has a substantial gap,
and hence the result is not available so far.)

\begin{cor}\label{cor-real}
Let $G$ be a group which satisfies the Baum-Connes conjecture
for homotopies of actions on $\K,$
and let $\om \in Z^2 (G, \TT)$ be a cocycle such that $[\om]$ is real.
Then $K_* (C_r^* (G, \om)) \cong K_* (C_r^* (G)).$
More generally, if $G$ is any second countable locally compact group
and $[\om] \in H^2 (G, \TT)$
is real, then $K_*^{\Top} (G, \om) \cong K_*^{\Top} (G).$
\end{cor}

Let $\Theta$ be any skew symmetric real $d \times d$ matrix.
Then $\Theta$ determines a $2$-cocycle
on $\Z^d$ via
\[
\om_{\Theta} (m, n) = \exp (\pi i \lk \Theta m, n \rk)
\]
for $m, n \in \Z^d.$
The corresponding $d$-dimensional noncommutative torus~\cite{Rf2}
is the twisted group \ca\  %
$A_{\Theta} = C^* (\Z^d, \om_{\Theta}).$
Equivalently,
$A_{\Theta}$ is the universal \ca\  %
generated by unitaries $u_1, u_2, \dots, u_d$
subject to the relations
\[
u_k u_j = \exp (2 \pi i \te_{j, k} ) u_j u_k
\]
for $1 \leq j, k \leq d.$
(Of course, if all $\te_{j, k}$ are integers, it is not really
noncommutative.)

Since $\om_{\Theta}$ is a real cocycle,
Corollary~\ref{cor-real} immediately gives the following
well known computation of the $K$-theory of a
higher dimensional noncommutative torus,
as found, for example, in \cite{Ell1} and~\cite{Rief-pro}.
Of course, we do not immediately get the additional information
given in \cite{Ell1} and~\cite{Rief-pro}.

\begin{cor}\label{KThOfNCTorus}
For any $d \geq 1$
and any skew symmetric real $d \times d$ matrix $\Te,$
we have
$K_j \big( C^* (\Z^d, \om_{\Theta}) \big) \cong \ZZ^{2^{d - 1}}$
for $j = 0, 1.$
\end{cor}

It is shown in Theorem~3.2 of~\cite{Bk}
(see Section~2 of~\cite{BB} for terminology)
that every class in $H^2 (\ZZ^d, \TT)$ has a
representative of the form
$\om_{\Theta} (n, m) = e^{ i \lk \Theta n, m \rk}$
with $\Theta$ a skew symmetric real $d \times d$ matrix.

\section{Twisted group C*-algebras of semidirect
    products}\label{sec-action}

As mentioned in the introduction,
the rotation algebra $A_{\theta}$ can be realized
as the twisted group algebra $C^* (\ZZ^2, \om_{\theta})$ with
\[
\om_{\theta} \left(
\left( \begin{smallmatrix} n \\ m \end{smallmatrix} \right),
\left( \begin{smallmatrix} n' \\ m' \end{smallmatrix} \right)
\right) = e^{\pi i \theta (n' m - n m')}.
\]
Indeed, if we define
$u_{\theta}
= \delta_{ \left(\begin{smallmatrix} 1 \\ 0 \end{smallmatrix} \right)}$
and
$v_{\theta}
= \delta_{ \left(\begin{smallmatrix} 0 \\ 1 \end{smallmatrix} \right)},$
then $u_{\theta}$ and $v_{\theta}$ are unitary generators
of $C^* (\ZZ^2, \om_{\theta})$ which satisfy
the commutation relation
$v_{\theta} u_{\theta} = e^{2 \pi i \theta} u_{\theta} v_{\theta}.$
This (and other computations below) follows because
$\left(\begin{smallmatrix} n \\m \end{smallmatrix} \right)
   \mapsto \delta_{\left(\begin{smallmatrix} n \\m \end{smallmatrix}
                    \right)}
   \in C^* (\ZZ^2, \om_{\theta})$
is an $\om_{\theta}$-representation of $\ZZ^2.$

As explained in the introduction,
we want to study the crossed products
$A_{\theta} \rtimes_{\alpha} F$
with $F$ a finite subgroup of $\SL_2 (\Z)$
and the action $\alpha^{\theta}$ of $\SL_2 (\ZZ)$ on
on $A_{\theta}$ as given in~(\ref{eq-action}).
We will realize the crossed product
$A_{\theta} \rtimes_{\alpha^{\te}} H,$
for {\emph{any}}
subgroup $H \subseteq \SL_2 (\ZZ),$
as a twisted group algebra
$C^* (\ZZ^2 \rtimes H, \, \widetilde{\om}_{\theta}),$
where $H \subseteq \SL_2 (\ZZ)$ acts on $\ZZ^2$
via matrix multiplication and
$\widetilde{\om}_{\theta}$ is a suitable extension
of $\om_{\theta}$ to $\ZZ^2 \rtimes H.$
A short computation shows that the cocycle
$\om_{\theta} \in Z^2 (\ZZ^2, \TT)$
is invariant under the action of $\SL_2 (\ZZ)$ on $\ZZ^2,$ that is
\[
\om_{\theta}
\left(N \left(\begin{smallmatrix} n \\ m \end{smallmatrix} \right), \,
   N \left(\begin{smallmatrix} n' \\ m' \end{smallmatrix} \right)
                         \right)
= \om_{\theta}
     \left(\left(\begin{smallmatrix} n \\ m \end{smallmatrix} \right),
     \,
     \left(\begin{smallmatrix} n' \\ m' \end{smallmatrix} \right)
\right)
\]
for all $N \in \SL_2 (\ZZ).$
This implies that the action of $\SL_2 (\ZZ)$ on $\ZZ^2$ defines
an action
$\beta^{\theta} \colon \SL_2 (\ZZ)
      \to \Aut \big(C^* (\ZZ^2, \om_{\theta}) \big)$
via
\[
\beta^{\theta}_N (f)
\left(\begin{smallmatrix} n \\m \end{smallmatrix} \right)
   = f \left(N^{-1} \left(\begin{smallmatrix} n \\ m \end{smallmatrix}
                 \right) \right)
\]
for $f$ in the dense subalgebra $l^1 (\ZZ^2, \om_{\theta}).$
It follows from this formula
and the fact that
$\left(\begin{smallmatrix} n \\m \end{smallmatrix} \right)
\mapsto \delta_{\left(\begin{smallmatrix} n \\ m \end{smallmatrix}
                           \right)}
\in C^* (\ZZ^2, \om_{\theta})$
is an $\om_{\theta}$-representation that
\begin{align*}
\beta^{\theta}_N (u_{\theta})
& = \beta_N^{\theta}
     (\delta_{\left(\begin{smallmatrix} 1 \\ 0 \end{smallmatrix}
                           \right)})
    = \delta_{\left(\begin{smallmatrix} n_{1, 1} \\ n_{2, 1}
                 \end{smallmatrix} \right)}
    = {\overline{\om_{\theta}
        \left(\left(\begin{smallmatrix} n_{1, 1} \\ 0 \end{smallmatrix}
                        \right),
        \left(\begin{smallmatrix} 0 \\ n_{2 1} \end{smallmatrix}
            \right) \right)}}
      \delta_{\left(\begin{smallmatrix} n_{1, 1} \\ 0 \end{smallmatrix}
                            \right)}
      \delta_{\left(\begin{smallmatrix} 0 \\ n_{2, 1} \end{smallmatrix}
             \right)} \\
& = e^{\pi i n_{11}n_{21} \theta}
     (\delta_{\left(\begin{smallmatrix} 1 \\ 0 \end{smallmatrix}
              \right)})^{n_{1, 1}}
     (\delta_{\left(\begin{smallmatrix} 0 \\ 1
                \end{smallmatrix} \right)})^{n_{2, 2}}
    = e^{\pi i n_{1, 1} n_{2, 1} \theta}
           u_{\theta}^{n_{1, 1}} v_{\theta}^{n_{2, 1}}
\end{align*}
and, similarly,
$\beta^{\theta}_N (v_{\theta})
   = e^{\pi i n_{1, 2}n_{2, 2} \theta}
         u_{\theta}^{n_{1, 2}} v_{\theta}^{n_{2, 2}},$
so
$\beta^{\theta} = \alpha^{\theta}$ for all $\theta \in [0, 1].$

We recall some facts about semidirect products.
Let $M$ and $H$ be locally compact groups;
in view of our intended application
($M = \ZZ^2$ and $H \subseteq \SL_2 (\ZZ)$),
we write
the group operation on $M$ additively
and the operation on $H$ multiplicatively.
Assume further that $H$ acts continuously on $M$ by automorphisms;
we write the action as $(h, m) \mapsto h m$ for $h \in H$ and $m \in M.$
Then
the semidirect product $M \rtimes H$ is the set $M \times H,$
equipped with the product topology and the
multiplication
$(m, h) \cdot (m', h') = (m + h m', \, h h').$
We let $\mu \colon H \to (0, \infty)$ be the function such that
$\int_M f (m) \, d m = \mu (h) \int_M f (h^{-1} \cdot m) \, d m$
for all $f \in L^1 (M),$
where $d m$ and $d h$ denote left
Haar measures on $M$ and $H.$
Recall that left Haar measure on $M \rtimes H$
is given by the formula
\[
\int_{M \rtimes H} f (m, h) \, d (m, h)
    = \int_H \int_M f (m, h) \mu (h) \, d m \, d h.
\]

We will also need a continuously parametrized version of
the following lemma.
A general statement is awkward to formulate,
but what we actually need (see Remark~\ref{rem-crossed-theta} below)
has essentially the same proof,
so we just refer to the proof.

\begin{lem}\label{lem-twistext}
Suppose that $M \rtimes H$ is a semidirect product
of the second countable locally compact groups $M$ and $H$
and assume that $\om \in Z^2 (M, \TT)$
is invariant under the action of $H$ on $M,$
that is, $\om (h \cdot n, h \cdot m) = \om (n, m)$
for all $n, m \in M, h \in H.$
Then:
\begin{enumerate}
\item\label{lem-twistext-1}
There is a cocycle $\widetilde{\om} \in Z^2 (M \rtimes H, \, \TT)$
defined by
\[
\widetilde{\om} \big( (m, h), (m', h') \big) = \om (m, \, h \cdot m').
\]
\item\label{lem-twistext-2}
There is an action $\alpha \colon H \to \Aut (C^* (M, \om))$
(and similarly for $C_r^* (M, \om)$)
determined by
$\alpha_h (f) (m) = \mu (h) f (h^{-1} m)$
for $f \in L^1 (M, \om).$
\item\label{lem-twistext-3}
There are isomorphisms
\[
C^* (M \rtimes H, \, \widetilde{\om}) \cong C^* (M, \om) \rtimes_{\af} H
\andeqn
C_r^* (M \rtimes H, \, \widetilde{\om})
\cong C_r^* (M, \om) \rtimes_{\af, r} H
\]
given on the level of $L^1$-functions
by
$f \mapsto \Phi (f) \in L^1 \big( H, \, L^1 (M, \om) \big)$
with $\Phi (f) (h) = \mu (h) f (\, \cdot \, , h).$
\end{enumerate}
\end{lem}

\begin{proof}
Part~(\ref{lem-twistext-1}) follows from straightforward computations.
For full twisted group algebras,
Parts (\ref{lem-twistext-2}) and~(\ref{lem-twistext-3})
are very special cases of
Theorem~4.1 of~\cite{PR}.
Since the deduction from that result is quite tedious,
we sketch the proof of~(\ref{lem-twistext-2}) and~(\ref{lem-twistext-3})
in the case of reduced twisted group algebras.
Note first that we have a faithful representation
$L_{\om}^M \colon C_r^* (M, \om) \to B (L^2 (M))$
given by $L_{\om}^M (f) \xi = f*_{\om} \xi$
for $f \in L^1 (M, \om)$ and $\xi \in L^2 (M).$
It is easy to check on $L^1$-functions
that $\alpha_h \colon L^1 (M, \om) \to L^1 (M, \om)$
preserves all algebraic operations.
Define a unitary operator
$U_h \colon L^2 (M) \to L^2 (M)$
by $(U_h \xi) (m) = \sqrt{\mu (h)} \xi (h^{-1} \cdot m).$
We check, as a sample,
that $U_h^* L_{\om}^M (\alpha_h (f)) U_h = L_{\om}^M (f)$
for all $f \in L^1 (M, \om),$ which then implies
that $\alpha_h$ is norm preserving.
Let $\xi \in L^2 (M).$
Then for almost all $n \in M,$
we have, using $H$-invariance of $\om$ at the fourth step
and the definition of $\mu$ at the fifth step,
\begin{align*}
\big( U_h^* L_{\om}^M   &   (\alpha_h (f)) U_h \xi \big) (n)
   = \sqrt{\mu (h^{-1})} \big( \alpha_h (f)*_{\om} U_h (\xi) \big)
                    (h \cdot n) \\
& = \sqrt{\mu (h^{-1})} \int_M \alpha_h (f) (m) U_h (\xi)
     \big(h \cdot n - m \big) \om \big(m, \, h \cdot n - m \big) \, d m
                           \\
& = \mu (h) \int_M f (h^{-1} \cdot m) \xi \big( n - h^{-1} \cdot m \big)
      \om \big( m, \, h \cdot n - m \big) \, d m \\
& = \mu (h) \int_M f (h^{-1} \cdot m) \xi \big( n - h^{-1} \cdot m \big)
       \om \big( h^{-1} \cdot m, \, n - h^{-1} \cdot m \big) \, d m \\
& = \int_M f (m) \xi (n - m) \om (m, \, n - m) \, d m
   = \big( L_{\om}^M (f) \xi \big) (n).
\end{align*}
This proves~(\ref{lem-twistext-2}).
To check the isomorphism of~(\ref{lem-twistext-3})
in the reduced case one first checks,
using Fubini's theorem,
that the map
$\Phi \colon L^1 (M \rtimes H, \, \widetilde{\om})
                      \to L^1 (H, \, L^1 (M, \om))$
is a bijection which preserves
all algebraic operations.
To see that it extends to an isomorphism of the reduced \ca s,
consider the unitary
$V \colon L^2 (M \rtimes H) \to L^2 (H, \, L^2 (M))$ given
by $(V \xi) (h) (n) = \xi (h \cdot n, \, h).$
Let $L_{{\widetilde{\om}}}^{M \rtimes H}$
denote the ${\widetilde{\om}}$-regular representation of $M \rtimes H.$
Let
$\Lambda \colon L^1 (H, \, C_r^* (M, \om))
                   \to B \big(L^2 \big( H, \, L^2 (M) \big) \big)$
denote the integrated from of the regular representation of
$(C_r^* (M, \om), \, H, \, \alpha)$
induced from the faithful representation
$L_{\om}^M \colon C_r^* (M, \om) \to B (L^2 (M)).$
It is given by
\[
\big( \Lambda (g) \eta \big) (h)
= \int_H L_{\om}^M (\alpha_{h^{-1}} (g (l)) \eta (l^{-1}h) \, d l
\]
for $g \in L^1 (H, \, C_r^* (M, \om))$
and $\eta \in L^2 (H, \, L^2 (M)).$
Using this formula, one can check that
\[
V L_{{\widetilde{\om}}}^{M \rtimes H} (f) \xi = \Lambda (\Phi (f)) V \xi
\]
for all $f \in L^1 (M \rtimes H, \, {\widetilde{\om}})$ and
$\xi \in L^2 (M \rtimes H).$
This implies that $\Phi$ preserves the C*-norm.
\end{proof}

\begin{cor}\label{cor-crossed-theta}
Suppose that $H$ is a subgroup of $\SL_2 (\ZZ).$
Define a cocycle on $\ZZ^2 \rtimes H$ by
\[
{\widetilde{\om}}_{\theta}
\left( \big(\left(\begin{smallmatrix} n \\ m \end{smallmatrix} \right),
        \,  N \big), \,
     \big(\left(\begin{smallmatrix} n' \\ m' \end{smallmatrix} \right),
         \,  N' \big) \right)
= \om_{\theta} \left( \left( \begin{smallmatrix} n \\
            m \end{smallmatrix} \right),
            \, N \cdot
             \left(\begin{smallmatrix} n' \\
               m' \end{smallmatrix} \right) \right)
\]
for
$\left(\begin{smallmatrix} n \\ m \end{smallmatrix} \right),
\left(\begin{smallmatrix} n' \\ m' \end{smallmatrix} \right) \in \ZZ^2$
and $N, N' \in H.$
Then there is a canonical isomorphism
\[
A_{\theta} \rtimes_{\alpha^{\theta}, r} H
   = C^* (\ZZ^2, \om_{\theta}) \rtimes_{\alpha^{\theta}, r} H
   \cong C_r^* (\ZZ^2 \rtimes H, \, \widetilde{\om}_{\theta})
\]
(and similarly for the full crossed products).
Moreover, for all $\theta \in [0, 1],$
we have graded isomorphisms
\[
K_* (A_{\theta} \rtimes_{\alpha^{\theta}, r} H)
\cong K_* (C_r^* (\ZZ^2 \rtimes H, \, \widetilde{\om}_{\theta}))
\cong K_* (C_r^* (\ZZ^2 \rtimes H)).
\]
\end{cor}

\begin{proof}
The isomorphism of \ca s follows directly from Lemma~\ref{lem-twistext}.
For the isomorphisms on $K$-theory,
first note that, because $\SL_2 (\ZZ)$
satisfies the Baum-Connes Conjecture with coefficients,
Theorem~2.5 and Corollary~3.14 of~\cite{CE2}
imply that $\ZZ^2 \rtimes H$ does too.
Then use Corollary~\ref{cor-real}
and the fact that $[\widetilde{\om}_{\theta}]$
is real in the sense of Definition~\ref{P_real}.
\end{proof}

\begin{remark}\label{rem-crossed-theta}
In our application,
to the classification of crossed products
$A_{\theta} \rtimes_{\af, r} H$
for certain subgroups $H \subseteq \SL_2 (\ZZ),$
we will need the $K$-groups together with a basis for each.
This problem can be simplified substantially,
by working with continuous families (in a suitable sense)
of projections or unitaries
in matrix algebras over the twisted group algebras
$C_r^* (\ZZ^2 \rtimes H, \, \widetilde{\om}_{\theta}).$
Thus, for any interval $[a, b],$
consider the cocycle
${\widetilde{\Om}} \colon (\ZZ^2 \rtimes H) \times (\ZZ^2 \rtimes H)
         \to C \big([a, b], \, \TT \big)$ given by
${\widetilde{\Om}} ( \, \cdot \, ,  \, \cdot \, ) (\theta)
           = \widetilde{\om}_{\theta}$
for $\theta \in [a, b].$
By a continuous family, we simply mean a projection (or unitary)
in some matrix algebra
$M_n \big( C \big([0, 1] \big)
       \rtimes_{{\widetilde{\Om}}, r} (\ZZ^2 \rtimes H) \big).$
It follows from
Corollary~\ref{cor-twist} that the evaluation map
\[
\ev_{\theta}
\colon C \big([0, 1] \big)
      \rtimes_{{\widetilde{\Om}}, r} (\ZZ^2 \rtimes H)
\to C_r^* (\ZZ^2 \rtimes H, \, \widetilde{\om}_{\theta})
\]
induces an isomorphism in $K$-theory,
and hence it maps a basis to a basis.
Lemma~\ref{lem-twistext} provides an isomorphism
$C_r^* (\ZZ^2 \rtimes H, \, \widetilde{\om}_{\theta})
    \cong A_{\theta} \rtimes_{\alpha^{\theta}, r} H.$
An obvious generalization, proved by an easy modification of its proof,
provides an isomorphism
\[
C \big([a, b] \big) \rtimes_{\widetilde{\Om}, r} (\ZZ^2 \rtimes H)
    \cong
      \big( C \big([a, b] \big) \rtimes_{{\Om}, r} \ZZ^2 \big) \rtimes H
\]
which moreover respects the evaluation maps $\ev_{\te}$ on both
algebras.
In particular, for \pj s $p_1, p_2, \ldots, p_n$ in matrix
algebras over
$\big( C \big([a, b] \big) \rtimes_{{\Om}, r} \ZZ^2 \big) \rtimes H,$
the following are equivalent:
\begin{enumerate}
\item\label{rem-crossed-theta-1}
$[p_1], [p_2], \ldots, [p_n]$ form a basis for
$K_0 \big(
\big( C \big([a, b] \big) \rtimes_{{\Om}, r} \ZZ^2 \big)
                              \rtimes H \big).$
\item\label{rem-crossed-theta-2}
For some $\te \in [a, b],$ the evaluated classes
$[\ev_{\te} (p_1)], \, [\ev_{\te} (p_2)],
\, \ldots, \, [\ev_{\te} (p_n)]$
form a basis for
$K_0 (A_{\theta} \rtimes_{\alpha^{\theta}, r} H).$
\item\label{rem-crossed-theta-3}
For all $\te \in [a, b],$ the evaluated classes
$[\ev_{\te} (p_1)], \, [\ev_{\te} (p_2)],
\, \ldots, \, [\ev_{\te} (p_n)]$
form a basis for
$K_0 (A_{\theta} \rtimes_{\alpha^{\theta}, r} H).$
\end{enumerate}
\end{remark}

\begin{remark}\label{lem-PR}
We should also remark at this point that twisted group algebras
of semidirect products
$\Gamma = \ZZ^2 \rtimes_A \ZZ,$ with $A$
a fixed element of $\SL_2 (\ZZ)$ and the action of
$\ZZ$ on $\ZZ^2$ given by powers of $A,$
were extensively studied by Packer and Raeburn
in~\cite{PR1}.
They use the embedding of $\Gamma$ into
$\RR^2 \rtimes_A \RR$ together with Connes's Thom isomorphism
to show
that $K_* (C^* (\Gamma, \om)) \cong K_* (C^* (\Gamma))$
for any cocycle
$\om \in H^2 (\Gamma, \TT).$
(See Corollary~2.12 of~\cite{PR1}.)
Since they also show (Corollary~3.6 of~\cite{PR1})
that every cocycle of such a group $\Gamma$ is the
exponential of a real cocycle,
this result also follows from
Corollary~\ref{cor-real}
and the proof of Corollary~\ref{cor-crossed-theta}.
\end{remark}

\section{The untwisted case}\label{sec-untwist}

As mentioned above, we are interested in applying
the results of the previous section to
the case of finite groups $F \subseteq \SL_2 (\ZZ)$
which are conjugate
to one of the groups $\ZZ_2, \ZZ_3, \ZZ_4, \ZZ_6$
with generators as given in~(\ref{eq-gen}).
In this section we give an explicit computation of the $K$-theory
of the untwisted group algebras
$C^* (\ZZ^2 \rtimes F).$
We need some notation.

\begin{ntn}\label{P_SgpNtn}
Let $G$ be a {\emph{discrete}} group.
\begin{enumerate}
\item\label{P_SgpNtn_1}
For a subgroup $H \subseteq G,$
let $\ind_H^G \colon K_* (C^*_r (H)) \to K_* (C^*_r (G))$
be the induction homomorphisms induced by
the obvious maps
of the reduced group \ca s $C^*_r (H) \to C^*_r (G).$
\item\label{P_SgpNtn_2}
For a subgroup $H \subseteq G$ such that the index
$[G : H]$ is finite, the restriction homomorphism
$\res_G^H \colon K_n (C^*_r (G)) \to K_n (C^*_r (H))$
(sometimes written $\res_H^G$)
is given by restriction of scalars:
every finitely generated projective $C^*_r (G)$-module
is also finitely generated and projective as a module
over the subalgebra $C^*_r (H).$
\item\label{P_SgpNtn_3}
We denote
by $\widetilde{K}_n (C^*_r (G))$
the cokernel of the map
$K_n (\CC) = K_n (C^* (\{1\}) ) \to K_n (C^*_r (G))$
induced by the inclusion of the trivial
subgroup $\{1\} \to G,$
and we let
$\pr \colon K_n (C^*_r (G)) \to \widetilde{K}_n (C^*_r (G))$
be the canonical projection.
\end{enumerate}
\end{ntn}

The following lemma, known as the double
coset formula, will be used several times.
It actually
holds in greater generality,
but all we need here is the version below, which
follows from the arguments
in the proof of
Proposition~5.6(b) in Chapter~III of~\cite{Br}.

\begin{lem}\label{lem-doublecoset}
Suppose that $G$ is a discrete
group and $H$ and $C$ are
subgroups of $G$ such that
$C$ has finite index in $G.$
For $g \in G,$ let
$\varphi_g \colon C_r^* (H) \to C_r^* (g H g^{-1})$
denote the isomorphism induced by conjugation.
Then
\[
\res_G^C\circ \ind_H^G=
\sum_{[g]\in C\backslash G/H}
\left(\ind_{gHg^{-1}\cap C}^C\;\circ \;
\res_{g H g^{-1}}^{g H g^{-1} \cap C} \; \circ \; (\varphi_g)_* \right)
\]
as maps from $K_0 (C_r^* (H) )$ to $K_0 (C_r^* (C)).$
\end{lem}

\begin{remark}\label{P_Split}
Let $H$ be a finite group.
Then the \hm\  %
$\ind_{\{ 1 \}}^H
      \colon K_0 (\CC) = K_0 (C^* (\{ 1 \})) \to K_0 (C^* (H))$
has a left inverse $r \colon K_0 (C^* (H)) \to K_0 (\CC).$
Letting $\CC$ carry the trivial representation of $H,$
it sends the class $[V]$
of a finite-dimensional $H$-representation $V$ to
$\dim_{\CC} (\CC \otimes_{C^* (H)} V) \cdot [\CC].$
(The number $\dim_{\CC} (\CC \otimes_{C^* (H)} V)$ is just
the multiplicity of the trivial representation in
the decomposition of $V$ into irreducibles.)
It induces a splitting
\begin{equation}\label{eq-s}
s \colon \widetilde{K}_0 (C^* (H)) \to K_0 (C^* (H))
\end{equation}
of
$\pr \colon K_0 (C^* (H)) \to \widetilde{K}_0 (C^* (H)).$
Namely, $s$ sends $x \in \widetilde{K}_0 (C^* (H))$
to $y - \big( \ind_{\{1\}}^H \circ r \big) (y)$
for any choice of $y \in K_0 (C^*_r (G))$ with $\pr (y) = x.$
\end{remark}

\begin{ntn}\label{P_bHN}
Let $G$ be a discrete group, and let
$H \subseteq G$ be a finite subgroup.
We denote by $b_H$ the composition
\[
b_H \colon \widetilde{K}_0 (C^* (H))
\stackrel{s}{\longrightarrow}
K_0 (C^* (H)) \xrightarrow{\ind_H^G} K_0 (C_r^* (G)).
\]
\end{ntn}

\begin{thm}\label{thm-Wolfgang}
Suppose that $F \subseteq \SL_2 (\ZZ)$
is one of the groups $\ZZ_2, \ZZ_3, \ZZ_4, \ZZ_6$
with generators as in~(\ref{eq-gen}).
Then there exist only finitely many conjugacy classes
of maximal finite subgroups $M \subseteq \ZZ^2 \rtimes F.$
Let $M_0, M_1, \ldots, M_l$
be a list of representatives for these conjugacy classes.
Let
$b_j \colon \widetilde{K}_0 (C^* (M_j)) \to K_0 (C^* (\ZZ^2 \rtimes F))$
for $j = 0, 1, \ldots, l$
be the map $b_{M_j}$ of Notation~\ref{P_bHN}.
Then the following sequence is exact:
\begin{multline*}
0 \to K_0 (C^* (\{1\})) \oplus
      \left( \bigoplus_{j = 0}^l \widetilde{K}_0 (C^* (M_j)) \right)
   \xrightarrow{\ind_{\{1\}}^G \oplus
     \left( \bigoplus_{j = 0}^l b_j \right)} K_0 (C^* (\ZZ^2 \rtimes F))
\\
\xrightarrow {\pr \circ \res_{\ZZ^2 \rtimes F}^{\ZZ^2}}
              \widetilde{K}_0 (C^* (\ZZ^2)) \to 0.
\end{multline*}
Moreover we have $K_1 (C^* (\ZZ^2 \rtimes F)) = \{ 0 \}$
for all $F$ as above.
\end{thm}

\begin{remark}\label{rem-map}
The image of the homomorphism
\[
K_0 (C^* (\{1\})) \oplus
    \left( \bigoplus_{j = 0}^l \widetilde{K}_0 (C^* (M_j)) \right)
\xrightarrow{\ind_{\{1\}}^G \oplus
    \left( \bigoplus_{j = 0}^l b_j \right)} K_0 (C^* (\ZZ^2 \rtimes F))
\]
does not depend on the special choice of the section
{\mbox{$s_j \colon \widetilde{K}_0 (C^* (M_j)) \to K_0 (C^* (M_j))$}}
in~(\ref{eq-s}) which was used in the definition of the maps
$b_j \colon \widetilde{K}_0 (C^* (M_j))
                  \to K_0 (C^* (\ZZ^2 \rtimes F)).$
So the theorem remains
true if we replace it by any other section
$\widetilde{K}_0 (C^* (M_j)) \to K_0 (C^* (M_j)).$
\end{remark}

We give two proofs of Theorem~\ref{thm-Wolfgang}, one depending on the
constructions of Sam Walters~\cite{W1, W2}
and some computations from Section~\ref{sec-gen} below,
and the other more topological,
although it has much in common with the first.
Before doing so, however, we should describe what one can get
from known results on the structure of
$A_{\theta} \rtimes_{\af} F$ for $\te$ rational.
Suppose that $\te = p / q$ in lowest terms, with $q > 0.$
Then the fixed point algebras $A_{\te}^F$ and the
crossed products $A_{\te} \rtimes_{\alpha} F$ have been
computed explicitly,
in a form from which it is easy to determine the $K$-theory.
In particular, this covers the case $C^* (\ZZ^2 \rtimes F).$
For $F = \ZZ_2,$ see Theorems~1.2 and~1.3 of~\cite{BEEK2};
for $F = \ZZ_3,$ see Theorem~3.4.1 of~\cite{P_FW6} for the fixed
point algebra and the Theorem on page~244 for the crossed product;
for $F = \ZZ_4,$ see Theorem~3.2.5 of~\cite{FW1} for the fixed
point algebra and the Theorem~6.2.1 for the crossed product;
and for $F = \ZZ_6,$ see Theorem~3.2.8 of~\cite{P_FW5} for the fixed
point algebra and the Theorem on page~2 for the crossed product.
(The proofs for the crossed products are similar to the proofs
for the fixed point algebras,
and in the cases $F = \ZZ_3$ and $F = \ZZ_6$
they are not actually given in the papers.)

We describe briefly what happens for $F = \ZZ_6$; the other cases
are similar, but sometimes simpler.
Identify the algebras $(M_{3 q})^2 = M_{3 q} \oplus M_{3 q},$
$(M_{2 q})^3,$
and $(M_q)^6$ with algebras of
block diagonal matrices in $M_{6 q}.$
Choose three distinct points $x_0, x_1, x_2 \in S^2.$
Then
\begin{equation}\label{FWIsom}
\begin{split}
& A_{p / q} \rtimes_{\alpha} \ZZ_6
\\
& \hspace*{1em} {\mbox{}}
   \cong \big\{ f \in C (S^2, M_{6 q}) \colon
     {\mbox{$f (x_0) \in (M_{q})^6,$ $f (x_1) \in (M_{2 q})^3,$
                 and $f (x_2) \in (M_{3 q})^2$}} \big\}.
\end{split}
\end{equation}
The fixed point algebra has a similar description,
in which the functions take values in $M_q,$
and the subalgebras $(M_{3 q})^2,$ $(M_{2 q})^3,$ and $(M_{q})^6$
are replaced by algebras with complicated formulas for the
summands and, if $q \leq 6,$ there are sometimes fewer summands.
In particular, if $q = 1,$ then the fixed point algebra is $C (S^2).$
Note that the isomorphism class of an algebra
given by such a description does not depend on the choice of
the points in $S^2$ or the particular identifications of the
direct sums with subalgebras of $M_{6 q}$ or $M_q.$

It follows immediately that $\TT^2 / \ZZ_6 \cong S^2$
(and similarly in the other cases),
a fact that will be proved by much less computational methods
in the course of the proof of Theorem~\ref{thm-Wolfgang}.
It is also relatively easy to see that,
for $\te = p / q,$ we have
\begin{equation}\label{P_RatKThy}
K_0 (A_{\te} \rtimes_{\alpha} \ZZ_6) \cong \ZZ^{10}
\andeqn
K_1 (A_{\te} \rtimes_{\alpha} \ZZ_6) = 0,
\end{equation}
and to write down nine \pj s
in the algebra on the right hand side of~(\ref{FWIsom})
such that their classes,
together with the Bott class, form a basis for
$K_0 (A_{p / q} \rtimes_{\alpha} \ZZ_6).$
(It was this computation for the fixed point algebra,
and the similar results
for $\ZZ_2,$ $\ZZ_3,$ and $\ZZ_4,$
that led to the conjecture that the fixed point algebras are~AF.)
It now follows from Remark~\ref{rem-crossed-theta},
as in the proof of Theorem~\ref{thm-generators},
that~(\ref{P_RatKThy}) holds
for all values of $\te.$
This is enough to make the proof of Theorem~\ref{T:CrPrdIsAF}
go through,
and conclude that $A_{\te} \rtimes_{\alpha} \ZZ_6$ is an AF~algebra.
But it does not allow us to decide which AF~algebra.
Probably the most serious difficulty
is that we do not know the trace of the last generator of
$K_0 (A_{\te} \rtimes_{\alpha} \ZZ_6)$
at irrational values of $\te.$
To find it, we presumably must compute,
for some rational value of $\te$ (such as $\te = 1$),
the $K_0$-class of the module $\E_{\theta}^F$
of Notation~\ref{P_ModNtn} in terms of
the algebra on the right hand side of~(\ref{FWIsom}).

\begin{proof}[First proof of Theorem~\ref{thm-Wolfgang}]
Put $G = \ZZ^2 \rtimes F.$
The group $G$ is amenable.
Thus, $C^* (G) \to C_r^* (G)$ is an isomorphism, and we
write $C^* (G)$ throughout.
Being amenable,
$G$ satisfies the Baum-Connes Conjecture by Corollary~9.2 of~\cite{HK},
that is, the assembly map
\[
\mu \colon K_n^G ({\underline{E}} G) \xrightarrow{\cong} K_n (C^* (G))
\]
is an isomorphism for all $n \in \ZZ.$
Let
$p \colon K_0^G ({\underline{E}} G)
        \to K_0 (G \backslash {\underline{E}} G))$
be the natural map.
The sequence $1 \to \ZZ^2 \to G \to F \to 1$ has the property that
the conjugation action of $F$ on $\ZZ^2$
is free away from the origin $0 \in \ZZ^2.$
We now claim that there is an isomorphism
\begin{equation}\label{P_BCIsom1}
K_1 (C^* (G)) \cong K_1 (G \backslash {\underline{E}} G)
\end{equation}
and an exact sequence
\begin{eqnarray}\label{DL_first_sequence}
0 %
\to \bigoplus_{j = 1}^l \widetilde{K}_0 (C^* (M_j))
\xrightarrow{ \bigoplus_{j = 0}^l b_j} K_0 (C^* (G))
\xrightarrow{p \circ \mu^{-1}} K_0 (G \backslash {\underline{E}} G)
\to 0.
\end{eqnarray}
To get these,
apply Theorem~5.1(a) and Remark~5.2 of~\cite{Davis-Lueck (2003)}
to the extension
$1 \to \ZZ^2 \to G \to F \to 1,$
or apply the more general
Theorem~1.6 of~\cite{Lc}, taking there $K = \{ 1 \}$ and $G = Q.$
(Theorem~1.6 of~\cite{Lc} actually gives a long exact sequence,
which breaks up into short exact sequences after tensoring with $\QQ.$
However, every group appearing in the long exact sequence
is torsion free.)

Clearly $G \backslash {\underline{E}} G$
is the same as $F \backslash \TT^2$ for the obvious
$F$-action on the two-dimensional torus
$\TT^2 = \ZZ^2 \backslash {\underline{E}} G = \ZZ^2 \backslash E \ZZ^2.$
Since $F$ is a subgroup of $\SL_2 (\ZZ),$
its action on $\TT^2$ is orientation preserving.
Also,
there are only finitely many points of $\TT^2$ at which the action is
not free.
Analyzing these, for example as on page~407 of~\cite{Sc},
one can show that
$F \backslash \TT^2$ is a compact $2$-dimensional manifold.
The rational cohomology
$H^* (G \backslash {\underline{E}} G; \, \QQ)$
agrees with $H^* (\TT^2; \QQ)^F.$
Since the action is orientation preserving,
$F$ acts trivially on $H^p (\TT^2; \QQ)$
for $p = 0, 2.$
Since $F$ acts freely on $\ZZ^2 = H_1 (\TT^2; \ZZ)$
away from $\{ 0 \},$
we conclude
$H^1 (\TT^2; \QQ)^F
   \cong {\mathrm{Hom}}_{\ZZ} \big( H_1 (\TT^2; \ZZ)^F, \, \QQ \big)
   = \{ 0 \}.$
Hence $G \backslash {\underline{E}} G = F \backslash \TT^2$
has the rational cohomology of $S^2$
and hence is homeomorphic to $S^2.$
This implies that $K_0 (G \backslash {\underline{E}} G) \cong \ZZ^2$
and $K_1 (G \backslash {\underline{E}} G) = 0.$
We conclude from~\eqref{P_BCIsom1} that
$K_1 (C^* (G)) = 0.$

Using the fact that
$G \backslash {\underline{E}} G$ is $2$-dimensional,
there is an edge homomorphism
$\edge \colon K_0 (G \backslash {\underline{E}} G)
            \to H_2 (G \backslash {\underline{E}} G)$
which comes from the
Atiyah-Hirzebruch spectral sequence converging to
$K_{p + q} (G \backslash {\underline{E}} G)$
with $E^2$-term
$E^2_{p, q} = H_p (G \backslash {\underline{E}} G, \, K_q (\pt)).$
(See, for example, Section~13.6 of~\cite{Wh} for
the Atiyah-Hirzebruch spectral sequence.)
Let $f \colon G \backslash {\underline{E}} G \to \pt$
be the projection onto the one-point space $\pt.$
An easy spectral sequence argument shows that
\[
K_0 (f) \oplus \edge \colon K_0 (G \backslash {\underline{E}} G)
\xrightarrow{\cong}
   K_0 (\pt) \oplus H_2 (G \backslash {\underline{E}} G; \, \ZZ)
\]
is an isomorphism.
The composition
\[
K_0 (C^* (\{1\}))
\xrightarrow{\ind_{1}^G} K_0 (C^* (G))
\xrightarrow{\mu^{-1}} K_0^G ({\underline{E}} G)
\stackrel{p}{\longrightarrow} K_0 (G \backslash {\underline{E}} G)
\stackrel{f}{\longrightarrow} K_0 (\pt)
\]
is an isomorphism, namely the inverse of the assembly map
$K_0 (\pt) \to K_0 (\CC)$ of the assembly map for the trivial
group.
Denote by $e_G$ the composition
\begin{equation}\label{P_eq_eG}
e_G \colon K_0 (C^* (G))
\xrightarrow{\mu^{-1}} K^G_0 ({\underline{E}} G)
\stackrel{p}{\longrightarrow} K_0 (G \backslash {\underline{E}} G)
\xrightarrow{\edge} H_2 (G \backslash {\underline{E}} G; \, \ZZ).
\end{equation}
Now the exact sequence
\eqref{DL_first_sequence} yields the exact sequence
\begin{multline}
0 \to K_0 (C^* (\{1\})) \oplus
   \bigg( \bigoplus_{j = 1}^l \widetilde{K}_0 (C^* (M_j)) \bigg)
\xrightarrow{\ind_{\{ 1 \}}^G \oplus
                  \left( \bigoplus_{j = 0}^l b_j \right)}
        K_0 (C^* (G))
\\
\xrightarrow{e_G} H_2 (G \backslash {\underline{E}} G; \, \ZZ) \to 0.
\label{DL_second_sequence}
\end{multline}

Let $\res_{G}^{\ZZ^2} \colon K_0 (C^* (G)) \to
K_0 (C^* (\ZZ^2))$ be the homomorphism induced by restriction to
the subgroup of finite index.
By the double coset formula of Lemma~\ref{lem-doublecoset}
and because each conjugate of each finite group $M_j$ has trivial
intersection with $\ZZ^2$ in $G,$
the composition
\[
K_0 (C^* (\{1\})) \oplus
      \left( \bigoplus_{j = 0}^l \widetilde{K}_0 (C^* (M_j)) \right)
\xrightarrow{\ind_{\{1\}}^G \oplus
           \left( \bigoplus_{j = 0}^l b_j \right)}
     K_0 (C^* (G))
\xrightarrow{\res_{G}^{\ZZ^2}} K_0 (C^* (\ZZ^2))
\]
factors through the map
$\ind_{\{1\}}^{\ZZ^2} \colon K_0 (C^* (\{ 1 \})) \to K_0 (C^* (\ZZ^2)).$
Combining this with
the exact sequence~(\ref{DL_second_sequence}) yields a commutative
diagram
\[
\begin{CD}
0 && 0
\\
@VVV @VVV
\\
K_0 (C^* (\{1\})) \oplus \left( \bigoplus_{j = 0}^l
   \widetilde{K}_0 (C^* (M_j)) \right)
@> \id >>
K_0 (C^* (\{1\})) \oplus \left( \bigoplus_{j = 0}^l
   \widetilde{K}_0 (C^* (M_j)) \right)
\\
@V\ind_{\{1\}}^G \oplus \left( \bigoplus_{j = 0}^l b_j \right) VV
@VV \ind_{\{1\}}^G \oplus \left( \bigoplus_{j = 0}^l b_j \right) V
\\
K_0 (C^* (G)) @> \id >> K_0 (C^* (G))
\\
@V e_G VV @VV \pr \circ \res_{G}^{\ZZ^2} V
\\
H_2 (G \backslash {\underline{E}} G; \, \ZZ)
@> \omega >>
\widetilde{K}_0 (C^* (\ZZ^2))
\\
@VVV @VVV
\\
0 && 0
\end{CD}
\]
for an appropriate map $\omega.$
Lemma~\ref{lem-gen} below
(based on Sam Walters' construction of Fourier modules~\cite{W1, W2};
the proof does not depend on Theorem~\ref{thm-Wolfgang})
provides explicit elements in $K_0 (C^* (G))$
which are mapped to a generator in
$\widetilde{K}_0 (C^* (\ZZ^2))$ under $\pr \circ \res_G^{\ZZ^2}.$
This implies that $\pr \circ \res_G^{\ZZ^2}$ is surjective.
Hence $\omega$ is a surjection of infinite cyclic groups
and therefore an isomorphism.
Since the left column in the diagram above is exact, the
same is true for the right column.
\end{proof}

We now give a second proof,
not relying on Lemma~\ref{lem-gen}.
Parts of the argument are the same as in the first proof,
and we refer freely to it.

\begin{proof}[Second proof of Theorem~\ref{thm-Wolfgang}]
Analogously to the definition of $e_G$ in~(\ref{P_eq_eG}),
we define a surjection
\[
e_{\ZZ^2} \colon K_0 (C^* (\ZZ^2)) \to H_2 (\ZZ^2 \backslash E \ZZ^2)
\]
for the subgroup $\ZZ^2 \subseteq G.$
Let
$q \colon \TT^2
= \ZZ^2 \backslash {\underline{E}} \ZZ^2
\to F \backslash \TT^2 = G \backslash {\underline{E}} G$
be the projection.
Let $i \colon \ZZ^2 \to G$ be the inclusion.
Then the following diagram commutes:
\begin{eqnarray}\label{square for e_G and e_ (Z^2)}
\comsquare{K_0 (C^* (\ZZ^2))}{e_{\ZZ^2}}{H_2 (\TT^2; \ZZ)}{K_0 (i)}
{H_2 (q)}{K_0 (C^* (G))}{e_G}{H_2 (F \backslash \TT^2; \, \ZZ).}
\end{eqnarray}
The $F$-action on $\TT^2$
is orientation preserving and has at least one free orbit.
The quotient $F \backslash \TT^2$ is a two dimensional manifold
as in the first proof,
so $q \colon \TT^2 \to F \backslash \TT^2$
is a map of closed oriented surfaces
of degree $\card (F)$ and thus induces an injective map
$H_2 (q)$ of infinite cyclic groups
whose image has index $\card (F).$

An easy spectral sequence argument shows that
the image of
\[
\ind_{\{1\}}^{\ZZ^2} \colon K_0 (C^* (\{1\})) \to K_0 (C^* (\ZZ^2))
\]
is the kernel of the surjection
$e_{\ZZ^2} \colon K_0 (C^* (\ZZ^2)) \to H_2 (\TT^2; \ZZ).$
Hence we obtain
from \eqref{DL_second_sequence} a map
$\mu \colon H_2 (F \backslash \TT^2; \, \ZZ) \to H_2 (\TT^2; \ZZ)$
such that the following square commutes:
\begin{eqnarray}\label{square defining mu}
\comsquare{K_0 (C^* (G))}{e_G}{H_2 (F \backslash \TT^2; \, \ZZ)}
{\res_G^{\ZZ^2}}{\mu}{K_0 (C^* (\ZZ^2))}{e_{\ZZ^2}}{H_2 (\TT^2; \ZZ),}
\end{eqnarray}
and an isomorphism
\[
\nu \colon H_2 (\TT^2; \ZZ)
\xrightarrow{\cong} \widetilde{K}_0 (C^* (\ZZ^2))
\]
such that the composition $\nu \circ e_{\ZZ^2}$ agrees with the
projection $\pr \colon K_0 (C^* (\ZZ^2)) \to
\widetilde{K}_0 (C^* (\ZZ^2)).$
Hence the following diagram
commutes:
\begin{eqnarray}\label{some square}
\comsquare{K_0 (C^* (G))}{e_G}{H_2 (F \backslash \TT^2; \, \ZZ)}
{\res_G^{\ZZ^2}}{\nu \circ \mu}
{K_0 (C^* (\ZZ^2))}{\pr}{\widetilde{K}_0 (C^* (\ZZ^2)).}
\end{eqnarray}
In view of the exact sequence~\eqref{DL_second_sequence}
it remains to show
that $\mu$ is bijective.
Again by the double coset formula (Lemma~\ref{lem-doublecoset}), the
composition
\[
K_0 (C^* (\ZZ^2))
\xrightarrow{K_0 (i)} K_0 (C^* (G))
\xrightarrow{\res_{G}^{\ZZ^2}} K_0 (C^* (\ZZ^2))
\]
is $\sum_{j = 1}^{\card (F)} K_0 (\sigma)^j,$
where $\sigma \colon \ZZ^2 \to \ZZ^2$ is given by
multiplication by the generator of the finite cyclic group $F.$
The map $\sigma$ induces an orientation
preserving automorphism of $\TT^2$ and hence the identity $H_2 (\TT^2).$
Therefore the following diagram commutes:
\begin{eqnarray}\label{square relating res circ i and |F| cdot id}
\comsquare{K_0 (C^* (\ZZ^2))}{e_{\ZZ^2}}{H_2 (\TT^2; \ZZ)}
{\res_{G}^{\ZZ^2} \circ K_0 (i)}{\card (F) \cdot \id}
{K_0 (C^* (\ZZ^2))}{e_{\ZZ^2}}{H_2 (\TT^2; \ZZ).}
\end{eqnarray}
Fix a generator $x$ of the infinite cyclic group $H_2 (\TT^2; \ZZ).$
Choose $y \in K_0 (C^* (\ZZ^2))$
with $e_{\ZZ^2} (y) = x.$
Let $z \in K_0 (C^* (G))$ be the image
of $y$ under the homomorphism
$K_0 (i) \colon K_0 (C^* (\ZZ^2)) \to K_0 (C^* (G)).$
We conclude from \eqref{square for e_G and e_ (Z^2)}
and the fact
that $H_2 (\pr)$ is injective and has cokernel of order $\card (F)$ that
$e_G (z)$ is $\card (F) \cdot a$ for a suitable generator $a$
of the infinite cyclic group
$H_2 (F \backslash \TT^2; \, \ZZ).$
We get from~\eqref{square defining mu}
and \eqref{square relating res circ i and |F| cdot id}
that
\[
\card (F) \cdot \mu (a) = \mu \circ e_G (z) = \card (F) \cdot x.
\]
Since $a$ and $x$ are generators
of infinite cyclic groups, $\mu$ is bijective.
(Notice that the map $\omega$ agrees with the map $\nu \circ \mu$
in the first proof.)
This finishes the second proof of Theorem~\ref{thm-Wolfgang}.
\end{proof}

It is clear how to use Theorem~\ref{thm-Wolfgang}
to give an explicit basis for $K_0 (\ZZ^2 \rtimes F),$
and we do this in all four cases.
Note first that $K_0 (C^* (\ZZ_k))$
is the free abelian group of rank $k.$
If $t$ is a generator for $\ZZ_k,$
identified with its canonical image in $C^* (\ZZ_k),$
and $\zt = e^{2 \pi i / k},$
then an explicit basis is given by the classes of the projections
\[
p_j = \frac{1}{k} \sum_{l = 0}^{k - 1} \big( \zeta^j t \big)^l
\]
for $j = 0, \ldots, k - 1.$
One easily checks that $\sum_{j = 0}^{k - 1} p_j = 1 \in C^* (\ZZ_k),$
and therefore
we may replace any one of the above projections by $[1]$
to obtain a new basis.
In particular, it follows that $\widetilde{K}_0 (C^* (\ZZ_k)),$
as defined in Notation~\ref{P_SgpNtn}(\ref{P_SgpNtn_3}),
is a free abelian group of rank $k - 1$
with basis
\[
\big( [p_0], \, [p_1], \, \ldots, \, [p_{k - 2}] \big).
\]

In the following examples we give a complete set of
representatives of the conjugacy classes of maximal finite
subgroups.
The elementary calculation that this is indeed such a
system is carried out for $F = \ZZ_4$ in
Lemma~2.2 of~\cite{Lc}.
The analogous calculations in the
other cases are left to the reader.

\begin{ex}\label{ex-K0}
{\textbf{(a) The case $F = \ZZ_2.$}}
The group $\ZZ^2 \rtimes \ZZ_2$ has the presentation
\[
\ZZ^2 \rtimes \ZZ_2
= \big\lk u, v, t \colon
    t^2 = 1, \; u v = v u, \; t u t = u^{-1}, \; t v t = v^{-1} \big\rk.
\]
The maximal finite subgroups are represented, up to conjugacy,
by the four groups
\[
M_0 = \lk t \rk, \quad M_1 = \lk u t \rk,
\quad M_2 = \lk v t \rk,
\quad M_3 = \lk u v t \rk,
\]
which are all isomorphic to $\ZZ_2.$
Hence it follows from Theorem~\ref{thm-Wolfgang} that,
for any $S \in K_0 (C^* (\ZZ^2 \rtimes \ZZ_2))$
which maps onto a generator of
$\widetilde{K}_0 (C^* (\ZZ^2)) \cong \ZZ$
via $\pr \circ \res^{\ZZ^2}_{\ZZ^2 \rtimes \ZZ_2},$
\[
\big( [1], \;
\big[\tfrac{1}{2} (1 + t) \big], \;
\big[\tfrac{1}{2} (1 + u t) \big], \;
\big[\tfrac{1}{2} (1 + v t) \big], \;
\big[\tfrac{1}{2} (1 + u v t) \big], \;
S \big)
\]
is a basis for $K_0 (C^* (\ZZ^2 \rtimes \ZZ_2)) \cong \ZZ^6.$

\medskip
\noindent
{\textbf{(b) The case $F = \ZZ_3.$}}
Here we have the presentation
\[
\ZZ^2 \rtimes \ZZ_3
   = \big\lk u, v, t \colon t^3 = 1, \; u v = v u, \;
         t u t^{-1} = u^{-1} v, \; t v t^{-1} = u^{-1} \big\rk.
\]
The conjugacy classes of the maximal finite subgroups are represented by
\[
M_0 = \lk t \rk, \quad M_1 = \lk u t \rk, \quad M_2 = \lk u^2 t \rk,
\]
which are all isomorphic to $\ZZ_3.$
Set $\zt = \tfrac{1}{2} \big( - 1 + i \sqrt{3} \big),$
and define
\[
p_0 = \tfrac{1}{3} (1 + t + t^2),
\quad
q_0 = \tfrac{1}{3} \big( 1 + u t + (u t)^2 \big),
\quad
r_0 = \tfrac{1}{3} \big( 1 + u^2 t + (u^2 t)^2 \big),
\]
\[
p_1 = \tfrac{1}{3} \big( 1 + \zeta t + (\zeta t)^2 \big),
\quad
q_1 = \tfrac{1}{3} \big( 1 + \zeta u t + (\zeta u t)^2 \big),
\quad
r_1 = \tfrac{1}{3} \big( 1 + \zeta u^2 t + (\zeta u^2 t)^2 \big).
\]
Then, if $S \in K_0 (C^* (\ZZ^2 \rtimes \ZZ_3))$
maps to a generator of $\widetilde{K}_0 (C^* (\ZZ^2)) \cong \ZZ,$
we obtain from Theorem~\ref{thm-Wolfgang} that
\[
( [1], \; [p_0],
\; [p_1], \; [q_0], \; [q_1], \; [r_0], \; [r_1], \; S )
\]
is a basis for $K_0 (C^* (\ZZ^2 \rtimes \ZZ_3)) \cong \ZZ^8.$

\medskip
\noindent
{\textbf{(c) The case $F = \ZZ_4.$}}
In this case we have the presentation
\[
\ZZ^2 \rtimes \ZZ_4
= \big\lk u, v, t \colon t^4 = 1, \; u v = v u, \; t u t^{-1} = v,
   \; t v t^{-1} = u^{-1} \big\rk.
\]
The maximal finite subgroups are given up to conjugacy by
\[
M_0 = \lk t \rk, \quad M_1 = \lk u t \rk, \quad M_2 = \lk u t^2 \rk.
\]
The groups $M_0$ and $M_1$ are isomorphic to $\ZZ_4,$
and $M_2 \cong \ZZ_2.$
Define
\[
p_0 = \tfrac{1}{4} \big( 1 + t + t^2 + t^3 \big),
\quad
p_1 = \tfrac{1}{4} \big( 1 + i t - t^2 - i t^3 \big),
\quad
p_2 = \tfrac{1}{4} \big( 1 - t + t^2 - t^3 \big),
\]
\[
q_0 = \tfrac{1}{4} \big( 1 + u t + (u t)^2 + (u t)^3 \big),
\quad
q_1 = \tfrac{1}{4} \big( 1 + i u t - (u t)^2 -i (u t)^3 \big),
\]
\[
q_2 = \tfrac{1}{4} \big( 1 - u t + (u t)^2 - (ut)^3 \big),
\quad {\text{and}} \quad
r = \tfrac{1}{2} (1 + u t^2).
\]
Then, if $S \in K_0 (C^* (\ZZ^2 \rtimes \ZZ_4))$
maps to a generator of $\widetilde{K}_0 (C^* (\ZZ^2)) \cong \ZZ,$
we obtain from Theorem~\ref{thm-Wolfgang} the basis
\[
\big( [1], \; [p_0], \; [p_1], \; [p_2], \; [q_0], \; [q_1],
\; [q_2], \; [r], \; S \big)
\]
of $K_0 (C^* (\ZZ^2 \rtimes \ZZ_4) ) \cong \ZZ^9.$

\medskip
\noindent
{\textbf{(d) The case $F = \ZZ_6.$}}
Here we have the presentation
\[
\ZZ^2 \rtimes \ZZ_6
   = \big\lk u, v, t \colon t^6 = 1, \; u v = v u, \; t u t^{-1} = v, \;
        t v t^{-1} = u^{-1} v \big\rk.
\]
The maximal finite subgroups are given up to conjugacy by
\[
M_0 = \lk t \rk, \quad M_1 = \lk u t^2 \rk, \quad M_2 = \lk u t^3 \rk,
\]
with $M_0 \cong \ZZ_6,$ $M_1 \cong \ZZ_3,$ and $M_2 \cong \ZZ_2.$
Set $\zt = \tfrac{1}{2} \big( 1 + i \sqrt{3} \big),$
define
\[
p_j = \tfrac{1}{6} \big( 1 + (\zeta^j t) + (\zeta^j t)^2 + (\zeta^j t)^3
     + (\zeta^j t)^4 + (\zeta^j t)^5 \big)
\]
for $j = 0, 1, 2, 3, 4,$
and define
\[
q_0 = \tfrac{1}{3} \big( 1 + u t^2 + (u t^2)^2 \big),
\quad
q_1 = \tfrac{1}{3} \big( 1 + \zeta^2 u t^2 + (\zeta^2 u t^2)^2 \big),
\quad {\text{and}} \quad
r = \tfrac{1}{2} (1 + u t^3).
\]
Then, if $S \in K_0 (C^* (\ZZ^2 \rtimes \ZZ_6))$
maps to a generator of $\widetilde{K}_0 (C^* (\ZZ^2)) \cong \ZZ,$
we obtain from Theorem~\ref{thm-Wolfgang} that
\[
\big( [1], \; [p_0], \; [p_1], \; [p_2], \; [p_3], \; [p_4], \; [q_0],
\; [q_1], \; [r], \;S  \big)
\]
is a basis of $K_0 (C^* (\ZZ^2 \rtimes \ZZ_6)) \cong \ZZ^{10}.$
\end{ex}

Using the construction of the Fourier modules
due to Sam Walters (see~\cite{W1, W2}) we
later explicitly construct candidates for $S$
in all four cases considered above.

As a direct consequence of Example~\ref{ex-K0}
and Corollary~\ref{cor-crossed-theta} we now get:

\begin{cor}\label{cor-K0}
Let $F \subseteq \SL_2 (\ZZ)$
be one of the groups $\ZZ_2, \ZZ_3, \ZZ_4, \ZZ_6$ with generators
as described above
and let $F$ act on the rotation algebra $A_{\theta}$
for $\theta \in [0, 1]$
via the restriction of the action of $\SL_2 (\ZZ)$
as described in (\ref{eq-action}) and~(\ref{eq-gen}).
Then \mbox{$K_1 (A_{\theta} \rtimes_{\af} F) = 0$}
for all $F$ and $\theta,$ and for all $\theta \in [0, 1]$ we have

\[
K_0 (A_{\theta} \rtimes_{\af} \ZZ_2) \cong \ZZ^6, \quad
K_0 (A_{\theta} \rtimes_{\af} \ZZ_3) \cong \ZZ^8,
\]
\[
K_0 (A_{\theta} \rtimes_{\af} \ZZ_4) \cong \ZZ^9,
\quad {\text{and}} \quad
K_0 (A_{\theta} \rtimes_{\af} \ZZ_6) \cong \ZZ^{10}.
\]
\end{cor}

\section{A basis for $K_0 (A_{\theta} \rtimes_{\af} F)$}\label{sec-gen}

Let $\theta \in [0, 1].$
In this section we present an explicit
basis for $K_0 (A_{\theta} \rtimes_{\af} F)$
for $F = \ZZ_2, \ZZ_3, \ZZ_4, \ZZ_6.$

\begin{ntn}\label{P_exntn}
To simplify notation we write from now on
$e (x) = e^{2 \pi i x}$ for all $x \in \RR.$
\end{ntn}

Recall from Corollary~\ref{cor-crossed-theta} that
$A_{\theta} \rtimes_{\af} F
\cong C^* (\ZZ^2 \rtimes F, \, \widetilde{\om}_{\theta}).$

\begin{lem}\label{P_GenRel}
Let $\theta \in [0, 1],$
and let $F$ be one of $\ZZ_2, \ZZ_3, \ZZ_4, \ZZ_6.$
Then $A_{\theta} \rtimes_{\af} F$ is the universal \ca\  generated by
three unitaries
\[
u_{\theta}, \quad v_{\theta}, \quad t_{\te},
\]
subject to the relations
$v_{\theta} u_{\theta} = e (\theta) u_{\theta} v_{\theta}$ and
\[
\begin{array}{llll}
{\text{for $F = \ZZ_2$:}} \quad{\ }
    & t_{\te}^2 = t_{\te}, \, \, %
    & t_{\te} u_{\theta} t_{\te}^{-1} = u_{\theta}^{-1},
\quad\quad \quad\quad{\ }
    & t_{\te} v_{\theta} t_{\te}^{-1} = v_{\theta}^{-1} \\
{\text{for $F = \ZZ_3$:}}
    & t_{\te}^3 = t_{\te}, \quad
    & t_{\te} u_{\theta} t_{\te}^{-1}
       = e \big( - \tfrac{1}{2} \theta \big) u_{\theta}^{-1} v_{\theta},
    & t_{\te} v_{\theta} t_{\te}^{-1} = u_{\theta}^{-1} \\
{\text{for $F = \ZZ_4$:}}
    & t_{\te}^4 = t_{\te}, \quad
    & t_{\te} u_{\theta} t_{\te}^{-1} = v_{\theta},
    & t_{\te} v_{\theta} t_{\te}^{-1} = u_{\theta}^{-1} \\
{\text{for $F = \ZZ_6$:}}
    & t_{\te}^6 = t_{\te}, \quad
    & t_{\te} u_{\theta} t_{\te}^{-1} = v_{\theta},
    & t_{\te} v_{\theta} t_{\te}^{-1}
        = e \big(- \tfrac{1}{2} \theta \big) u_{\theta}^{-1} v_{\theta}.
\end{array}
\]
The \ca\  generated by $u_{\te}$ and $v_{\te}$ is $A_{\te}.$
Moreover, the \hm\  given on the generators by
\[
u_0 \mapsto u, \quad v_0 \mapsto v
\quad {\text{and}} \quad t_0 \mapsto t
\]
is an isomorphism
$A_0 \rtimes_{\af} F \to C^* (\ZZ^2 \rtimes F),$
and the \hm\  given on the generators by
\[
u_1 \mapsto -u, \quad v_1 \mapsto -v
\quad {\text{and}} \quad t_1 \mapsto t
\]
is an isomorphism
$A_1 \rtimes_{\af} F \to C^* (\ZZ^2 \rtimes F).$
\end{lem}

\begin{proof}
The first part follows from the description of the action of
$\SL_2 (\ZZ)$ on $A_{\theta}$ as given in
Corollary~\ref{cor-crossed-theta} and the realization of the groups
$F = \ZZ_2, \ZZ_3, \ZZ_4, \ZZ_6 \subseteq \SL_2 (\ZZ)$
with generators as in~(\ref{eq-gen}).
The last part is trivial for $\te = 0,$
and follows from $e \big(- \frac{1}{2} \theta \big) = -1$ for $\te = 1.$
\end{proof}

Although the cocycle $\widetilde{\om}_1$
is in the class of the trivial cocycle, it is not
trivial itself,
so it has to be taken into account in the computations below!

We now fix some $a \in (0, 1].$
Consider the $2$-cocycle
\[
\widetilde{\Om}_a
\colon (\ZZ^2 \rtimes F) \times (\ZZ^2 \rtimes F)
          \to C \big([a, 1], \, \TT \big)
\]
given by
$\widetilde{\Om}_a \big(\, \cdot\, ,\,  \cdot, \, \theta)
    = \widetilde{\om}_{\theta} ( \, \cdot \, ,  \, \cdot \, )$
for $\theta \in [a, 1].$
It is clear that it restricts at $\te$ to the cocycle
$\om_{\theta}$ on $\ZZ^2.$
For convenience, we write
$C^* \big( \ZZ^2 \rtimes F, \, \widetilde{\Om}_a \big)$
for the twisted crossed product
$C \big([a, 1] \big) \rtimes_{\widetilde{\Om}_a} (\ZZ^2 \rtimes F)$
and $C^* (\ZZ^2, \Om_a)$ for the twisted crossed
product $C \big([a, 1] \big) \rtimes_{{\Om}_a} \ZZ^2.$
As in Remark~\ref{rem-crossed-theta},
we may write
\begin{equation}\label{eq-crossed}
C^* \big( \ZZ^2 \rtimes F, \, \widetilde{\Om}_a \big)
            \cong C^* (\ZZ^2, \Om_a) \rtimes F
\end{equation}
where the action of $F$ on $C^* (\ZZ^2, \Om_a)$
is given fiberwise by the actions on the fibers
$A_{\theta} = C^* (\ZZ^2, \om_{\theta})$ for $\theta \in [a, 1].$

We want to construct suitable elements of
$K_0 \big( C^* \big( \ZZ^2 \rtimes F, \, \widetilde{\Om}_a \big) \big)$
which,
with respect to the identification
$C^* (\ZZ^2 \rtimes F, \, {\widetilde{\om}}_1)
                       \cong C^* (\ZZ^2 \rtimes F)$
of Lemma~\ref{P_GenRel},
project to classes in the fiber at~$1$
as given in Example~\ref{ex-K0}.

\begin{ntn}\label{P_UFam}
Let $a \in (0, 1].$
Since $G = \ZZ^2 \rtimes F$ is discrete,
there is a canonical map from $\ZZ^2 \rtimes F$ into the
group of unitaries of
$C^* \big( \ZZ^2 \rtimes F, \, \widetilde{\Om}_a \big)$
which sends a group element to the corresponding
Dirac function.
Let
$U_a, V_a, T_a
    \in C^* \big( \ZZ^2 \rtimes F, \, \widetilde{\Om}_a \big)$
be the images of the generators $u, v, t \in \ZZ^2 \rtimes F$
under this map.
We use the same notation for these unitaries regarded as elements
of suitable dense subalgebras,
and, for $U_a$ and $V_a,$ with $C^*\big( \ZZ^2, {\Om}_a \big)$
in place of $C^*\big( \ZZ^2 \rtimes F, \, \widetilde{\Om}_a \big).$
\end{ntn}

These unitaries clearly project to the generators
$u_{\theta}, v_{\theta}, t_{\te}
\in C^* (\ZZ^2 \rtimes F, \, \widetilde{\om}_{\theta})$
for $\theta \in [a, 1].$
Moreover, every function
$\varphi \in C \big( [a, 1] \big)$ can be considered canonically
as a central element of
$C^* \big( \ZZ^2 \rtimes F, \, \widetilde{\Om}_a \big)$
by identifying $\varphi$
with
$\delta_{1} \otimes \varphi \in
    C_c \big( \ZZ^2 \rtimes F, \, C ([a, 1]) \big)
    \subseteq C^* \big( \ZZ^2 \rtimes F, \, \widetilde{\Om}_a \big)$
(where $1$ denotes the
identity of $\ZZ^2 \rtimes F$).
In the lemma below, we identify the generators
$u, v, t$ of $C^* (\ZZ^2 \rtimes F)$
with the generators $-u_1, \, -v_1, \, t_1$ of
$C^* (\ZZ^2 \rtimes F, \, \widetilde{\om}_1)$ as
in the last part of Lemma~\ref{P_GenRel}.

\begin{lem}\label{lem-unitaries}
Recalling Notation~\ref{P_exntn},
define functions
$\varphi_1, \ph_2, \varphi_3 \in C \big([a, 1] \big)$ by
\[
\varphi_1 (\theta) = -e \big(\tfrac{1}{2} \theta \big),
\quad \varphi_2 (\theta) = e \big(\tfrac{1}{6} (2 + \theta) \big),
\quad {\text{and}}
\quad \varphi_3 (\theta) = i e \big(\tfrac{1}{4}\theta \big).
\]
Then:
\begin{enumerate}
\item
If $F = \ZZ_2,$ the elements
$T_a,$ $-U_a T_a,$ $-V_a T_a,$ and $\varphi_1 U_a V_a T_a$
are unitaries of order two which project onto the elements
$t,$ $u t,$ $v t,$ and $u v t$
in the fiber at~$1.$
\item
If $F = \ZZ_3,$ the elements
$T_a,$ $\varphi_2 U_a T_a,$ and $U_a^2 T_a$
are unitaries of order three which project onto the elements
$t,$ $u t,$ and $u^2 t$
in the fiber at~$1.$
\item
If $F = \ZZ_4,$ the elements $T_a,$ and $\varphi_3 U_a T_a$ are
unitaries of order four and $-U_a T_a^2$ is a unitary of order two
which project
onto the elements $t,$ $u t,$ and $u t^2$ in the fiber at~$1.$
\item
If $F = \ZZ_6,$ then $T_a$ is a unitary of order six,
$\varphi_2 U_a T_a^2$ is a unitary
of order three, and $-U_a T_a^3$ is a unitary of order two,
and they project onto the
elements $t,$ $u t^2,$ and $u t^3$ in the fiber at~$1.$
\end{enumerate}
\end{lem}

\begin{proof}
The proof is a straight-forward (but tedious) computation using the
relations for the unitaries $u_{\theta}, v_{\theta}, t_{\te}$
in each fiber $C^* (\ZZ^2 \rtimes F, \, \widetilde{\om}_{\theta}).$
We leave the details for the reader.
\end{proof}

As a consequence of Lemma~\ref{lem-unitaries},
we can now use exactly the same formulas
as for the explicit projections in $C^* (\ZZ^2 \rtimes F)$
presented in Example~\ref{ex-K0},
replacing the combinations of $u, v, t$ in the example
by appropriate combinations of $U_a, V_a, T_a,$
to obtain projections in
$C^* \big( \ZZ^2 \rtimes F, \, \widetilde{\Om}_a \big)$
which project onto
the ones in
$C^* (\ZZ^2 \rtimes F)
     \cong C^* (\ZZ^2 \rtimes F, \, \widetilde{\om}_1)$
via evaluation at~$1.$
Therefore, all we have to do is
to find a class
$S_a \in K_0 \big( C^* \big( \ZZ^2 \rtimes F,
                   \, \widetilde{\Om}_a \big) \big)$
which projects onto
a class $S \in K_0 (C^* (\ZZ^2 \rtimes F))$
via evaluation at $1$ with the properties
as in Example~\ref{ex-K0}.
Before we construct this element, it is convenient
to recall some basic properties
about equivariant $K$-theory of \ca s
with respect to finite group actions.

\begin{prop}\label{prop-GJ}
Suppose that $F$ is a finite group and $\alpha \colon F \to \Aut (A)$
is an action of $F$ on the \ca\  $A.$
Suppose further that
$\E$ is a finitely generated projective (right) $A$-module together
with a right action
$W \colon F \to \Aut (\E),$ written $(\xi, g) \mapsto \xi W_g,$
such that
$(\xi W_g) a = (\xi \alpha_g (a)) W_g$
for all $\xi \in \E,$ $a \in A,$ and $g \in F.$
Then $\E$ becomes a finitely generated projective
$A \rtimes_{\af} F$-module with action defined by
\[
\xi \cdot \bigg( \sum_{g \in F} a_g \dt_g \bigg)
    = \sum_{g \in F} (\xi a_g) W_g.
\]
Moreover, the restriction to $A$ of this module structure
is just the original $A$-module $\E,$ with
the action of $F$ forgotten.
\end{prop}

\begin{proof}
This is the definition of the map in the Green-Julg theorem
(Theorem 11.7.1 of~\cite{Bl} or Theorem 2.6.1 of~\cite{Ph0}).
The last statement is obvious.
\end{proof}

We will use Proposition~\ref{prop-GJ}
and the isomorphism~(\ref{eq-crossed}) to construct
finitely generated projective
$C^* \big( \ZZ^2 \rtimes F, \, \widetilde{\Om}_a \big)$-modules
which represent a suitable class $S_a.$
To do this, we exploit below some basic constructions
and ideas of  Alain Connes and Marc Rieffel in~\cite{Co, Rief-canc}
and Sam Walters in~\cite{W1, W2}.
Define a new cocycle
$\Om_{1 / a} \in Z^2 \big(\ZZ^2, \, C \big([a, 1], \, \TT \big) \big)$
by setting
$\Om_{1 / a} ( \, \cdot \, ,  \, \cdot \, ) (\theta) = \om_{1 / \theta}$
for all $\theta \in [a, 1].$
Set $A = C^* (\ZZ^2, \Om_a)$ and
$B = C^* \big( \ZZ^2, \, \Om_{1 / a} \big).$
Then the fiber $B_{\theta}$ of $B$ at $\theta \in [a, 1]$
is the rotation algebra
$A_{1 / \theta} = C^* (\ZZ^2, \, \om_{1 / \theta}).$
Consider the space $\cS [a, 1]$
consisting of all complex functions on $\RR \times [a, 1]$
which are smooth and rapidly decreasing in the first variable
and continuous in the second variable
in each derivative
of the first variable.
Denote by $\cS (\ZZ^2, {\Om}_a)$ the set of rapidly decreasing
$C \big([a, 1] \big)$-valued functions on $\ZZ^2,$
viewed as a (dense) subalgebra of $C^* (\ZZ^2, {\Om}_a),$
and denote by $\cS (\ZZ^2, \Om_{1 / a})$ the same set,
viewed as a (dense) subalgebra
of $C^* (\ZZ^2, \, \Om_{1 / a}).$
With $U_a$ and $V_a$ as in Notation~\ref{P_UFam},
we may represent each element $\varphi \in \cS (\ZZ^2, {\Om}_a)$
as an infinite sum
\begin{equation}\label{eq-sum}
\varphi = \sum_{n, m \in \ZZ} \beta_{n, m} \, U_a^n V_a^m
\end{equation}
with
$\beta_{n, m} \in C \big([a, 1] \big) \subseteq Z (C^* (\ZZ^2, \Om_a))$
and
$\left(\begin{smallmatrix} n \\ m \end{smallmatrix} \right)
     \to \| \beta_{n, m} \|_{\infty}$ rapidly decreasing on $\ZZ^2.$
Similarly,
we may write each element
$\psi \in \cS (\ZZ^2, \Om_{1 / a})$ as a sum
\begin{equation}\label{eq-sum1}
\psi = \sum_{n, m \in \ZZ}
      \alpha_{n, m} \, U_{1 / a}^n V_{1 / a}^m
\end{equation}
with $\alpha_{n, m} \in C \big([a, 1] \big).$
We define a right action of $\cS (\ZZ^2, {\Om}_a)$
on $\cS[a, 1]$ via the actions of
$U_a,$ $V_a,$ and $f \in C \big([a, 1] \big)$ given by
\begin{equation}\label{eq-action1}
(\xi \cdot U_a) (s, \theta)
   = \xi (s + \theta, \, \theta), \quad
(\xi \cdot V_a) (s, \theta)
   = e (s) \xi (s, \theta), \quad
{\text{and}} \quad
(\xi \cdot f) (s, \theta) = f (\theta) \xi (s, \theta).
\end{equation}
Moreover, we define an $\cS (\ZZ^2, \Om_a)$-valued inner product on
$\cS[a, 1]$ as follows.
For $n, m \in \ZZ,$ set
\begin{equation}\label{eq-act2}
\lk \xi, \eta \rk^a_{n, m} (\theta)
    = \theta \int_{- \infty}^{\infty}
      {\overline{\xi (x + n \theta, \, \theta)}}
      \eta (x, \theta) e (- m x) \, d x.
\end{equation}
Then set
\begin{equation}\label{eq-action2}
\lk \xi, \eta \rk_{\cS (\ZZ^2, {\Om}_a)}
     = \sum_{n, m \in \ZZ}
          \lk \xi, \eta \rk^a_{n, m} U_a^n V_a^m.
\end{equation}
There is also a left action
of $\cS (\ZZ^2, \Om_{1 / a})$ on $\cS[a, 1]$ given by
\begin{equation}\label{eq-action3}
U_{1 / a} \cdot \xi (s, \theta) = \xi (s + 1, \, \theta),
\quad
V_{1 / a} \cdot \xi (s, \theta)
    = e (- s / \theta) \xi (s, \theta),
\quad
\alpha \cdot \xi (s, \theta) = \alpha (\theta) \xi (s, \theta),
\end{equation}
for $\alpha \in C \big( [a, 1] \big),$ and a left inner product
\begin{equation}\label{eq-action4}
{_{\cS (\ZZ^2, \, {\Om}_{1 / a})} \lk \xi, \eta \rk}
     = \sum_{n, m \in \ZZ}
          \lk \xi, \eta \rk^{1 / a}_{n, m}
                 U_{1 / a}^n V_{1 / a}^m
\end{equation}
with
\[
\lk \xi, \eta \rk^{1 / a}_{n, m} (\theta)
    = \int_{- \infty}^{\infty}
           \xi (x - n, \, \theta)
           {\overline{\eta (x, \theta) }} e ( m x / \theta) \, d x.
\]

Evaluated at each fixed $\theta \in [a, 1],$
these actions and inner products determine
$B_{\theta} - A_{\theta}$ imprimitivity
bimodules $\E_{\theta} = \overline{\cS (\RR)}$
as in~\cite{Rief-canc} (Theorem~1.1 and the proof of Lemma~1.5).
The formulas given here are
from the beginning of Section~3 of~\cite{W1}.
(The paper~\cite{Rief-canc} uses a dense subalgebra of
$A_{\theta}$ consisting of functions on $\ZZ \times \TT.$)
It follows then easily from the fact
that all operations are $C \big( [a, 1] \big)$-linear
that $\cS[a, 1]$ completes to give a
$B - A$ imprimitivity bimodule $\E.$
Since both $B$ and $A$ are unital, it follows from
the arguments given after Notation~1.3 of~\cite{Rief-canc}
that $\E$ is a finitely generated projective
$A$-module with respect to the given action of $A$ on $\E.$
In order to construct a module over
$C^* \big( \ZZ^2 \rtimes F, \, \widetilde{\Om}_a \big) = A \rtimes F,$
it thus suffices to construct a suitable action of $F$ on $\E.$
The operator $W_t$ in the following proposition is
really the action of the unitary $T_a$ on the module.

\begin{prop}\label{prop-E}
Let $F$ be one of the groups $\ZZ_2, \ZZ_3, \ZZ_4, \ZZ_6$
with fixed generator $t.$
We define an action of $F$
on the dense subspace $\cS[a, 1]$ (described above) of the
$B - A$ imprimitivity bimodule $\E$ via the
action of the generator $t$ given by:
\begin{align*}
\begin{array}{ll}
(\xi W_t) (s, \theta) = \xi (- s, \, \theta)
& {\text{ (for $F = \ZZ_2$)}}  \\
(\xi W_t) (s, \theta)
   = e^{- \pi i / 12} \theta^{-1/2} e (  s^2 / (2 \theta))
     \int_{- \infty}^{\infty} \xi (x, \theta) e (s x / \theta) \, d x
& {\text{ (for $F = \ZZ_3$)}} \rule{0em}{2.7ex} \\
(\xi W_t) (s, \theta)
   = \theta^{-1/2} \int_{- \infty}^{\infty} \xi (x, \theta)
             e (s x / \theta) \, d x
& {\text{ (for $F = \ZZ_4$)}} \rule{0em}{2.7ex} \\
(\xi W_t) (s, \theta)
   = e^{ \pi i / 12} \theta^{-1/2}
     \int_{- \infty}^{\infty}
        \xi (x, \theta) e \big( (2 s x - x^2) / (2 \theta) \big) \, d x
          \quad\quad
& {\text{ (for $F = \ZZ_6$)}}. \rule{0em}{2.7ex}
\end{array}
\end{align*}
Then these actions extend to actions on $\E$
such that $\E$ becomes an $F$-equivariant
finitely generated projective $A$-module
as in Proposition~\ref{prop-GJ}.
\end{prop}

\begin{proof}
The operator which appears in the case $F = \ZZ_6$ is,
in the fiber over $\theta \in [a, 1],$
the hexic transform of Walters~\cite{W2}
with respect to the parameter $\mu = \frac{1}{2 \theta}.$
It follows from
Theorem~1 of~\cite{W2} that it has period six
and that its square is the operator which appears in the
case $F = \ZZ_3,$ which therefore has period three.
The operator which appears in case $F = \ZZ_4$ is,
in the fiber over $\theta,$
the one introduced by Walters
in Section~3 of~\cite{W1}, where it is shown that it has period four
with square equal to the operator
which appears in case $F = \ZZ_2.$
We have to show that, in all four cases, there exists
a unique extension of $W_t$ to $\E.$
This will be clear as soon as we verify the equation
\[
\lk \xi W_t, \, \eta W_t \rk_{A}
        = \alpha_{t^{-1}} (\lk \xi, \eta \rk_{A})
\]
for all $\xi, \eta \in \cS[a, 1]$
and in all cases $F = \ZZ_2, \ZZ_3, \ZZ_4, \ZZ_6.$
The cases $F = \ZZ_2$ and $F = \ZZ_3$
will follow directly from the cases
$F = \ZZ_4$ and $F = \ZZ_6.$
As a sample, we do the details for the case
$F = \ZZ_6$ here and leave the computation for the case $F = \ZZ_4$
to the reader.
(Note that these
relations in the case $F = \ZZ_4$
have been used implicitly in Section~3 of~\cite{W1}.)

Replacing $\xi$ by $\xi W_t^{-1}$
and applying $\alpha_t$ on both sides,
we reduce to checking the equation
\begin{equation}\label{eq-eta}
\lk \xi W_t^{-1}, \, \eta \rk_{A}
= \alpha_t (\lk \xi, \, \eta W_t \rk_{A}).
\end{equation}
The formula for the transformation $\xi \mapsto \xi W_t^{-1}$ is given
in~\cite{W2} (see the remark immediately after Theorem 1~of~\cite{W2},
and continue to take $\mu = \frac{1}{2 \theta}$),
and is
\begin{equation}\label{eq-Uinverse}
(\xi W_t^{-1}) (s, \theta)
= \frac{e^{- \pi i / 12}}{\sqrt{\theta}}
   e \big( \tfrac{1}{2} \te^{-1} s^2 \big)
     \int_{- \infty}^{\infty} \xi (x, \theta)
                     e (- s x / \theta) \, d x.
\end{equation}
Let $\varphi \in \cS (\ZZ^2, \Om_a)$ be given
$\ph = \sum_{n, m \in \ZZ} \beta_{n, m} \, U_a^n V_a^m$
as in~(\ref{eq-sum}).
We use the definition of the action $\alpha_t$
(for example, see the commutation relations
in Lemma~\ref{P_GenRel} for the case $F = \ZZ_6$)
and, at the third step,
the commutation relation
$v_{\theta} u_{\theta}^{-1} = e (- \theta) u_{\theta}^{-1} v_{\theta},$
to get, for $\theta \in [a, 1],$
\begin{align*}%
\big(\alpha_t (\varphi) \big) (\theta)
   & = \sum_{n, m \in \ZZ}
           \beta_{n, m} (\theta) \alpha^{\theta}_t (u_{\theta})^n
           \alpha^{\theta}_t (v_{\theta})^m \\
   & = \sum_{n, m} \beta_{n, m} (\theta)
            e \big( - \tfrac{1}{2} m \theta \big)
            v_{\theta}^n \big( u_{\theta}^{-1}v_{\theta} \big)^m \\
   & = \sum_{n, m} \beta_{n, m} (\theta)
            e \big( - \tfrac{1}{2} \theta (m^2 + 2 n m) \big)
                      u_{\theta}^{-m}v_{\theta}^{n + m} \\
   & = \sum_{k, l} \beta_{l + k, \, -k} (\theta)
            e \big( \tfrac{1}{2} \theta (2 k (l + k) - k^2) \big)
            u_{\theta}^k v_{\theta}^l.
\end{align*}

Using this formula together with~(\ref{eq-action2}),
we see that we have to check the identity
\begin{equation}
\label{eq-ident}
\lk \xi W_t^{-1}, \eta \rk_{k, l}^{a} (\theta)
     = \lk \xi, \eta W_t \rk_{k + l, \, -k}^a (\theta)
          e \big(\tfrac{1}{2} \theta [2 k (k + l) - k^2] \big)
\end{equation}
for all $k, l \in \ZZ.$
Combining~(\ref{eq-act2}) with~(\ref{eq-Uinverse})
gives for the left hand side:
\begin{equation}\label{eq-identleft}
\begin{split}
& \lk \xi W_t^{-1}, \eta \rk_{k, l}^{a} (\theta)
\\
& = e^{ \pi i / 12} \sqrt{\theta}
     \int_{- \infty}^{\infty} \int_{- \infty}^{\infty}
       {\overline{\xi (r, \theta)}} \eta (x, \theta)
        e \big( \tfrac{1}{2} \theta^{-1} \big[ 2 (x + k \theta) r
           - (x + k \theta)^2 - 2 \theta l x \big] \big) \, d r \, d x.
\end{split}
\end{equation}
For the right hand side of~(\ref{eq-ident}), first define
\[
\begin{split}
& \psi (r, x, \theta, k, l)
\\
&  \hspace*{1em} {\mbox{}}
= e \left( \frac{1}{2 \theta}
    \big[ 2 \big( r - (l + k) \theta \big) x - x^2
         + 2 \theta k \big( r - (l + k) \theta \big)
               + \theta^2 \big( 2 k (k + l) - k^2 \big) \big] \right).
\end{split}
\]
Calculate, using
at the first step~(\ref{eq-act2})
and the formula for $W_t$ as given in the proposition,
and using at the second step the change of variables
$r$ to $r - (k + l) \theta,$
to get:
\begin{equation*}\label{eq-identright}
\begin{split}
& \lk \xi, \eta W_t \rk_{k + l, \, -k}^a (\theta)
          e \big(\tfrac{1}{2} \theta [2 k (k + l) - k^2] \big) \\
&  \hspace*{2em} {\mbox{}}
   = e^{ \pi i / 12}
    \sqrt{\theta} \int_{- \infty}^{\infty} \int_{- \infty}^{\infty}
      {\overline{\xi (r + (k + l) \theta, \, \theta)}} \eta (x, \theta)
              \\
& \hspace*{10em}
      e \big( \tfrac{1}{2} \theta^{-1} (2 r x - x^2 + 2 \theta k r)
                      \big)
      e \big( \tfrac{1}{2}\theta [2 k (l + k) - k^2] \big)
                 \, d r \, d x \\
&  \hspace*{2em} {\mbox{}}
= e^{ \pi i / 12} \sqrt{\theta} \int_{- \infty}^{\infty}
          \int_{- \infty}^{\infty}
             {\overline{\xi (r, \theta)}} \eta (x, \theta)
                   \psi (r, x, \theta, k, l) \, d r \, d x.
\end{split}
\end{equation*}
So the result follows from the identity
\begin{align*}
2 (x + k \theta ) r
    & - (x + k \theta)^2 - 2 \theta l x \\
    & = 2 \big( r - (l + k) \theta \big) x - x^2
         + 2 \theta k \big( r - (l + k) \theta \big)
               + \theta^2 \big( 2 k (k + l) - k^2 \big).
\end{align*}

We still have to check the equation
\[
(\xi W_t) \varphi = (\xi \alpha_t (\varphi)) W_t
\]
for all $\xi \in \cS[a, 1]$ and $\varphi \in \cS (\ZZ^2, \Om_a),$
which will then imply that $\E$ is $F$-equivariant
in the sense of Proposition~\ref{prop-GJ}.
Of course it is enough to check this identity on the generating
set $C \big([a, 1] \big) \cup \{U_a, V_a \}.$
Since $F$ acts trivially on $C \big([a, 1] \big),$
this amounts to verifying the identities
\[
(\xi W_t) U_a = (\xi \alpha_t (U_a)) W_t,
\quad
(\xi W_t) V_a = (\xi \alpha_t (V_a)) W_t,
\quad {\text{and}} \quad
(\xi W_t) f = (\xi f) W_t,
\]
for all $f \in C \big([a, 1] \big).$
The last equation is easy,
so we only have to consider the first two equations.
For $F = \ZZ_4,$ these equations
were used (but not explicitly checked) in Section~3 of~\cite{W1}.
The case
$F = \ZZ_2$ follows directly from the case $F = \ZZ_4,$
since the operator in the quadratic
case is the square of the operator in the quartic case.
We do the computations for the case
$F = \ZZ_6$ below.
The case $F = \ZZ_3$ will then follow from this in the same way as
the case $\ZZ_2$ follows from the case $\ZZ_4.$

The computations in the case that $t$
is the generator of $\ZZ_6$ are basically the same as
the ones given by Walters in Section~3 of~\cite{W2}.
By analogy with Lemma~\ref{P_GenRel}, and with
$g \in C \big( [a, 1] \big)$ defined by $g (\te) = e (- \te / 2),$
we have
\[
\alpha_t (U_a) = V_a
\quad {\text{and}} \quad
\alpha_t (V_a) = g U_a^{-1} V_a.
\]
Using this, we have to check the identities
\[
(\xi W_t) U_a = (\xi V_a) W_t
\andeqn
(\xi W_t) V_a = (\xi g U_a^{-1} V_a) W_t.
\]
The first of these identities is very easy,
and we omit the proof.
For the second identity, we compute
\begin{equation}\label{eq-1}
\big( (\xi W_t) V_a \big) (s, \theta)
= e (s) (\xi W_t) (s, \theta)
= e (s) \frac{e^{ \pi i / 12}}{\sqrt{\theta}}
      \int_{- \infty}^{\infty} \xi (x, \theta)
           e \big( \tfrac{1}{2} \theta^{-1} (2 s x - x^2) \big) \, d x.
\end{equation}
On the other side,
changing the variable $x$ to $x + \te$ at the last step, we have
\begin{align*}
\big( (\xi g
& U_a^{-1} V_a) W_t \big) (s, \theta)
   = e \big(- \tfrac{1}{2} \theta \big)
       \frac{e^{ \pi i / 12}}{\sqrt{\theta}}
       \int_{- \infty}^{\infty} \big( \xi U_a^{-1} V_a \big) (x, \theta)
          e \big( \tfrac{1}{2} \theta^{-1} (2 s x - x^2) \big) \, d x
                     \\
& = e \big(- \tfrac{1}{2} \theta \big)
       \frac{e^{ \pi i / 12}}{\sqrt{\theta}}
       \int_{- \infty}^{\infty} e (x) \xi (x - \theta, \, \theta)
           e \big( \tfrac{1}{2} \theta^{-1} (2 s x - x^2) \big) \, d x
                     \\
& = e \big(- \tfrac{1}{2} \theta \big)
       \frac{e^{ \pi i / 12}}{\sqrt{\theta}}
       \int_{- \infty}^{\infty} e (x + \theta) \xi (x, \theta)
            e \big( \tfrac{1}{2} \theta^{-1}
                   [2 s (x + \theta) - (x + \theta)^2] \big)
          \, d x.
\end{align*}
The last expression is equal to $\big( (\xi W_t) V_a \big) (s, \theta)$
by~(\ref{eq-1}) and the easy to verify identity
\[
e \big(- \tfrac{1}{2} \theta \big) e (x + \theta)
      e \big( \tfrac{1}{2} \theta^{-1}
                   [2 s (x + \theta) - (x + \theta)^2] \big)
= e (s) e \big( \tfrac{1}{2} \theta^{-1} (2 s x - x^2) \big),
\]
which holds for all $s, x$ and $\theta.$
\end{proof}

\begin{ntn}\label{P_ModNtn}
Let $F$ be one of $\ZZ_2, \ZZ_3, \ZZ_4, \ZZ_6,$
and let $\E$ be the $F$ equivariant $A$-module
of Proposition~\ref{prop-E}.
We denote by $\E^F$
the finitely generated
$A \rtimes F$ module constructed from $\E$
by the procedure described in Proposition~\ref{prop-GJ}.
We write $\E_{\theta}^F$
and $\E_{\theta}$ for the evaluations of
$\E^F$ and $\E$ at $\theta \in [a, 1].$
\end{ntn}

Thus, with $\ev_{\te}$ denoting evaluation at $\te$
(on both $A = C^* (\ZZ^2, \Om_a)$ and $A \rtimes F$),
in $K$-theory we have
\[
(\ev_{\te})_* ( [ \E ] ) = [ \E_{\te} ] \in K_0 (A_{\te})
\andeqn
(\ev_{\te})_* \big( \big[ \E^F \big] \big)
   = \big[ \E_{\te}^F \big]  \in K_0 (A_{\te} \rtimes_{\af} F).
\]

We now look at the special
value $\theta = 1,$ for which we have
$A_1 \rtimes_{\af} F \cong C^* (\ZZ^2 \rtimes F).$

\begin{lem}\label{lem-gen}
Let $F$ be one of the groups $\ZZ_2, \ZZ_3, \ZZ_4, \ZZ_6$ and let
$[\E_1^F] \in K_0 (C^* (\ZZ^2 \rtimes F))$
be as in Notation~\ref{P_ModNtn} for $\te = 1.$
Then $\big[ \E_1^F \big]$ is mapped onto a generator
of the group $\widetilde{K}_0 (C^* (\ZZ^2)) \cong \ZZ$
under the projection
\[
\pr \circ \res^{\ZZ^2}_{\ZZ^2 \rtimes F}
\colon K_0 (C^* (\ZZ^2 \rtimes F)) \to \widetilde{K}_0 (C^* (\ZZ^2))
\]
which appears in the sequence of Theorem~\ref{thm-Wolfgang}.
\end{lem}

\begin{proof}
Since $\res^{\ZZ^2}_{\ZZ^2 \rtimes F} ([\E_1^F]) = [\E_1],$
it suffices to show that
$( [1], [\E_1] )$ is a basis for $K_0 (C^* (\ZZ^2)).$
By Remark~\ref{rem-crossed-theta}
(applied to the trivial group $H = \{ 1 \}$)
it suffices to check that there exists one $\theta \in [a, 1]$
such that $( [1], [\E_{\theta}] )$ is a basis
for $K_0 (A_{\theta}).$
Choose $\theta\in [a,1] \smallsetminus \QQ$
and let $\ta \colon A_{\theta} \to \RR$
be the canonical \tst\  on $A_{\te}.$
In Theorem~1.4 of~\cite{Rief-canc} it is shown
that $\tau_* ([\E_{\theta}]) = \theta.$
It follows from the Appendix in~\cite{PV}
that $\ta_*$ is an isomorphism from $K_0 (A_{\theta})$
to $\ZZ + \te \ZZ,$
so the result follows.
\end{proof}

We are now ready to give an explicit basis
for $K_0 (A_{\theta} \rtimes F)$
for all $\theta \in (0, 1]$
and for $F = \ZZ_2, \ZZ_3, \ZZ_4, \ZZ_6.$

\begin{thm}\label{thm-generators}
Let $\te \in (0, 1].$
Let $u_{\theta},$ $v_{\theta,}$ and $t_{\theta}$ be the
generators of
$A_{\theta} \rtimes_{\af} F
= C^* (\ZZ^2 \rtimes F, \, \widetilde{\om}_{\theta})$
of Lemma~\ref{P_GenRel}.
Choose $a$ with $0 < a \leq \te,$ and
let $\big[ \E^F_{\theta} \big] \in K_0 (A_{\theta} \rtimes_{\af} F)$
be as in Notation~\ref{P_ModNtn}.
For $n = 2, 3, 4, 6$
let $\ta_{n, \te}$ denote the canonical \tst\  on
$A_{\theta} \rtimes_{\af} \ZZ_n.$
Then:
\begin{enumerate}
\item\label{thm-generators-Z_2}
If $F = \ZZ_2,$ define
\[
p^{\theta}
   = \tfrac{1}{2} (1 + t_{\theta}), \quad \quad
r^{\theta}
   = \tfrac{1}{2} \big( 1 - e \big( \tfrac{1}{2} \theta \big)
          u_{\theta} v_{\theta} t_{\theta} \big),
\]
\[
q_0^{\theta}
   = \tfrac{1}{2} (1 - u_{\theta} t_{\theta}),
\quad \andeqn \quad
q_1^{\theta}
   = \tfrac{1}{2} (1 - v_{\theta} t_{\theta}).
\]
Then the classes
\[
[1], \; \big[ p^{\theta} \big], \; \big[ q_0^{\theta} \big],
\; \big[ q_1^{\theta} \big],
\; \big[ r^{\theta} \big], \; \big[ \E_{\theta}^{\ZZ_2} \big]
\]
form a basis for $K_0 (A_{\theta} \rtimes_{\af} \ZZ_2) \cong \ZZ^6.$
Moreover,
$(\tau_{2, \theta})_*$
takes the following values (in this order)
on these classes:
\[
1, \; \frac{1}{2}, \; \frac{1}{2}, \; \frac{1}{2},
   \;  \frac{1}{2}, \; \frac{\theta}{2}.
\]
\item\label{thm-generators-Z_3}
If $F = \ZZ_3,$ set $\zeta = \tfrac{1}{2} \big( -1 + i \sqrt{3} \big),$
and define
\[
p_0^{\theta}
   = \tfrac{1}{3} \big( 1 + t_{\theta} + t_{\theta}^2 \big), \quad \quad
p_1^{\theta}
   = \tfrac{1}{3}
     \big( 1 + \zeta t_{\theta} + (\zeta t_{\theta})^2 \big),
\]
\[
q_0^{\theta}
   = \tfrac{1}{3} \big( 1 + e \big(\tfrac{1}{6} [2 + \theta] \big)
                                 u_{\theta} t_{\theta}
     + e \big( \tfrac{1}{3} [2 + \theta] \big)
                         (u_{\theta} t_{\theta})^2 \big),
\]
\[
q_1^{\theta}
   = \tfrac{1}{3} \big( 1
      + e \big( \tfrac{1}{6} [2 + \theta] \big)
                         \zeta u_{\theta} t_{\theta}
      + e \big( \tfrac{1}{3} [2 + \theta] \big)
                     ( \zeta u_{\theta} t_{\theta} )^2 \big),
\]
\[
r_0^{\theta}
     = \tfrac{1}{3} \big( 1 + u_{\theta}^2 t_{\theta}
        + (u_{\theta}^2 t_{\theta})^2  \big), \quad \andeqn \quad
r_1^{\theta}
    = \tfrac{1}{3} \big( 1 + \zeta u_{\theta}^2 t_{\theta}
       + (\zeta u_{\theta}^2 t_{\theta})^2 \big).
\]
Then the classes
\[
[1], \; \big[ p_0^{\theta} \big], \; \big[ p_1^{\theta} \big],
\; \big[ q_0^{\theta} \big],
\; \big[ q_1^{\theta} \big],
\; \big[ r_0^{\theta} \big], \; \big[ r_1^{\theta} \big],
\; \big[ \E_{\theta}^{\ZZ_3} \big]
\]
form a basis for $K_0 (A_{\theta} \rtimes_{\af} \ZZ_3) \cong \ZZ^8.$
Moreover,
$(\tau_{3, \theta})_*$
takes the following values
on these classes:
\[
1, \; \frac{1}{3}, \; \frac{1}{3}, \; \frac{1}{3},
\; \frac{1}{3}, \; \frac{1}{3}, \; \frac{1}{3},
   \; \frac{\theta}{3}.
\]
\item\label{thm-generators-Z_4}
If $F = \ZZ^4,$ define
\[
p_0^{\theta}
     = \tfrac{1}{4}
          \big( 1 + t_{\theta} + t_{\theta}^2 + t_{\theta}^3 \big),
\quad \quad
p_1^{\theta}
    = \tfrac{1}{4}
      \big( 1 + i t_{\theta} - t_{\theta}^2 -i t_{\theta}^3 \big),
\]
\[
p_2^{\theta}
    = \tfrac{1}{4}
      \big( 1 - t_{\theta} + t_{\theta}^2 - t_{\theta}^3 \big),
\quad \quad
r^{\theta} = \tfrac{1}{2}  \big( 1 - u_{\theta} t_{\theta}^2 \big),
\]
\[
q_0^{\theta}
    = \tfrac{1}{4}
         \big( 1 + i e \big( \tfrac{1}{4} \theta \big)
                  u_{\theta} t_{\theta}
        - \big[ e \big( \tfrac{1}{4} \theta \big)
                 u_{\theta} t_{\theta} \big]^2
        - i \big[ e \big( \tfrac{1}{4} \theta \big)
                 u_{\theta} t_{\theta} \big]^3 \big),
\]
\[
q_1^{\theta} = \tfrac{1}{4}
     \big( 1 - e \big( \tfrac{1}{4} \theta \big) u_{\theta} t_{\theta}
     + \big[e \big( \tfrac{1}{4} \theta \big)
              u_{\theta} t_{\theta} \big]^2
     - \big[e \big( \tfrac{1}{4} \theta \big)
              u_{\theta} t_{\theta}\big]^3 \big),
\]
and
\[
q_2^{\theta}
      = \tfrac{1}{4} (1 - i e \big( \tfrac{1}{4} \theta \big)
             u_{\theta} t_{\theta}
       - \big[e \big( \tfrac{1}{4} \theta \big) u_{\theta} t_{\theta})^2
        + i \big[e \big( \tfrac{1}{4} \theta \big)
                  u_{\theta} t_{\theta} \big]^3 \big).
\]
Then the classes
\[
[1], \; \big[ p_0^{\theta} \big], \; \big[ p_1^{\theta} \big],
\; \big[ p_2^{\theta} \big],
\;  \big[ q_0^{\theta} \big], \;  \big[ q_1^{\theta} \big],
\; \big[ q_2^{\theta} \big],
\; \big[ r^{\theta} \big], \; \big[ \E_{\theta}^{\ZZ_4} \big]
\]
form a basis for $K_0 ( A_{\theta} \rtimes_{\af} \ZZ_4) \cong \ZZ^9.$
Moreover,
$(\tau_{4, \theta})_*$
takes the following values
on these classes:
\[
1, \; \frac{1}{4}, \; \frac{1}{4}, \; \frac{1}{4},
\; \frac{1}{4}, \; \frac{1}{4}, \; \frac{1}{4}, \; \frac{1}{2}, \;
\frac{\theta}{4}.
\]
\item\label{thm-generators-Z_6}
If $F = \ZZ_6,$
set $\zeta = \tfrac{1}{2} \big( 1 + i \sqrt{3} \big),$
define
\[
p_j^{\theta}
    = \tfrac{1}{6}
       \big( 1 + (\zeta^j t_{\theta}) + (\zeta^j t_{\theta})^2
            + (\zeta^j t_{\theta})^3 + (\zeta^j t_{\theta})^4
            + (\zeta^j t_{\theta})^5 \big)
\]
for $0 \leq j \leq 4,$ and define
\[
q_0^{\theta} = \tfrac{1}{3}
\big( 1 + e \big( \tfrac{1}{6} (2 + \theta) \big)
                             u_{\theta} t_{\theta}^2
     +  \big[ e \big( \tfrac{1}{6} (2 + \theta) \big)
                    u_{\theta} t_{\theta}^2 \big]^2 \big),
\]
\[
q_1^{\te} = \tfrac{1}{3}
   \big( 1 + \zeta^2 e \big( \tfrac{1}{6} (2 + \theta) \big)
                            u_{\theta} t_{\theta}^2
    + \big[ \zeta^2 e \big( \tfrac{1}{6} (2 + \theta) \big)
      u_{\theta} t_{\theta}^2 \big]^2 \big),
\]
and
\[
r^{\theta} = \tfrac{1}{2} \big( 1 - u_{\theta} t_{\theta}^3 \big).
\]
Then the classes
\[
[1], \; \big[ p_0^{\theta} \big], \; \big[ p_1^{\theta} \big],
\; \big[ p_2^{\theta} \big],
\; \big[ p_3^{\theta} \big], \; \big[ p_4^{\theta} \big],
\; \big[ q_0^{\theta} \big],
\; \big[ q_1^{\theta} \big], \; \big[ r^{\theta} \big],
\; \big[ \E_{\theta}^{\ZZ_6} \big]
\]
form a basis for $K_0 ( A_{\theta} \rtimes_{\af} \ZZ_6) \cong \ZZ^{10}.$
Moreover,
$(\tau_{6, \theta})_*$
takes the following values
on these classes:
\[
1, \; \frac{1}{6}, \; \frac{1}{6}, \; \frac{1}{6}, \; \frac{1}{6},
\; \frac{1}{6}, \; \frac{1}{3}, \; \frac{1}{3}, \; \frac{1}{2},
\; \frac{\theta}{6}.
\]
\end{enumerate}
\end{thm}

\begin{proof} %
In each case, we define \pj s
in $C^* \big( \ZZ^2 \rtimes_{\af} F, \, \widetilde{\Om}_a \big)$
as follows.
In the formulas in Example~\ref{ex-K0}
(excluding $S,$ which we treat below), substitute for
$t,$ $u t,$ etc.\  the unitaries given in Lemma~\ref{lem-unitaries}
whose images in the fiber at~$1$ are $t,$ $u t,$ etc.
Thus, for example, for $F = \ZZ_6$ we replace $u t^2$ by
$\ph_2 U_a T_a^2,$ and in place of the \pj\  $q_0$ we obtain
$Q_0 = \tfrac{1}{3}
         \big( 1 + \ph_2 U_a T_a^2 + ( \ph_2 U_a T_a^2 )^2 \big).$
In each case, by Lemma~\ref{lem-unitaries},
the evaluation at $\te \in [0, 1]$ of the \pj\  obtained this way
is the corresponding \pj\  in the statement of the theorem.
For example, for $F = \ZZ_6$ we get $\ev_{\te} (Q_0) = q_0^{\te}.$
To get our replacement for $S,$
choose a \pj\  in some matrix algebra
over $C^* \big( \ZZ^2 \rtimes F, \, \widetilde{\Om}_a \big)$
corresponding to the module $\E^F$
as in Notation~\ref{P_ModNtn}.
Thus, for each choice of $F,$ we have a collection ${\mathcal{P}}$ of
\pj s in $C^* \big( \ZZ^2 \rtimes F, \, \widetilde{\Om}_a \big)$
or a matrix algebra over this algebra.

Evaluate at $\te = 1.$
The classes of the \pj s $\ev_1 (p)$ for $p \in {\mathcal{P}}$
are exactly
the classes listed in Example~\ref{ex-K0},
except with $S$ replaced by $\big[ \E^F_{\theta} \big].$
By Example~\ref{ex-K0} and Lemma~\ref{lem-gen},
these classes form a basis
for $K_0 (C^* (\ZZ^2 \rtimes F)).$
By Remark~\ref{rem-crossed-theta},
for any $\te \in [a, 1],$
the classes of the \pj s $\ev_{\te} (p)$ for $p \in {\mathcal{P}}$
form a basis for $K_0 (A_{\theta} \rtimes_{\af} F).$

The computation of the trace is clear
from the description for all classes
except the classes $\big[ \E_{\theta}^F \big].$
The values of the traces of the classes
$\big[ \E_{\theta}^F \big]$ are computed in Proposition~3.3 of~\cite{W1}
for $F = \ZZ_4$ and
in Theorem~2 of~\cite{W2} for $F = \ZZ_3, \ZZ_6.$
The case $F = \ZZ_2$ is done similarly.
\end{proof}

\section{The tracial Rokhlin property}\label{Sec:TRP}

In this section, we prove that for $\te \in \R \smallsetminus \Q$
the actions of $\Z_2,$ $\Z_3,$ $\Z_4,$ and $\Z_6$ on $A_{\te},$
given by~(\ref{eq-action}) for the groups
generated as in~(\ref{eq-gen}),
have the tracial Rokhlin property,
and that for every nondegenerate skew-symmetric matrix $\Te,$
the flip action of $\Z_2$ on $A_{\Te}$
has the tracial Rokhlin property.
We do this by showing that,
under suitable conditions,
the \tRp\  is equivalent to outerness
of the corresponding action on the
type II$_1$ factor obtained as the weak operator closure of
$A$ in the Gelfand-Naimark-Segal representation associated with
the \tst.
A preliminary result, stated entirely in terms of \ca s,
gives a criterion for the \tRp\  in terms of trace norms.
For use elsewhere,
we state and prove this result for \ca s with an arbitrary
simplex of \tst s.

We will make extensive use of the $L^2$-norm (or seminorm)
associated with a \tst\  $\ta$ of a \ca\  $A,$ given by
$\| a \|_{2, \ta} = \ta (a^* a)^{1/2}.$
See the discussion before Lemma V.2.20 of~\cite{Tk}
for more on this seminorm in the von Neumann algebra context.
All the properties we need are immediate from its identification
with the seminorm in which one completes $A$ to obtain the
Hilbert space $H_{\ta}$ for the
Gelfand-Naimark-Segal representation associated with $\ta,$
and from the relation $\ta (b a) = \ta (a b).$
In particular, we always have
$\| a b c \|_{2, \ta}
    \leq \| a \| \cdot \| b \|_{2, \ta} \cdot \| c \|.$

Since the \ca s in this section will almost all be finite,
we recall the \tRp\  for actions of finite groups
on finite infinite dimensional simple unital \ca s.
The following result is Lemma~1.16 of~\cite{PhtRp1a}.

\begin{prop}\label{TRPCond}
Let $A$ be a finite infinite dimensional \suca,
and let $\af \colon G \to \Aut (A)$
be an action of a finite group $G$ on $A.$
Then $\af$ has the
tracial Rokhlin property \ifo\  for every finite set
$S \subseteq A,$ every $\ep > 0,$
and every nonzero positive element $x \in A,$
there are \mops\  $e_g \in A$ for $g \in G$ such that:
\begin{itemize}
\item[(1)]
$\| \af_g (e_h) - e_{g h} \| < \ep$ for all $g, h \in G.$
\item[(2)]
$\| e_g a - a e_g \| < \ep$ for all $g \in G$ and all $a \in S.$
\item[(3)]
With $e = \sum_{g \in G} e_g,$ the \pj\  $1 - e$ is \mvnt\  to a
\pj\  in the \hsa\  of $A$ generated by $x.$
\end{itemize}
\end{prop}

We begin with a simple reformulation of the \tRp\  which is
valid in the presence of good comparison properties for \pj s.
Recall that, if $A$ is a unital \ca,
then we say that the
order on \pj s over $A$ is determined by traces
if Blackadar's Second Fundamental Comparability Question
(1.3.1 in~\cite{Bl3}) holds for all matrix algebras over $A.$
That is, whenever $n \in \N$ and $p, q \in M_n (A)$ are \pj s such that
$\ta (p) < \ta (q)$ for all \tst s $\ta$ on $A,$ then
$p \precsim q.$

\begin{lem}\label{DfnUsingTrace}
Let $A$ be a finite infinite dimensional
simple separable unital \ca\   with Property~(SP)
(every nonzero hereditary subalgebra contains a nonzero \pj)
and such that the order on \pj s over $A$ is determined by traces.
Let $\af \colon G \to \Aut (A)$
be an action of a finite group $G$ on $A.$
Then $\af$ has the \tRp\  \ifo\  %
for every finite set $S \subseteq A$ and every $\ep > 0,$ there exist
orthogonal \pj s $e_g \in A$ for $g \in G$ such that:
\begin{itemize}
\item[(1)]
$\| \af_g (e_h) - e_{g h} \| < \ep$ for all $g, h \in G.$
\item[(2)]
$\| e_g a - a e_g \| < \ep$ for all $g \in G$ and all $a \in S.$
\item[(3)]
With $e = \sum_{g \in G} e_g,$ we have $\ta (1 - e) < \ep$
for all $\ta \in T (A).$
\end{itemize}
\end{lem}

\begin{proof}
Suppose that the condition of Proposition~\ref{TRPCond} holds,
and let $S$ and $\ep$ be given.
Choose $n \in \N$ with $n > 1 / \ep.$
By Lemma~1.10 of~\cite{PhtRp1a}, %
there exist $n$ nonzero \mops\  in $A.$
Let $x$ be one of them.
Then $\ta (x) < \ep$ for all $\ta \in T (A).$
Apply Proposition~\ref{TRPCond} with this $x$
and with $S$ and $\ep$ as given.
Conversely, assume the condition of the lemma,
and let $S,$ $\ep,$ and $x$ be given.
Using Property~(SP), choose a nonzero \pj\  $q \in \ov{x A x},$
and apply the condition of the lemma with $\ep$ replaced by
$\min \left( \ep, \, \inf_{\ta \in T (A)} \ta (q) \right).$
The assumption that the order on \pj s over $A$ is determined by traces
implies that $1 - e \precsim q,$
giving~(3) of Proposition~\ref{TRPCond}.
\end{proof}

\begin{thm}\label{TraceVersion}
Let $A$ be an infinite dimensional
simple separable unital \ca\  with tracial rank zero.
Let $\af \colon G \to \Aut (A)$
be an action of a finite group $G$ on $A.$
Then $\af$ has the \tRp\  \ifo\  %
for every finite set $S \subseteq A$ and every $\ep > 0,$ there exist
orthogonal \pj s $e_g \in A$ for $g \in G$ such that:
\begin{itemize}
\item[(1)]
$\| \af_g (e_h) - e_{g h} \|_{2, \ta} < \ep$
for all $g, h \in G$ and all $\ta \in T (A).$
\item[(2)]
$\| [e_g, a] \|_{2, \ta} < \ep$
for all $g \in G,$ all $a \in S,$ and all $\ta \in T (A).$
\item[(3)]
$\sum_{g \in G} e_g = 1.$
\end{itemize}
The equivalence is also valid if one substitutes for
Condition~(3) the following condition:
\begin{itemize}
\item[(3$'$)]
$\ta \left( 1 - \sum_{g \in G} e_g \right) < \ep$
for all $\ta \in T (A).$
\end{itemize}
\end{thm}

\begin{proof}
We first prove that if $\af$ has the \tRp,
then the version of the condition of the theorem using
Conditions~(1), (2), and~(3$'$) holds.
Let $S \subseteq A$ be finite and let $\ep > 0.$
Since $A$ has tracial rank zero,
Corollary~5.7 and Theorems~5.8 and~6.8 of~\cite{LnTTR}
imply that the order on \pj s over $A$ is determined by traces.
Therefore we may apply Lemma~\ref{DfnUsingTrace} with $S$ and $\ep$
as given.
Since $\| a \|_{2, \ta} \leq \| a \|$ for all $\ta \in T (A),$
it is easy to show that the resulting \pj s
satisfy~(1), (2), and~(3$'$).

Next, assume that the version with Conditions~(1), (2), and~(3$'$)
holds.
Let $S \subseteq A$ be finite and let $\ep > 0.$
We may assume that $\| a \| \leq 1$ for all $a \in S.$
Apply the hypothesis with $S$ as given and with
$\min \left( \ts{ \frac{1}{3}} \ep, \ts{ \frac{1}{9}} \ep^2 \right)$
in place of $\ep.$
Call the resulting \pj s $p_g$ for $g \in G.$
Set $p = \sum_{g \in G} p_g,$ set $e_1 = p_1 + 1 - p,$
and set $e_g = p_g$ for $g \in G \smallsetminus \{ 1 \}.$
Since
$\| e_1 - p_1 \|_{2, \ta} = \ta (1 - p)^{1/2} < \ts{ \frac{1}{3}} \ep$
for all $\ta \in T (A),$ it is easy to show that
the $e_g$ satisfy~(1), (2), and~(3).

Now we prove that if~(1), (2), and~(3) hold,
then $\af$ has the \tRp.
We follow the proof of Theorem~2.14 of~\cite{OP2}
(but note that the finite set there is called $F$).
We describe the choices and constructions carefully,
but omit many of the details in the verification of the estimates.
We verify the condition of Lemma~\ref{DfnUsingTrace}.
So let $S \subseteq A$ be a finite set, and let $\ep > 0.$
\Wolog\  $\ep < 1.$

Choose $\ep_0 > 0$ with
\[
\ep_0 < \frac{\ep}{3 \card (G)}
\]
and so small that whenever
$p_g,$ for $g \in G,$ are \pj s in a \ca\  $B$ which satisfy
$\| p_g p_h \| < 4 \ep_0$ for $g \neq h,$
then there are orthogonal \pj s $e_g \in B,$ for $g \in G,$
such that $e_1 = p_1$ and $\| e_g - p_g \| < \ts{ \frac{1}{3}} \ep$
for $g \in G.$

Apply Lemma~2.8 of~\cite{OP2}
with $\ep_0$ in place of $\ep$ and with $n = \card (G) - 1,$
and let $\ep_1 > 0$ be the resulting value of $\dt.$
We also require $\ep_1 \leq \ep_0.$
Then set
\[
\ep_2 = \min \left(1, \, \frac{\ep_1^2}{16}, \, \frac{\ep_1}{8},
       \, \frac{\ep}{18} \right).
\]
Set $T = \bigcup_{g \in G} \af_g (S).$
Apply Lemma~2.13 of~\cite{OP2} with $\ep_2$ in place of $\ep,$
with $T$ in place of $F,$
and with $G$ in place of $S.$
We obtain \pj s $q, q_0 \in A,$
unital \fd\  subalgebras $E \subseteq q A q$
and $E_0 \subseteq q_0 A q_0,$
and automorphisms $\ph_g \in \Aut (A)$
for $g \in G,$
such that:
\begin{itemize}
\item[(1)]
$\ph_1 = \id_A$ and
$\| \ph_g - \af_g \| < \ep_2$ for all $g \in G.$
\item[(2)]
For every  $g \in G$ and $x \in E,$
we have $q_0 \ph_g (x) = \ph_g (x) q_0$
and $q_0 \ph_g (x) q_0 \in E_0.$
\item[(3)]
For every $a \in T,$ we have $\| q a - a q \| < \ep_2$ and
$\dist (q a q, \, E) < \ep_2.$
\item[(4)]
$\ta (1 - q), \, \ta (1 - q_0) < \ep_2$ for all $\ta \in T (A).$
\end{itemize}

Apply Lemma~2.12 of~\cite{OP2} with $\ep_2$
in place of $\ep$ and $E_0 + \C (1 - q_0)$ in place of $E,$
obtaining $\dt > 0.$
Also require $\dt \leq \ep_2.$

Apply the hypothesis
with $\dt$ in place of $\ep,$
and with a system of matrix units for $E_0 + \C (1 - q_0)$
in place of $S,$
getting \pj s $p_g$ for $g \in G.$
Let $B_0 = A \cap [E_0 + \C (1 - q_0)]',$
the subalgebra of $A$ consisting of all elements which commute
with everything in $E_0 + \C (1 - q_0).$
Apply the choice of $\dt$ using Lemma~2.12 of~\cite{OP2} to $p_1,$
obtaining a \pj\  $f \in B_0$
which satisfies
$\| f - p_1 \|_{2, \ta} < \ep_2$ for all $\ta \in T (A).$
Since $q_0$ is in the center of $E_0 + \C (1 - q_0),$
the element $f_1 = q_0 f$ is also a \pj\  in $B_0.$
For $g \in G,$ since $q_0 \ph_{g^{-1}} (E) q_0 \subseteq E_0,$
it follows that $f_1$ commutes with all elements of
$q_0 \ph_{g^{-1}} (E) q_0$
and hence with all elements of $\ph_{g^{-1}} (E).$
Therefore $f_g = \ph_g (f_1)$ commutes with all elements of $E,$
including $q.$
So $f_g$ also commutes with $1 - q.$

Following the analogous estimates
in the proof of Theorem~2.14 of~\cite{OP2},
we get
$\| f_1 - p_1 \|_{2, \ta} < \ts{ \frac{1}{4}} \ep_1 + \ep_2$
for $\ta \in T (A)$ and
$\| f_g - \af_g (f_1) \| < \ep_2$ for $g \in G.$
We then get
(with a slight difference:
the sum in the estimate in~\cite{OP2} is not needed here),
\[
\| f_g - p_g \|_{2, \ta}
< \ep_2 + \left( \ts{ \frac{1}{4}} \ep_1 + \ep_2 \right) + \ep_2
\]
for $\ta \in T (A).$
Continuing as there, for $g \in G \smallsetminus \{ 1 \}$ we get
\[
\| f_1 f_g \|_{2, \ta}
   = \| f_1 f_g - p_1 p_g \|_{2, \ta}
   < \ts{ \frac{1}{2}} \ep_1 + 4 \ep_2
   \leq \ep_1.
\]

We saw that $f_g \in B = A \cap [E + \C (1 - q)]'$ for $g \in G.$
This algebra has real rank zero because $E + \C (1 - q)$ is \fd.
Therefore Lemma~2.8 of~\cite{OP2},
applied to $B$ with $\{ \ta |_B \colon \ta \in T (A) \}$
in place of $T,$
and the choice of $\ep_1,$ provide
a \pj\  $r \in A \cap [E + \C (1 - q)]'$
such that $r \leq f_1,$
such that $\| r f_g \| < \ep_0$ for $g \in G,$
and such that $\ta (r) > \ta (f_1) - \ep_0$
for all $\ta \in T (A).$

Now use Corollary~2.4 of~\cite{OP2} to find a \pj\  %
$e_1 \in A \cap [E + \C (1 - q)]'$ such that
$e_1 \leq q,$
such that $\| r e_1 - e_1 \| < \ep_0,$
and such that $[e_1] \geq [1] - ([1] - [q]) - ([1] - [r])$
in $K_0 (A \cap [E + \C (1 - q)]').$
The last inequality implies that
\[
\ta (e_1) \geq \ta (r) - \ta (1 - q)
   > \ta (f_1) - \ep_0 - \ep_2
   \geq \ta (f_1) - 2 \ep_0
\]
for all $\ta \in T (A).$

We now show that the $\af_g (e_1),$ for $g \in G,$
are approximately orthogonal.
It suffices to estimate $\| e_1 \af_g (e_1) \|$
for $g \in G.$
We first use $\ph_g (r) \leq \ph_g (f_1) = f_g$ to get
$\| r \ph_g (r) \| \leq \| r f_g \| < \ep_0.$
Then, following the proof of Theorem~2.14 of~\cite{OP2},
for $g \in G$ we get
$\| e_1 \ph_g (e_1) \| < 3 \ep_0$ and
\[
\| e_1 \af_g (e_1) \|
\leq \| \af_g - \ph_g \| + \| e_1 \ph_g (e_1) \|
< 4 \ep_0.
\]

Now use the choice of $\ep_0$ to find orthogonal \pj s
$e_g \in A \cap [E + \C (1 - q)]'$ for $g \in G \smallsetminus \{ 1 \},$
all orthogonal to $e_1,$ such that
$\| e_g - \af_g (e_1) \| < \ts{ \frac{1}{3}} \ep$
for $g \in G \smallsetminus \{ 1 \}.$
Then the $e_g,$ for $g \in G,$ satisfy
Conditions (1) through~(3) of Lemma~\ref{DfnUsingTrace},
by estimates essentially the same as in the last part of
the proof of Theorem~2.14 of~\cite{OP2}.
\end{proof}

\begin{ntn}\label{N:DComm}
If $A$ is a simple \ca\  with unique \tst\  $\ta,$
and $\af \in \Aut (A),$
we write $\af''$ for the
automorphism of $\pi_{\ta} (A)''$ determined by $\af.$
\end{ntn}

\begin{thm}\label{OuterImpTRP}
Let $A$ be a simple separable unital \ca\  with tracial rank zero,
and suppose that $A$ has a unique \tst\  $\ta.$
Let $\pi_{\ta} \colon A \to B (H_{\ta})$
be the Gelfand-Naimark-Segal representation associated with $\ta.$
Let $G$ be a finite group, and let $\af \colon G \to \Aut (A)$
be an action of $G$ on $A.$
Then $\af$ has the \tRp\  \ifo\  %
$\af_g''$ is an outer automorphism of
$\pi_{\ta} (A)''$ for every $g \in G \smallsetminus \{ 1 \}.$
\end{thm}

\begin{proof}
Assume that $\af_g''$ is outer for
every $g \in G \smallsetminus \{ 1 \}.$
We verify the hypotheses of Theorem~\ref{TraceVersion}.
Thus let $\ep > 0$ and let $S \subseteq A$ be a finite subset.
\Wolog\  $\| a \| \leq 1$ for all $a \in S.$
Set $n = \card (G),$
and set $\ep_0 = (4 n + 1)^{-1} \ep.$
Choose $\dt > 0$ as in Lemma~2.9 of~\cite{OP2}
with $n$ as given
and with $\ep_0$ in place of $\ep.$

We regard $A$ as a subalgebra of $N = \pi_{\ta} (A)'',$
and we let $\ta$ also denote the extension
of $\ta$ to $N.$
The algebra $N$ is hyperfinite by Lemma~2.16 of~\cite{OP2}.

We apply Lemma~5.2.1 of~\cite{Jn}.
The term ``equivariant s.m.u.'' [system of matrix units]
is defined in Section~1.5 of~\cite{Jn},
and $N (\ps)$ is defined in Section~1.2 of~\cite{Jn}.
(Also note that pages~44 and~45 of~\cite{Jn} are switched.)
We take the set $X$ of Lemma~5.2.1 of~\cite{Jn} to be $G$
with the left translation action.
The lemma provides, ignoring the off diagonal matrix units,
\pj s $p_g \in \pi_{\ta} (A)''$
such that $\af_g'' (p_h) = p_{g h}$ for $g, h \in G,$
and such that $\| [p_g, a] \|_{2, \ta} < \dt$
for $g \in G$ and $a \in S.$

For $g \in G$ use Lemma~2.15 of~\cite{OP2} to find a \pj\  %
$q_g \in A$ such that
$\| q_g - p_g \|_{2, \ta}
   < \min \left( \frac{1}{2} \dt, \, \ep_0 \right).$
Then $\| q_g q_h \|_{2, \ta} < \dt$ for $g \neq h,$
so the choice of $\dt$ using Lemma~2.9 of~\cite{OP2} provides
\mops\  $e_g \in A$ for $g \in G \smallsetminus \{ 1 \}$ such that
$\| e_g - q_g \|_{2, \ta} < \ep_0,$
and therefore $\| e_g - p_g \|_{2, \ta} < 2 \ep_0.$
Set $e_1 = 1 - \sum_{g \in G \smallsetminus \{ 1 \}} e_g.$
Since $\sum_{g \in G} p_g = 1,$ we get
\[
\| e_1 - p_1 \|_{2, \ta}
   \leq \sum_{g \in G \smallsetminus \{ 1 \} } \| e_g - p_g \|_{2, \ta}
   < 2 (n - 1) \ep_0.
\]
So $\| e_g - p_g \|_{2, \ta} < 2 n \ep_0$
for all $g \in G.$

For $g, h \in G,$ using $\af_g'' (p_h) = p_{g h}$
we now get
\[
\| \af_g (e_h) - e_{g h} \|_{2, \ta}
   \leq \| \af_g'' (e_h - p_h) \|_{2, \ta}
           + \| p_{g h} - e_{g h} \|_{2, \ta}
   < 2 n \ep_0 + 2 n \ep_0
   < \ep.
\]
Also, if $a \in S$ and $g \in G$ then
\[
\| [a, e_g] \|_{2, \ta}
   \leq 2 \| e_g - p_g \|_{2, \ta} \cdot \| a \|
             + \| [a, p_g] \|_{2, \ta}
    < 2 \cdot 2 n \ep_0 + \ep_0
    \leq \ep.
\]
This completes the proof that outerness in the trace
representation implies the \tRp.

The proof of the converse is similar to the proof
of Lemma~1.5 of~\cite{PhtRp1a}.
Set $n = \card (G)$ as before,
and set
\[
\ep = \min \left( \frac{1}{2}, \, \frac{1}{4 \sqrt{2 n}} \right).
\]
Let $g_0 \in G \smallsetminus \{ 1 \},$ and let $u \in \pi_{\ta} (A)''$
be a unitary; we show that $\af_{g_0}'' \neq \Ad (u).$
Let $\xi_{\ta}$ be the standard cyclic vector for $\pi_{\ta},$
so that $\| x \|_{2, \ta} = \| x \xi_{\ta} \|$
for all $x \in \pi_{\ta} (A)''.$
The Kaplansky Density Theorem implies that
the unit ball of $\pi_{\ta} (A)$ is strong operator dense in
the unit ball of $\pi_{\ta} (A)''.$
Therefore there is $a \in A$ with $\| a \| \leq 1$ such that
$\| u - a \|_{2, \ta} < \ep.$
Choose a nonzero \pj\  $f \in A$ with $\ta (f) < \ep,$
and choose \mops\  $e_g \in A$ for $g \in G$
such that $\| e_g a - a e_g \| < \ep$ for $g \in G,$
such that $\| \af_g (e_h) - e_{g h} \| < \ep$ for $g, h \in G,$
and such that $1 - \sum_{g \in G} e_g \precsim f.$
Since $\ta$ is the unique \tst,
it is $G$-invariant,
so that
\[
\ta (e_1) = {\textstyle{ \frac{1}{n} \sum_{g \in G} \ta (e_g)  }}
           \geq \tfrac{1}{n} (1 - \ta (f))
           \geq \tfrac{1}{2 n}.
\]

Now
\begin{align*}
\| e_1 - e_{g_0} \|_{2, \ta}
& \leq 2 \| u - a \|_{2, \ta} + \| a e_1 - e_1 a \|
           + \| u e_1 u^* - \af_{g_0} (e_1) \|
           + \| \af_{g_0} (e_1) - e_{g_0} \|     \\
& < 2 \ep + \ep + \| u e_1 u^* - \af_{g_0} (e_1) \| + \ep.
\end{align*}
On the other hand, using orthogonality,
\[
\| e_1 - e_{g_0} \|_{2, \ta}
   \geq \| e_1 \|_{2, \ta}
   \geq \frac{1}{\sqrt{2 n}}.
\]
Therefore
\[
\| u e_1 u^* - \af_{g_0} (e_1) \| > \frac{1}{\sqrt{2 n}} - 4 \ep
\geq 0.
\]
Thus $u e_1 u^* \neq \af_{g_0} (e_1),$ and $\af_{g_0}'' \neq \Ad (u).$
\end{proof}

In case the group is $\Z,$
Theorem~2.18 of~\cite{OP2} gives two more equivalent conditions
for an action $\af \colon G \to \Aut (A)$ on a simple unital
\ca\  with tracial rank zero and unique \tst\  to have the \tRp,
namely that the crossed product $A \rtimes_{\af} G$ have real rank zero,
and that it have a unique \tst.
For finite groups, the \tRp\  implies these conditions,
but neither condition implies the \tRp.

The \tRp\  implies that the crossed product has real rank zero
under these hypotheses,
because it implies that the crossed product has tracial rank zero
(Theorem~2.6 of~\cite{PhtRp1a}).
The converse is easily seen to be false,
by considering the trivial action of any finite group on
any simple unital \ca\  %
with tracial rank zero and a unique \tst.
Example~2.9 of~\cite{PhtRp1b} shows that
it does not help to require that the action be outer.

We now consider the condition that $A \rtimes_{\af} G$
have a unique \tst.
One can give a direct argument that this follows from the \tRp.
The line of reasoning we give here uses some machinery
and an additional useful lemma.
For this lemma, we do not require the algebra to be finite,
so the \tRp\  is as in Definition~1.2 of~\cite{PhtRp1a}.

\begin{lem}\label{L:SubGp}
Let $A$ be an infinite dimensional \suca,
and let $\af \colon G \to \Aut (A)$
be an action of a finite group $G$ on $A$ which has the \tRp.
Let $H$ be a subgroup of $G.$
Then $\af |_H$ has the \tRp.
\end{lem}

\begin{proof}
Set $m = \card (G) / \card (H).$
Choose a set $S$ of right coset representatives for $H$ in $G,$
so that $\card (S) = m$ and every element of $G$
can be written uniquely as $h s$ with $h \in H$ and $s \in S.$

Let $S \subseteq A$ be finite,
let $\ep > 0,$
and let $x \in A$ be a positive element with $\| x \| = 1.$
Set $\ep_0 = \ep / m.$
The \tRp\  for $\af$ provides
\mops\  $p_g \in A$ for $g \in G$ such that:
\begin{itemize}
\item[(1)]
$\| \af_g (p_h) - p_{g h} \| < \ep_0$ for all $g, h \in G.$
\item[(2)]
$\| p_g a - a p_g \| < \ep_0$ for all $g \in G$ and all $a \in S.$
\item[(3)]
With $p = \sum_{g \in G} p_g,$ the \pj\  $1 - p$ is \mvnt\  to a
\pj\  in the \hsa\  of $A$ generated by $x.$
\item[(4)]
With $p$ as in~(3), we have $\| p x p \| > 1 - \ep.$
\end{itemize}
For $h \in H$ define $e_h = \sum_{s \in S} p_{h s}.$
These are clearly \mops.
We verify the analogs of Conditions (1) through~(4).
For $h, k \in H,$ we have
\[
\| \af_k (e_h) - e_{k h} \|
   \leq \sum_{s \in S} \| \af_k (p_{h s}) - p_{k h s} \|
   < m \ep_0
   = \ep.
\]
For $h \in H$ and $a \in S,$ we have
\[
\| e_h a - a e_h \|
   \leq \sum_{s \in S} \| p_{h s} a - a p_{h s} \|
   < m \ep_0
   = \ep.
\]
The \pj\  $e = \sum_{h \in H} e_h$ is equal to $p,$
so that $1 - e = 1 - p$ is \mvnt\  to a
\pj\  in the \hsa\  of $A$ generated by $x.$
Finally, $\| e x e \| = \| p x p \| > 1 - \ep.$
\end{proof}

\begin{prop}\label{P:TRPAndTraces}
Let $A$ be an infinite dimensional \suca,
and let $\af \colon G \to \Aut (A)$
be an action of a finite group $G$ on $A$ which has the \tRp.
Let $r \colon T (A \rtimes_{\af} G) \to T (A)$ be the restriction map
between the \tst\  spaces.
Then $r$ is an affine homeomorphism from $T (A \rtimes_{\af} G)$
to the subspace $T (A)^{G}$ of $G$-invariant \tst s on $A.$
\end{prop}

\begin{proof}
It is easy to see that the image of $r$ is exactly $T (A)^{G},$
and that $r$ is \ct\  and affine.
It therefore suffices to show that $r$ is injective,
that is, that a \tst\  $\ta$ on $A \rtimes_{\af} G$ is
uniquely determined by its values on $A.$
We prove this by proving that if $a \in A,$
if $g \in G \smallsetminus \{ 1 \},$
and if $u_g \in A \rtimes_{\af} G$ is the corresponding unitary,
then $\ta (a u_g) = 0.$

Let $H$ be the subgroup of $G$ generated by $g,$
and let $\bt = \af |_H.$
Set $\sm = \ta |_{A \rtimes_{\bt} H}.$
Then $\bt$ has the \tRp\  by Lemma~\ref{L:SubGp}.
Therefore ${\widehat{\bt}}$ is tracially approximately representable,
by Theorem~3.12(1) of~\cite{PhtRp1a}.
Choose $\mu \in {\widehat{H}}$ such that $\mu (g) \neq 1.$
It follows from Proposition~6.1 of~\cite{PhtRp1a}
that $\sm \circ {\widehat{\bt}}_{\mu} = \sm.$
Since ${\widehat{\bt}}_{\mu} (a u_g) = {\overline{\mu (g)}} a u_g,$
this implies that $\sm (a u_g) = 0.$
So $\ta (a u_g) = 0.$
\end{proof}

In particular, if $A$ has a unique \tst\  %
and $\af \colon G \to \Aut (A)$ has the \tRp,
then $A \rtimes_{\af} G$ has a unique \tst.

Again, the converse is false.

\begin{ex}\label{E:UniqueTrace}
First take $A = M_2$ and $G = (\Z_2)^2.$
Let the two generators act by conjugation by the unitaries
$\left( \begin{smallmatrix} 1 & 0 \\ 0 & -1 \end{smallmatrix} \right)$
and
$\left( \begin{smallmatrix} 0 & 1 \\ 1 & 0 \end{smallmatrix} \right).$
Because the matrices commute up to a scalar,
one gets a well defined action $\af \colon G \to \Aut (M_2),$
such that $\af_g$ is inner for every $g \in G.$
It is well known (see, for example, Example~4.2.3 of~\cite{Ph0})
that $A \rtimes_{\af} G \cong M_4,$ which has a unique \tst.

This example is finite dimensional,
and the \tRp\  was only defined for infinite dimensional \ca s.
To get an infinite dimensional example,
tensor the above example with, say, the trivial action of $G$
on the $2^{\infty}$ UHF algebra.
\end{ex}

The following lemma is well known,
and is given in this form as Lemma~3.1 of~\cite{OP2}.

\begin{lem}\label{L2Sum}
Let $A$ be a separable unital \ca\  with a faithful \tst\  $\ta.$
Let $(y_i)_{i \in I}$ be a family of unitaries such that
$\ta (y_i^* y_j) = 0$ for $i \neq j$
and whose linear span is dense in $A.$
Let $\pi \colon A \to B (H)$ be the
Gelfand-Naimark-Segal representation associated with $\ta.$
We identify $A$ with its image in $\pi (A)''.$
Then every $a \in \pi (A)''$ has a unique representation
as $a = \sum_{i \in I} \ld_i y_i,$
with convergence in $\| \cdot \|_{2, \ta}$ and whose coefficients
satisfy $\sum_{i \in I} | \ld_i |^2 = \| a \|_{2, \ta}^2.$
If $a$ is unitary then $\sum_{i \in I} | \ld_i |^2 = 1.$
\end{lem}

Recall that a real skew symmetric $d \times d$ matrix
$\Te$ is
{\emph{nondegenerate}} if whenever $x \in \Z^d$ satisfies
$\exp (2 \pi i \langle x, \Te y \rangle ) = 1$ for all $y \in \Z^d,$
then $x = 0.$
This definition is essentially from Section~1.1 of~\cite{Sl}.
It is well known that
the noncommutative torus $A_{\Te}$
(defined after Corollary~\ref{cor-real})
is simple \ifo\  $\Te$
is nondegenerate; the main part is Theorem~3.7 of~\cite{Sl},
and the complete statement is Theorem~1.9 of~\cite{PhtRp2}.

\begin{lem}\label{OuterInTr}
Let $\Te$ be a nondegenerate real skew symmetric $d \times d$ matrix.
Let $1 \leq k \leq d,$ and
let $\af \in \Aut (A_{\Te})$ be an automorphism such that
$\af (u_k) = \rh u_1^{m_1} u_2^{m_2} \cdots u_d^{m_d}$
for some $\rh \in \C$ with $| \rh | = 1,$
and for some $m \in \Z^d$
not equal to the $k$th standard basis vector $\dt_k.$
Then $\af''$ (as in Notation~\ref{N:DComm}) is outer.
\end{lem}

\begin{proof}
\Wolog\  $k = 1.$

Let $\ta$ be the unique \tst\  on $A_{\Te},$
and identify $A_{\Te}$ with its image in $\pi_{\ta} (A_{\Te})''.$

Suppose $\af''$ is inner,
so $\af'' = \Ad (w)$ for some unitary $w \in \pi_{\ta} (A_{\Te})''.$
We apply Lemma~\ref{L2Sum} with $I = \Z^d$
and $y_n = u_1^{n_1} u_2^{n_2} \cdots u_d^{n_d}$ for $n \in \Z^d,$
and write
\[
w^* = \sum_{n \in \Z^d} \ld_{n} u_1^{n_1} u_2^{n_2} \cdots u_d^{n_d}
\]
with convergence in $\| \cdot \|_{2, \ta}.$
Using the commutation relations,
there are numbers $\zt_n$ for $n \in \Z^d,$ with $| \zt_n | = 1,$
such that
\[
\big( u_1^{n_1} u_2^{n_2} \cdots u_d^{n_d} \big)
      \big( u_1^{m_1} u_2^{m_2} \cdots u_d^{m_d} \big)
= \zt_n u_1^{n_1 + m_1} u_2^{n_2 + m_2} \cdots u_d^{n_d + m_d}.
\]
Then, also with convergence in $\| \cdot \|_{2, \ta},$
\begin{align*}
\sum_{n \in \Z^d}
        \ld_{n} u_1^{n_1 + 1} u_2^{n_2} \cdots u_d^{n_d}
& = u_1 w^*
   = w^* \af (u_1)   \\
& = \sum_{n \in \Z^d} \rh
      \zt_n \ld_n u_1^{n_1 + m_1} u_2^{n_2 + m_2} \cdots u_d^{n_d + m_d}
              \\
& = \sum_{n \in \Z^d} \rh
    \zt_{n - m + \dt_1} \ld_{n - m + \dt_1}
             u_1^{n_1 + 1} u_2^{n_2} \cdots u_d^{n_d}.
\end{align*}
We have $\ld_{l} \neq 0$ for some $l \in \Z^d.$
Since $\rh$ and the $\zt_n$ have absolute value~$1,$
uniqueness of the series representation implies
\[
| \ld_{l} |
= | \ld_{l - m + \dt_1} |
= | \ld_{l - 2 m + 2 \dt_1} |
= \cdots.
\]
Since $m \neq \dt_1,$
this contradicts $\sum_{n \in \Z^d} | \ld_{n} |^2 < \I.$
\end{proof}

\begin{cor}\label{C:TRP}
Let $F \subseteq \SL_2 (\ZZ)$
be one of the groups $\ZZ_2, \ZZ_3, \ZZ_4, \ZZ_6,$
and let $\af \colon F \to \Aut (A_{\theta})$
be given by~(\ref{eq-action}) for the groups
generated as in~(\ref{eq-gen}).
Suppose $\te \not\in \Q.$
Then $\af$ has the \tRp.
\end{cor}

\begin{proof}
In each case,
Lemma~\ref{OuterInTr} implies that $\af_g''$ is outer
for $g \in G \smallsetminus \{ 1 \}.$
Apply Theorem~\ref{OuterImpTRP}.
\end{proof}

Let $\Z_2$ act on $\Z^d$ via the flip $n \mapsto - n.$
It is immediately clear that for any
$d \times d$ skew symmetric real matrix $\Theta,$
the cocycle $\om_{\Theta}$
is invariant under this action.
It thus induces an action, also called the flip,
of $\Z_2$ on the higher dimensional noncommutative torus $A_{\Theta}.$
On the standard generators, this action sends
$u_k$ to $u_k^*$ for $1 \leq k \leq d.$

\begin{cor}\label{L:FlipTRP}
Let $\Te$ be a nondegenerate real $d \times d$ skew symmetric matrix.
Then the flip action on $A_{\Te}$ has the \tRp.
\end{cor}

\begin{proof}
Lemma~\ref{OuterInTr} implies that the flip automorphism is outer.
Theorem~3.5 of~\cite{PhtRp2} implies that
$A_{\Te}$ is a simple separable unital \ca\  with
tracial rank zero.
Apply Theorem~\ref{OuterImpTRP}.
\end{proof}

\section{The structure of the crossed products}\label{Sec:CP}

In this section, we prove that for $\te \in \R \smallsetminus \Q$
the crossed products
$A_{\te} \rtimes_{\af} \Z_2,$
$A_{\te} \rtimes_{\af} \Z_3,$
$A_{\te} \rtimes_{\af} \Z_4,$
and $A_{\te}\rtimes_{\af} \Z_6,$
by the action~(\ref{eq-action}) for the groups
generated as in~(\ref{eq-gen}),
are AF~algebras,
and we determine when two of them are isomorphic.
For every nondegenerate $\Te,$
we also prove that the crossed product $A_{\Te} \rtimes \Z_2$
by the flip
is an AF~algebra.
We further prove that all the fixed point algebras are~AF.
As mentioned in the introduction,
this was already known for $A_{\te}\rtimes_{\af} \Z_2$ for all
irrational $\te$~\cite{P_BK},
for $A_{\te}\rtimes_{\af} \Z_4$ for ``most'' irrational $\te$~\cite{W3},
and for $A_{\Te}\rtimes \Z_2$
for ``most'' skew symmetric real matrices $\Te$~\cite{Bc}.

We need to show that the crossed products
satisfy the Universal Coefficient Theorem,
as in Theorem~1.17 of~\cite{RSUCT}.
A much stronger result than the one formulated in the proposition
below has been obtained by Meyer and Nest in
Proposition~8.5 of~\cite{MN}.
Since the arguments used in the special case we consider
here are much easier, we include the short proof.

\begin{prop}\label{prop-bootstrap}
Let $G$ be an amenable
(or, more generally, an a-$T$-menable) group
which can be embedded as a closed
subgroup of some almost connected group $L.$
Then the full and reduced
crossed products $A \rtimes_{\alpha} G$
and $A \rtimes_{\alpha, r} G$ of any type~I algebra $A$ by $G$
are $KK$-equivalent
to commutative \ca s.
\end{prop}

\begin{proof}
Let $C$ denote the maximal compact subgroup of $L$
and let $T^* (L/C)$
denote the cotangent space of the quotient $L/C$
(which is a manifold).
Then it follows
from the work of Kasparov (see Theorem~2 of~\cite{Kas1};
Definitions 4.1 and~5.1 and Theorems~5.2 and~5.7 of~\cite{Kas2})
that there exist Dirac and dual-Dirac
elements
\[
D \in KK_0^L \big( C_0 (T^* (L/C)), \, \CC \big)
\andeqn
\beta \in KK_0^L \big( \CC, \, C_0 (T^* (L/C)) \big)
\]
such that
\[
D \otimes_{\CC} \beta = 1
\in KK_0^L \big( C_0 (T^* (L/C)), \, C_0 (T^* (L/C)) \big)
\]
and such that
$\gamma = \beta \otimes_{C_0 (T^* (L/C))}D \in KK_0^L (\CC, \CC)$
is the so-called $\gamma$-element of $L.$
The restriction $\gamma_G$ of $\gamma$
from $L$ to $G$ is the $\gamma$-element of $G.$
(For example, see Remark~6.4 of~\cite{CE1}.)
Since $G$ is a-$T$-menable, we know from the work
of Higson and Kasparov (Theorem~1.1 of~\cite{HK})
that $\gamma_G = 1 \in KK_0^G (\CC, \CC).$
But this implies that $\CC$ is $KK^G$-equivalent to
the commutative algebra $C_0 (T^* (L/C)),$
on which $G$ acts properly.
Now proceed as in the proof of Lemma~5.4 of~\cite{CEO2}.
\end{proof}

The algebras in the previous proposition thus
satisfy the Universal Coefficient Theorem.
In particular, we have:

\begin{cor}\label{C:UCT}
Let $G$ be an amenable group
which can be embedded as a closed
subgroup of some almost connected group $L.$
Let $\om \colon G \times G \to \TT$
be a Borel $2$-cocycle on~$G.$
Then $C^* (G, \om)$ satisfies  the Universal Coefficient Theorem
(Theorem~1.17 of~\cite{RSUCT}).
\end{cor}

\begin{proof}
As in the discussion at the beginning of Section~\ref{sec-twist},
there is action $\alpha \colon G\to \Aut(\K)$
such that $C^* (G, \om) \otimes \K \cong \K \rtimes_{\alpha}G.$
Thus Proposition~\ref{prop-bootstrap} implies that $C^* (G, \om)$ is
$KK$-equivalent to a commutative \ca.
It therefore
satisfies the Universal Coefficient Theorem,
by Proposition~7.1 of~\cite{RSUCT}.
\end{proof}

\begin{thm}\label{T:CrPrdIsAF}
Let $F \subseteq \SL_2 (\ZZ)$
be one of the groups $\ZZ_2, \ZZ_3, \ZZ_4, \ZZ_6,$
and let $\af \colon F \to A_{\theta}$
be given by~(\ref{eq-action}) for the groups
generated as in~(\ref{eq-gen}).
Suppose $\te \not\in \Q.$
Then $A_{\te} \rtimes_{\alpha} F$ is an AF~algebra.
\end{thm}

\begin{proof}
We know that $A_{\te}$ is a simple separable unital \ca\  with
tracial rank zero,
say by the Elliott-Evans Theorem~\cite{EE}.
(See Theorem~3.5 of~\cite{PhtRp2}.)
By Corollary~\ref{C:TRP},
we may apply Corollary~1.6 and Theorem~2.6 of~\cite{PhtRp1a}
to conclude that $A_{\te} \rtimes_{\alpha} F$
is a simple separable unital \ca\  with tracial rank zero.
Corollary~\ref{C:UCT} implies that $A_{\te} \rtimes_{\alpha} F$
satisfies the Universal Coefficient Theorem.
Therefore Huaxin Lin's classification theorem
(Theorem~5.2 of~\cite{Ln15}),
in the form presented in Proposition~3.7 of~\cite{PhtRp2},
together with the $K$-theory computation of Corollary~\ref{cor-K0},
implies the conclusion.
\end{proof}

The case $F = \Z_2$ is not new~\cite{P_BK}.

Using the results of Section~\ref{sec-gen},
one can determine exactly which AF~algebra the crossed product is.
As an example, we show that, for fixed $F,$
the crossed product algebras are isomorphic \ifo\  the
original irrational rotation algebras are isomorphic.

\begin{thm}\label{T:IsomOfCrPrd}
Let $F \subseteq \SL_2 (\ZZ)$
be one of the groups $\ZZ_2, \ZZ_3, \ZZ_4, \ZZ_6,$
and let $\af \colon F \to A_{\theta}$
be given by~(\ref{eq-action}) for the groups
generated as in~(\ref{eq-gen}).
Let $\te_1, \te_2 \in \R \smallsetminus \Q.$
Then
$A_{\te_1} \rtimes_{\alpha} F \cong A_{\te_2} \rtimes_{\alpha} F$
\ifo\  $\te_1 = \pm \te_2$ mod $\Z.$
\end{thm}

\begin{proof}
We give the proof for the case $F = \Z_2,$
using Theorem~\ref{thm-generators}(\ref{thm-generators-Z_2}).
The proofs for the other parts are the same,
using instead Parts~(\ref{thm-generators-Z_3}),
(\ref{thm-generators-Z_4}),
and~(\ref{thm-generators-Z_6}) of Theorem~\ref{thm-generators},
and choosing in place of $p^{\te}$ \pj s with traces
$\frac{1}{3},$ $\frac{1}{4},$ and $\frac{1}{6}.$

We clearly have
$A_{\te} \rtimes_{\alpha} \Z_2 \cong A_{\te+n} \rtimes_{\alpha} \Z_2$
for $n \in \Z,$
so \wolog\  $\te_1, \te_2 \in (0, 1).$
{}From Theorem~\ref{thm-generators}(\ref{thm-generators-Z_2}),
it is clear that for $\te \in (0, 1) \smallsetminus \Q$
the canonical \tst\  %
$\ta_{2, \te}$ on $A_{\te} \rtimes_{\alpha} \Z_2$ satisfies
\[
(\ta_{2, \te})_* (K_0 (A_{\te} \rtimes_{\alpha} \Z_2))
= \tfrac{1}{2} \Z + \tfrac{1}{2} \te \Z.
\]
It is well known
(and also follows directly from Proposition~\ref{P:TRPAndTraces}
and Corollary~\ref{C:TRP})
that $\ta_{2, \te}$ is the unique \tst\  %
on $A_{\te} \rtimes_{\alpha} \Z_2.$
Accordingly,
if $\te_1 \neq \pm \te_2$ mod $\Z,$
then
\[
2 (\ta_{2, \te_1})_* (K_0 (A_{\te_1} \rtimes_{\alpha} \Z_2))
   = \Z + \te_1 \Z
   \neq \Z + \te_2 \Z
   = 2 (\ta_{2, \te_2})_* (K_0 (A_{\te_2} \rtimes_{\alpha} \Z_2))
\]
whence
$A_{\te_1} \rtimes_{\alpha} \Z_2
     \not\cong A_{\te_2} \rtimes_{\alpha} \Z_2.$

If $\te_1, \te_2 \in (0, 1)$ with $\te_1 = \pm \te_2$ mod $\Z,$
then $\te_2 = \te_1$ or $\te_2 = 1 - \te_1.$
Thus, it remains to prove that
$A_{\te} \rtimes_{\alpha} \Z_2
    \cong A_{1-\te} \rtimes_{\alpha} \Z_2$
for $\te \in (0, 1) \smallsetminus \Q.$
By Theorem~\ref{T:CrPrdIsAF}
and George Elliott's classification theorem for AF~algebras
(Theorem~4.3 of~\cite{Ell-1}),
it suffices to exhibit an isomorphism
$f \colon K_0 ( A_{\te} \rtimes_{\alpha} \Z_2)
      \to K_0 ( A_{1-\te} \rtimes_{\alpha} \Z_2)$
of scaled ordered $K_0$-groups.
Since each algebra has a unique \tst,
and since the order on the $K_0$-group of a simple unital AF~algebra
is determined by its \tst\  %
(this is true much more generally:
see Corollary~5.7 and Theorems~5.8 and~6.8 of~\cite{LnTTR}),
it suffices to find a group isomorphism $f$ as above such that
$(\ta_{2, \, 1 - \te})_* \circ f = (\ta_{2, \te})_*$
and $f ([1]) = [1].$

Define $f$ on the basis elements listed in
Theorem~\ref{thm-generators}(\ref{thm-generators-Z_2})
by
$f \big( \big[ \E_{\te_1}^{\Z_2} \big] \big)
         = [p^{\te}] - \big[ \E_{\te_2}^{\Z_2} \big],$
and by sending each of the other basis elements
for $K_0 (A_{\te} \rtimes_{\alpha} \Z_2)$
to the corresponding basis element
for $K_0 (A_{1-\te} \rtimes_{\alpha} \Z_2),$
that is, $f ([1]) = [1],$ $f ([p^{\te}]) = [p^{1 - \te}],$ etc.
The computations of the trace on the generators given in
Theorem~\ref{thm-generators}(\ref{thm-generators-Z_2})
show that $(\ta_{2, \, 1 - \te})_* \circ f = (\ta_{2, \te})_*,$
as desired.
\end{proof}

In the cases $F = \Z_2$ and $F = \Z_4,$
it also is easy to use the isomorphism
$\ph \colon A_{\te} \to A_{1 - \te},$
given by $u_{\te} \mapsto v_{1 - \te}$
and $v_{\te} \mapsto u_{1 - \te},$
to directly exhibit an isomorphism
$A_{\te} \rtimes_{\alpha} F
    \cong A_{1-\te} \rtimes_{\alpha} F.$
Presumably this can be done in the other two cases as well.

Theorem~\ref{T:CrPrdIsAF} also has the following corollary.

\begin{cor}\label{C:FixIsAF}
Under the hypotheses of Theorem~\ref{T:CrPrdIsAF},
the fixed point algebra $A_{\te}^F$ is~AF.
\end{cor}

\begin{proof}
By the Proposition in~\cite{Rs},
the fixed point algebra $A_{\te}^F$ is isomorphic to
a corner of $A_{\te} \rtimes_{\alpha} F.$
\end{proof}

We now turn to the crossed product of a simple higher dimensional
noncommutative torus by the flip,
defined before Corollary~\ref{L:FlipTRP}.
The following theorem generalizes Theorem~3.1 of~\cite{Bc},
and completely answers
a question raised in the introduction to~\cite{P_FW2}.

\begin{thm}\label{T:HDFlipIsAF}
Let $\Te$ be a nondegenerate real $d \times d$ skew symmetric matrix.
Let $\ph \colon \Z_2 \to \Aut (A_{\Te})$ be the flip action.
Then $A_{\Te} \rtimes_{\ph} \Z_2$ is an AF~algebra.
\end{thm}

\begin{proof}
By Theorem~3.5 of~\cite{PhtRp2},
we know that $A_{\Te}$ is a simple separable unital \ca\  with
tracial rank zero.
By Corollary~\ref{L:FlipTRP},
we may apply Corollary~1.6 and Theorem~2.6 of~\cite{PhtRp1a}
to conclude that $A_{\Te}\rtimes_{\ph}\Z_2$
is a simple separable unital \ca\  with tracial rank zero.

Lemma~\ref{lem-twistext}
implies that the cocycle $\om_{\Theta}$ extends to  a cocycle
${\widetilde{\om}}_{\Theta}$ on the semidirect product
$\ZZ^d \rtimes \ZZ_2$
such that
$C^* (\ZZ^d \rtimes \ZZ_2, \, \widetilde{\om}_{\Theta})
          \cong A_{\Te}\rtimes_{\ph}\Z_2.$
Since $\ZZ^d \rtimes \ZZ_2$ can be embedded into $\RR^d \rtimes \ZZ_2,$
Corollary~\ref{C:UCT} therefore implies that $A_{\Te}\rtimes_{\ph}\Z_2$
satisfies the Universal Coefficient Theorem.

Since $\om_{\Theta}$ and $\widetilde{\om}_{\Theta}$ are real cocycles
(Definition~\ref{P_real}),
it follows from
Corollary~\ref{cor-real} that
\[
K_0 (A_{\Theta} \rtimes_{\ph} \ZZ_2)
= K_0 (C^* (\ZZ^d \rtimes \ZZ_2, \, \widetilde{\om}_{\Theta}))
\cong K_0 (C^* (\ZZ^d \rtimes \ZZ_2)),
\]
and
$K_1 (A_{\Theta} \rtimes_{\ph} \ZZ_2)
         \cong K_1 (C^* (\ZZ^d \rtimes \ZZ_2))$
for all $\Theta.$
But in~\cite{P_FW2} it was shown that
$K_0 (A_{\Theta} \rtimes_{\ph} \ZZ_2) \cong \ZZ^{3 \cdot 2^{d - 1}}$
and $K_1 (A_{\Theta} \rtimes_{\ph} \ZZ_2) = 0$
for some special values for $\Theta.$
This must then be true for {\emph{all}} $\Theta.$
Therefore Lin's classification theorem
(Theorem~5.2 of~\cite{Ln15}),
in the form presented in Proposition~3.7 of~\cite{PhtRp2},
implies the conclusion.
\end{proof}

\begin{cor}\label{C:FlipFixIsAF}
Let $\Te$ be a nondegenerate real $d \times d$ skew symmetric matrix.
Let $\ph \colon \Z_2 \to \Aut (A_{\Te})$ be the flip action.
Then the fixed point algebra $A_{\Theta}^{\ZZ_2}$ is~AF.
\end{cor}

\begin{proof}
The proof is the same as for Corollary~\ref{C:FixIsAF}.
\end{proof}

\end{document}